\documentclass[a4paper, 12pt, abstraction]{scrartcl}

\usepackage{amsmath}
\usepackage{amssymb}
\usepackage{graphicx}
\usepackage{epstopdf}
\usepackage{inputenc}
\usepackage{geometry}

\usepackage{caption}
\usepackage{subcaption}
\usepackage{cite}
\usepackage{calc}
\usepackage{float}
\usepackage{grffile} 
\usepackage{enumitem}

\usepackage{color}
\usepackage{xcolor}

\usepackage{mathtools}
\DeclarePairedDelimiter{\ceil}{\lceil}{\rceil} 

\usepackage{comment} 

\newcommand{\U}{\mathcal{U}}

\newcommand{\R}{\mathbb{R}}
\newcommand{\X}{\mathcal{X}}

\usepackage{algpseudocode}
\usepackage{algorithm}

\newcommand{\ddre}{d.d.\,Re}

\usepackage{multirow}
\usepackage{rotating}

\usepackage{authblk}

\usepackage{lscape} 

\usepackage[normalem]{ulem} 

\begin{document}

\title{A Largest Empty Hypersphere Metaheuristic for Robust Optimisation with Implementation Uncertainty\thanks{Funded through EPSRC grants EP/L504804/1 and EP/M506369/1.}}
\author[1]{Martin Hughes\footnote{Corresponding author. Email: \texttt{m.hughes6@lancaster.ac.uk}} }
\author[2]{Marc Goerigk}
\author[1]{Michael Wright}
\affil[1]{Department of Management Science, Lancaster University, United Kingdom}
\affil[2]{University of Siegen, Germany}

\date{} 

\maketitle


\abstract{We consider box-constrained robust optimisation problems with implementation uncertainty. In this setting, the solution that a decision maker wants to implement may become perturbed. The aim is to find a solution that optimises the worst possible performance over all possible perturbances.

Previously, only few generic search methods have been developed for this setting. We introduce a new approach for a global search, based on placing a largest empty hypersphere. We do not assume any knowledge on the structure of the original objective function, making this approach also viable for simulation-optimisation settings. In computational experiments we demonstrate a strong performance of our approach in comparison with state-of-the-art methods, which makes it possible to solve even high-dimensional problems.
}


\section{Introduction}
\label{sec:Introduction}

The use of models to support informed decision making is ubiquitous. However, the size and nature of the decision variable solution space, and the model runtime, may make a comprehensive -- exhaustive or simply extensive -- evaluation of the problem space computationally infeasible. In such cases an efficient approach is required to search for global optima.

Mathematical programs are one form of model that are explicitly formulated as optimisation problems, where the model representation imposes assumptions on the structure of the decision variable space and objective function. Such models are well suited to efficient solution, and identification of global optima may be theoretically guaranteed when feasible solutions exist. However many real-world problems are not suited to expression as a mathematical program (e.g., a solution is evaluated by using a simulation tool). From an optimisation perspective models where no assumptions are made about the model structure can be thought of as a black-box, where decision variables values are input and outputs generated for interpretation as an objective. In this case optimisation search techniques such as metaheuristics are required, i.e., general rule-based search techniques that can be applied to any model.

An additional widespread feature of many real-world problems is the consideration of uncertainty which may impact on model outputs, and so on corresponding objective function values. One strategy is to simply ignore any uncertainty and perform a standard search, possibly assessing and reporting on the sensitivity of the optimum after it has been identified. However it has been established that optimal solutions which are sensitive to parameter variations within known bounds of uncertainty may substantially degrade the optimum objective function value, meaning that solutions sought without explicitly taking account of uncertainty are susceptible to significant sub-optimality, see \cite{BenTalElGhaouiNemirovski2009, GoerigkSchobel2016}. In the face of uncertainty the focus of attention for an optimisation analysis shifts from the identification of a solution that just performs well in the expected case, to a solution that performs well over a range of scenarios.

In this paper we develop a new algorithm for box-constrained robust black-box global optimisation problems taking account of implementation uncertainty, i.e., the solution that a decision maker wants to implement may be slightly perturbed in practice, and the aim is to find a solution that performs best under the worst-case perturbation. Our method is based on an exploration technique that uses largest empty hyperspheres (LEHs) to identify regions that can still contain improving robust solutions. In a computational study we compare our method with a local search approach from the literature (see \cite{BertsimasNohadaniTeo2010}) and a standard particle swarm approach. We find that our approach considerably outperforms these methods, especially for higher-dimensional problems.

\paragraph{Structure of this paper.} We begin with a review of the literature on metaheuristics for robust optimisation in Section~\ref{sec:literature} before outlining the formal description of robust min max problems in Section~\ref{sec:RobustOptimisation}. We also consider some of the details of the established local robust search technique due to \cite{BertsimasNohadaniTeo2010}. In Section~\ref{sec:LargestEmptyHypersphere} we introduce a novel approach, an exploration-focused movement through the search space identifying areas that are free of previously identified poor points. We include a discussion and descriptions of the algorithms used to identify empty regions of the decision variable search space. The approach is then tested against alternative heuristics in Section~\ref{sec:ExperimentsResults}, on test problems of varying dimension. The experimental set up is described and the results of this analysis presented. Finally we summarise and consider further extensions of this work in Section~\ref{sec:SummaryConcusionsFurtherWork}.

\section{Literature review}
\label{sec:literature}

\subsection{Robust optimisation}

Different approaches to model uncertainty in decision making problems have been explored in the literature. Within robust optimisation, a frequent distinction is made between parameter uncertainty (where the problem data is not known exactly) and implementation uncertainty (where a decision cannot be put into practice with full accuracy). Implementation uncertainty is also known as decision uncertainty \cite{BenTalElGhaouiNemirovski2009, Ghazali2009, BertsimasNohadaniTeo2010}. 

A common approach to the incorporation of uncertainty for black-box problems is stochastic optimisation. Here knowledge of the probability distributions of the uncertain parameters is assumed and some statistical measure of the fitness of a solution assessed, e.g. using Monte Carlo simulation to estimate the expected fitness. This may be the expected value, or a more elaborate model such as the variance in the fitness of a solution, or even a multi-objective optimisation setting, see \cite{PaenkeBrankeJin2006, HomemdeMelloBayraksan2014}.

An alternative to a stochastic approach is robust optimisation, whose modern form was first developed in \cite{KouvelisYu1997} and \cite{BenTalNemirovski1998}. Whereas with stochastic optimisation a knowledge of probability distributions over all possible scenarios is typically assumed, in robust optimisation it is only assumed that some set is identified containing all possible uncertainty scenarios (potentially infinite in number). A classic robust approach is then to find a solution across all scenarios that is always feasible (strictly robust) and optimises its performance in the worst case. This is known as min max. For a given point in the decision variable space there is an `inner' objective to identify the maximal (worst case) function value in the local uncertainty neighbourhood, and an overall `outer' objective to identify the minimum such maximal value.

The field of robust optimisation has been primarily aligned with mathematical programming approaches. There the methodology is based around the definition of reasonable uncertainty sets and the reformulation of computationally tractable mathematical programming problems. For specific forms of convex optimisation problems, the problem incorporating uncertainty can be re-formulated to another tractable, convex problem, see \cite{BertsimasNohadaniTeo2007, GohSim2010}. To overcome concerns that the strictly robust worst case approach may be overly conservative, the concept of robustness can be expanded in terms of both the uncertainty set considered and the robustness measure \cite{GoerigkSchobel2016}. On the assumption that it is overly pessimistic to assume that all implementation errors take their worst value simultaneously \cite{BertsimasSim2004} consider an approach where the uncertainty set is reduced, and a robust model defined where the optimal solution is required to remain feasible for uncertainty applied to only a subset of the decision variables at any given time. Min max regret, see \cite{AissiBazganVanderpooten2009}, is an alternative to min max, seeking to minimise the maximum deviation between the value of the solution and the optimal value of a scenario, over all scenarios. \cite{BenTalBertsimasBrown2010} considers soft robustness, which utilises a nested family of uncertainty sets. The distributionally robust optimisation approach, see \cite{GohSim2010}, attempts to bridge robust and stochastic techniques by utilizing uncertainty defined as a family of probability distributions, seeking optimal solutions for the worst case probability distribution. \cite{ChasseinGoerigk2016} use a bi-objective approach to balance average and worst case performance by simultaneously optimising both. 

Robust optimisation in a mathematical programming context has been application-driven, so considerable work has been undertaken in applying robustness techniques to specific problems or formulations, see \cite{BeyerSendhoff2007, GoerigkSchobel2016}. There has also been some cross-over into the application of specific heuristics, for example see \cite{GoldenLaporteTaillard1997, ValleMartinezdaCunhaMateus2011}. However application to general problems has been less well addressed \cite{GoerigkSchobel2016}. Furthermore robust approaches applied to black-box models are much less widely considered than approaches for mathematical programming problems, see \cite{MarzatWalterPietLahanier2013, GoerigkSchobel2016, MarzatWalterPietLahanier2016}. Recently, robust optimisation with implementation uncertainty has also been extended to multi-objective optimisation, see \cite{eichfelder2017decision}.

\subsection{Metaheuristic for robust optimisation}

The min max approach has been tackled with standard metaheuristic techniques applied to both the inner maximisation and outer minimisation problems. In co-evolutionary approaches two populations (or swarms) evolve separately but are linked. The fitness of individuals in one group is informed by the performance of individuals in the other, see \cite{CramerSudhoffZivi2009}. \cite{Herrmann1999, Jensen2004} use such a two-population genetic algorithm (GA) approach, whilst \cite{ShiKrohling2002, MasudaKuriharaAiyoshi2011} consider two-swarm co-evolutionary particle swarm optimisation (PSO) techniques for min max problems. A brute force co-evolutionary approach is to employ complete inner maximisation searches to generate robust values for each individual in each generation of the outer minimisation, however this is expensive in terms of model runs (i.e.,  function evaluations), see \cite{MarzatWalterPietLahanier2016}. More practical co-evolutionary approaches, for example using only small numbers of populations for the outer search and the inner (uncertainty) search which share information between populations from generation to generation, or following several generations, require the application of additional simplifications and assumptions, see \cite{CramerSudhoffZivi2009, MasudaKuriharaAiyoshi2011}.

One general area of research is the use of emulation to reduce the potential burden of computational run times and the number of model-function evaluations, see \cite{VuDAmbrosioHamadiLiberti2016}. \cite{ZhouZhang2010} use a surrogate-assisted evolutionary algorithm to tackle the inner search for black-box min max problems. \cite{MarzatWalterPietLahanier2013, MarzatWalterPietLahanier2016} employs Kriging meta-modelling coupled with an expected improvement (EI) metric, as well as a relaxation of the inner maximisation search. The EI metric is used to efficiently choose points in the decision variable space where nominal (expensive) function evaluation should be undertaken, see \cite{JonesSchonlauWelch1998}, here with a view to most efficiently improving the estimate of the robust global minimum. The relaxation involves iteratively performing the outer minimisation on a limited inner uncertainty neighbourhood followed by an extensive inner maximisation search in the region of the identified outer minimum. This continues whilst the inner search sufficiently deteriorates the outer solution, with the inner maximum point being added to the limited inner uncertainty set with each iteration.

A second approach due to \cite{urRehmanLangelaarvanKeulen2014, urRehmanLangelaar2017} also uses Kriging and an EI metric, building on a meta-model of the expensive nominal problem by applying a robust analysis directly to the Kriging model and exploiting the fact that many more inexpensive function evaluations can be performed on this meta-model. A modified EI metric is calculated for the worst case cost function of the meta-model, to efficiently guide the search in the nominal expensive function space. In \cite{urRehmanLangelaar2017} the approach is applied to a constrained non-convex 2 dimensional problem due \cite{BertsimasNohadaniTeo2010, BertsimasNohadaniTeo2010nonconvex}, the unconstrained version of which is also considered here. The Kriging-based approach is shown to significantly outperform the approaches outlined here, in terms of the number of expensive function evaluations required to converge towards the robust optimum. In general we would expect the approach from \cite{urRehmanLangelaarvanKeulen2014, urRehmanLangelaar2017} to outperform the approaches considered here, in terms of efficiency when applied to  low dimensional non-convex problems. However the primary challenge with meta-model based approaches is their application to higher dimensional problems. The test cases considered in \cite{MarzatWalterPietLahanier2013, MarzatWalterPietLahanier2016, urRehmanLangelaarvanKeulen2014, urRehmanLangelaar2017} have either been restricted to low dimensional non-convex problems, or simpler convex and convex-concave problems of up to 10 dimensions.

One local black-box min max approach is due to \cite{BertsimasNohadaniTeo2007, BertsimasNohadaniTeo2010, BertsimasNohadaniTeo2010nonconvex}. Here a search is undertaken by iteratively moving along 'descent directions'. Uncertainty around individual points is assessed using local gradient ascents, based on which undesirable 'high cost points' (hcps) are identified. Steps are taken in directions which point away from these hcps, until no direction can be found.

Our approach is inspired by both elements of the descent directions technique and the concept of relaxation of the inner maximisation search. We extend the idea of locally moving away from identified hcps to a global perspective, seeking regions of the solution space currently empty of such undesirable points. Furthermore the nature of our outer approach enables the curtailing of an inner maximisation search if it is determined that the current point under consideration cannot improve on the current best robust global solution.


\section{Notation and previous results}
\label{sec:RobustOptimisation}

\subsection{Problem description}
\label{sec:ProblemDescription}

We consider a general optimisation problem of the form
\begin{align*}
	\quad \min\ & f(\pmb{x}) \\
	\text{s.t. } & \pmb{x} \in \X
\end{align*}
where $\pmb{x}=(x_{1}, x_{2}, \ldots, x_{n})^T$ denotes the $n$-dimensional vector of decision variables, $f: \R^n \to \R$ is the objective function, and $\X \subseteq \R^n$ is the set of feasible solutions. We write $[n]:=\{1,\ldots,n\}$. In this paper, we assume box constraints $\X = \prod_{i\in[n]} [l_i,u_i]$. Any other potential feasibility constraints are assumed to be ensured through a penalty in the objective. 

In implementation uncertainty, we assume that a desired solution $\pmb{x}$ might not be possible to put into practice with full accuracy. Instead, a ''close'' solution $\tilde{\pmb{x}}=\pmb{x}+\Delta\pmb{x}$ may be realised. The aim is to find a robust $\pmb{x}$ such that for any such solution $\tilde{\pmb{x}}$ from the neighbourhood of $\pmb{x}$, the worst case performance is optimised.

More formally, we follow the setting of \cite{BertsimasNohadaniTeo2010} and consider the so-called uncertainty set
\[
	\U:=\{ \Delta \pmb{x}\in\R^n \mid \| \Delta \pmb{x} \| \leq \Gamma \}
\]
where $\Gamma > 0$ defines the magnitude of the uncertainty, and $\|\cdot\|$ refers to the Euclidean norm. The worst case costs of a solution $\pmb{x}\in\X$ are then given as 
\[ g(\pmb{x}):=\max_{\Delta \pmb{x} \in \U} f(\pmb{x} + \Delta \pmb{x}) \]
and so the robust optimisation problem is given by:
\[ \min_{\pmb{x}\in\X} g(\pmb{x}) = \min_{\pmb{x}\in\X} \max_{\Delta\pmb{x} \in \U} f(\pmb{x} + \Delta \pmb{x}) \tag{ROP}\]
We therefore have an inner maximisation and outer minimisation problem, such that the identification of the robust global optimum is based on finding the (outer) minimum worst case cost objective function value in the decision variable space, and that objective is determined by the (inner) maximisation of the nominal objective function in the uncertainty neighbourhood around each point in the decision variable space. This type of problem is also known as min max.

Note that $\pmb{x}+\Delta\pmb{x}$ may not be in $\X$, for which reason we assume that the definition of $f$ is not limited to $\X$. However, if it is desired that $\pmb{x}+\Delta\pmb{x}\in\X$ for all $\Delta\pmb{x}\in\U$, then this can be ensured by reducing the size of the feasible search space by $\Gamma$.

\begin{figure}[htbp]
	\centering
	
	\begin{subfigure}{.5\textwidth}
		\centering
		\includegraphics[width=2.8in, height=3.0in]{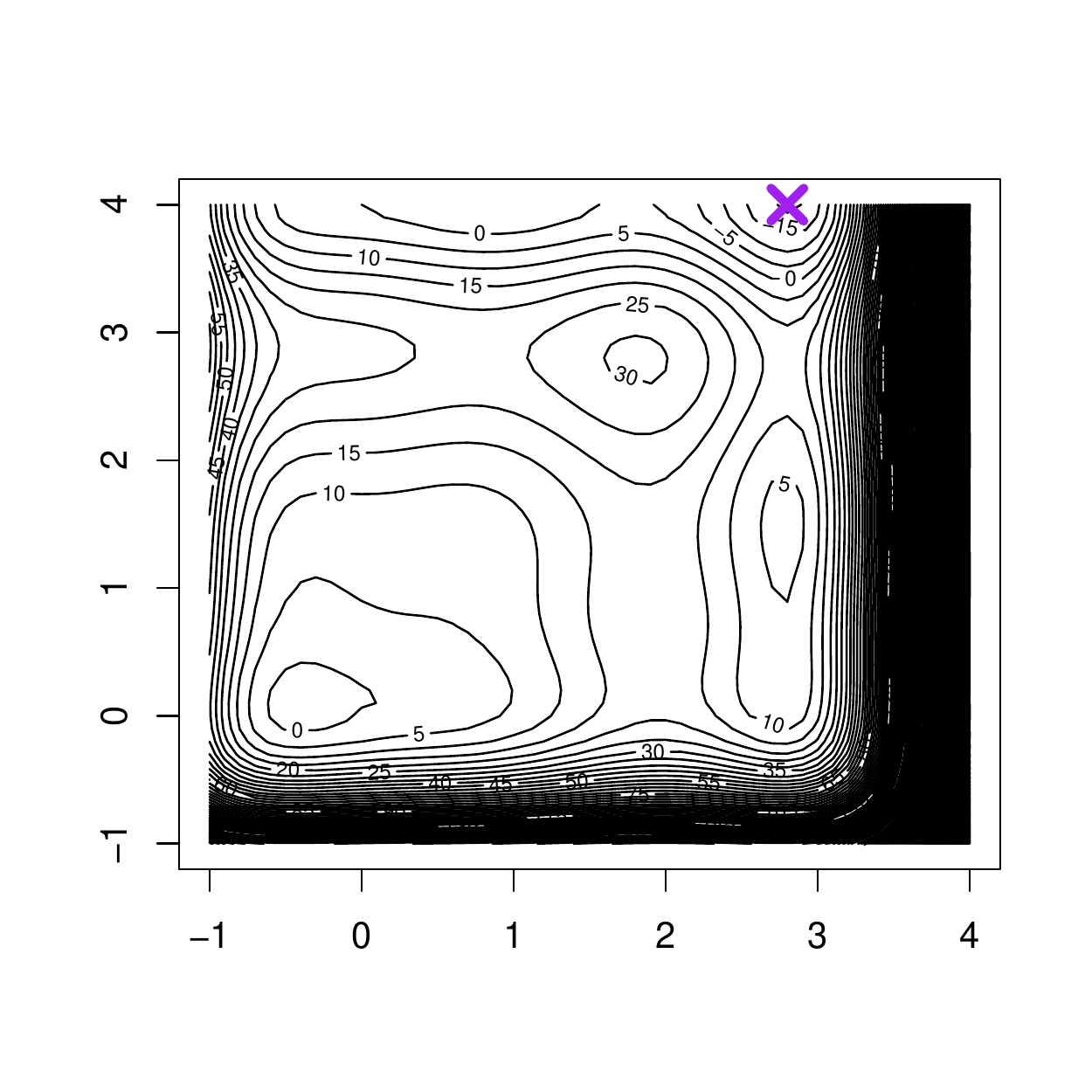} 	
    \vspace*{-7mm} 
		\caption{Nominal problem}
		\label{BertsimasNominal}
	\end{subfigure}%
	\begin{subfigure}{.5\textwidth}
		\centering
		\includegraphics[width=2.8in, height=3.0in]{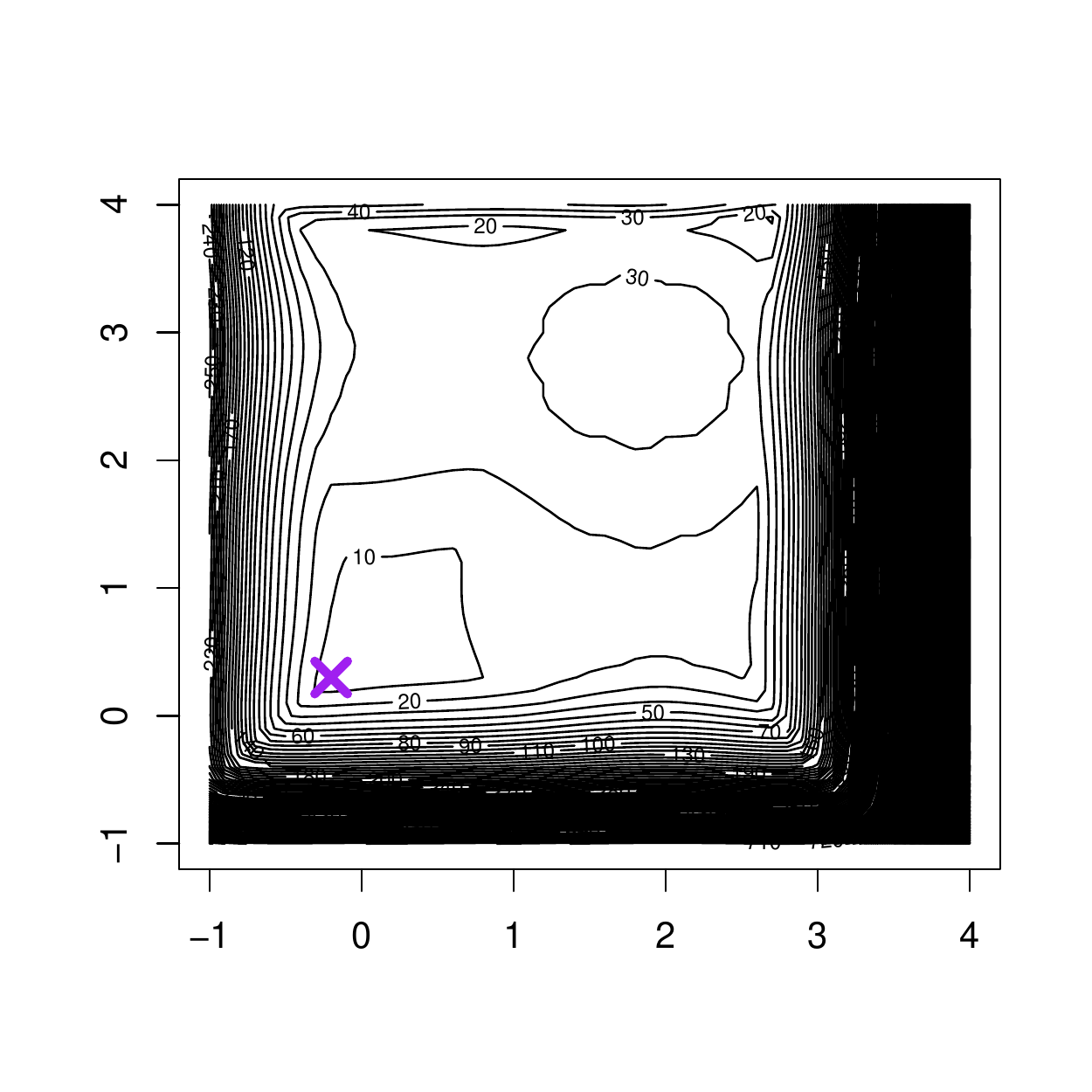}
    \vspace*{-7mm} 
    \caption{Worst case problem with $\Gamma$=0.5}
		\label{BertsimasWorstCase}
	\end{subfigure}

	\caption{Nominal and worst case cost for (poly2D) from \cite{BertsimasNohadaniTeo2010}. Marked in purple are the respective optima.}
	\label{fig:NominalWorstCaseBertsimasPoly}
\end{figure}

As an example for our problem setting, consider the 2-dimensional polynomial function due to \cite{BertsimasNohadaniTeo2010}:
\begin{align*}
f(x, y) = &2x^6 - 12.2x^5 + 21.2x^4 + 6.2x - 6.4x^3 - 4.7x^2 - y^6 \\
					- &11y^5 + 43.3y^4 - 10y - 74.8y^3 + 56.9y^2 - 4.1xy \\
					 - &0.1y^2x^2 + 0.4y^2x + 0.4x^2y \tag{poly2D}
\end{align*}
For a feasible solution space within bounds $[-1, 4]$ in each dimension, and uncertainty defined by a $\Gamma$-radius value of $0.5$, the nominal and worst case plots for (poly2D) are shown in Figure~\ref{fig:NominalWorstCaseBertsimasPoly}. In min max the problem is one of finding the global minimum for the worst case cost function. If uncertainty is ignored the problem is just one of finding the global minimum of the (nominal) objective as shown in Figure~\ref{BertsimasNominal}, whereas including uncertainty the problem becomes one of finding the (worst case cost) objective as shown in Figure~\ref{BertsimasWorstCase}. In both cases the search proceeds based on generating nominal objective values but for the worst case cost we must further undertake some assessment of the impact of uncertainty on those objective outputs. 

\vspace*{-6mm} 

Here the global optimum value for the nominal problem is -20.8 at (2.8, 4.0). The worst case plot is estimated by randomly sampling large numbers of points in the $\Gamma$-uncertainty neighbourhood around each plotted point. The worst case cost at each point is then approximated as the highest value of $f(x)$ identified within each $\Gamma$-uncertainty neighbourhood. The global optimum for the worst case problem is approximately 4.3 at (-0.18, 0.29). The significant shift in the nominal versus robust optima, both in terms of its location and the optimum objective, emphasises the potential impact of considering uncertainty in decision variable values. The difference between the nominal and robust optimal objective function values is the `price of robustness', see \cite{BertsimasSim2004}.

\subsection{Local robust search using descent directions}
\label{sec:LocalRobustSearchDescentDirections}

\begin{figure}[htbp]
	\centering	
	\begin{subfigure}{.44\textwidth}
		\centering
		\includegraphics[width=.7\textwidth]{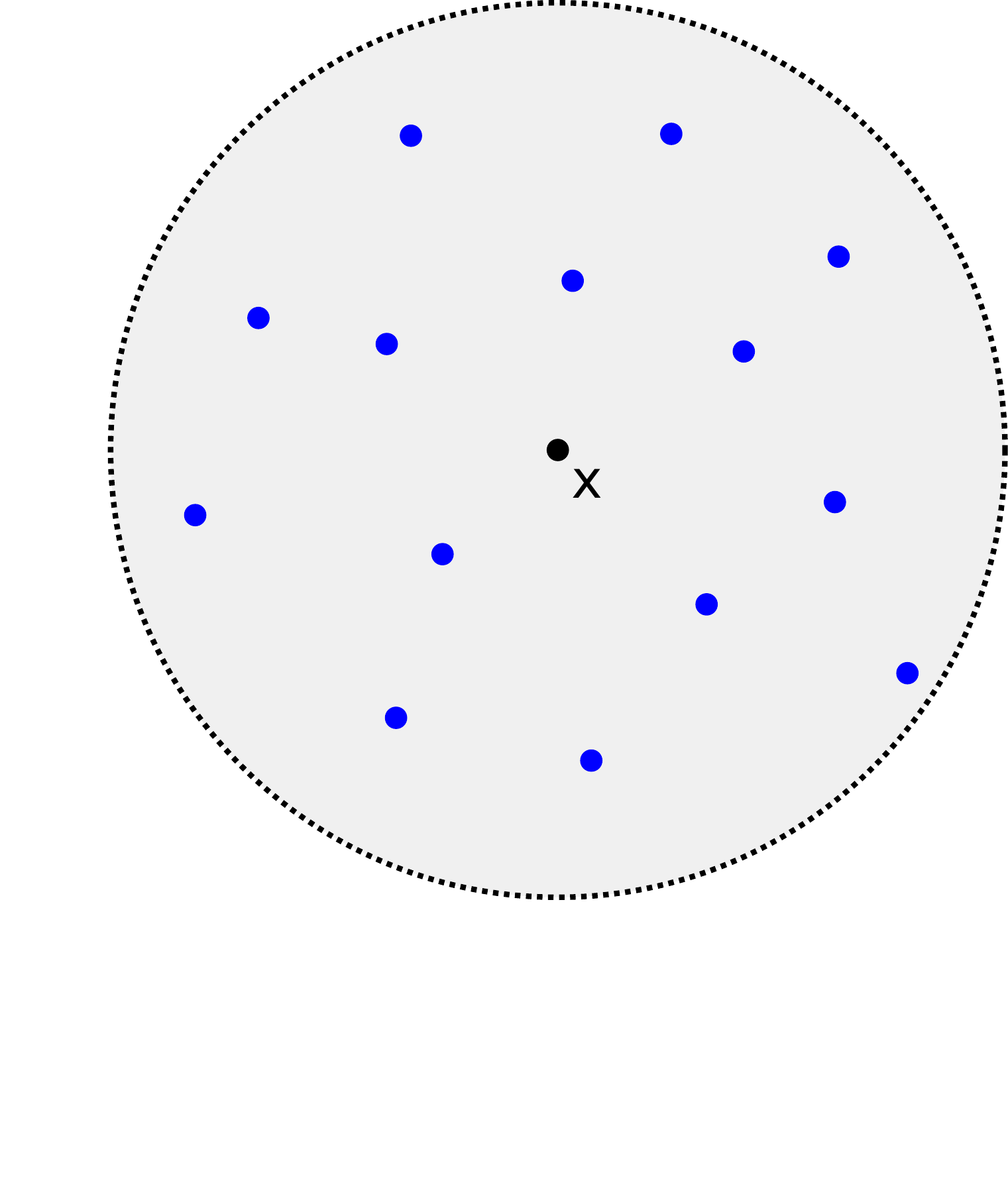} 
  	\caption{Candidate point $\pmb{x}$ (centre), and points evaluated for the inner maximisation problem (blue).}
		\label{fig:algo2}
    \vspace*{9mm} 
	\end{subfigure}%
	\hspace{.1\textwidth}
	\begin{subfigure}{.44\textwidth}
		\centering
		\includegraphics[width=.7\textwidth]{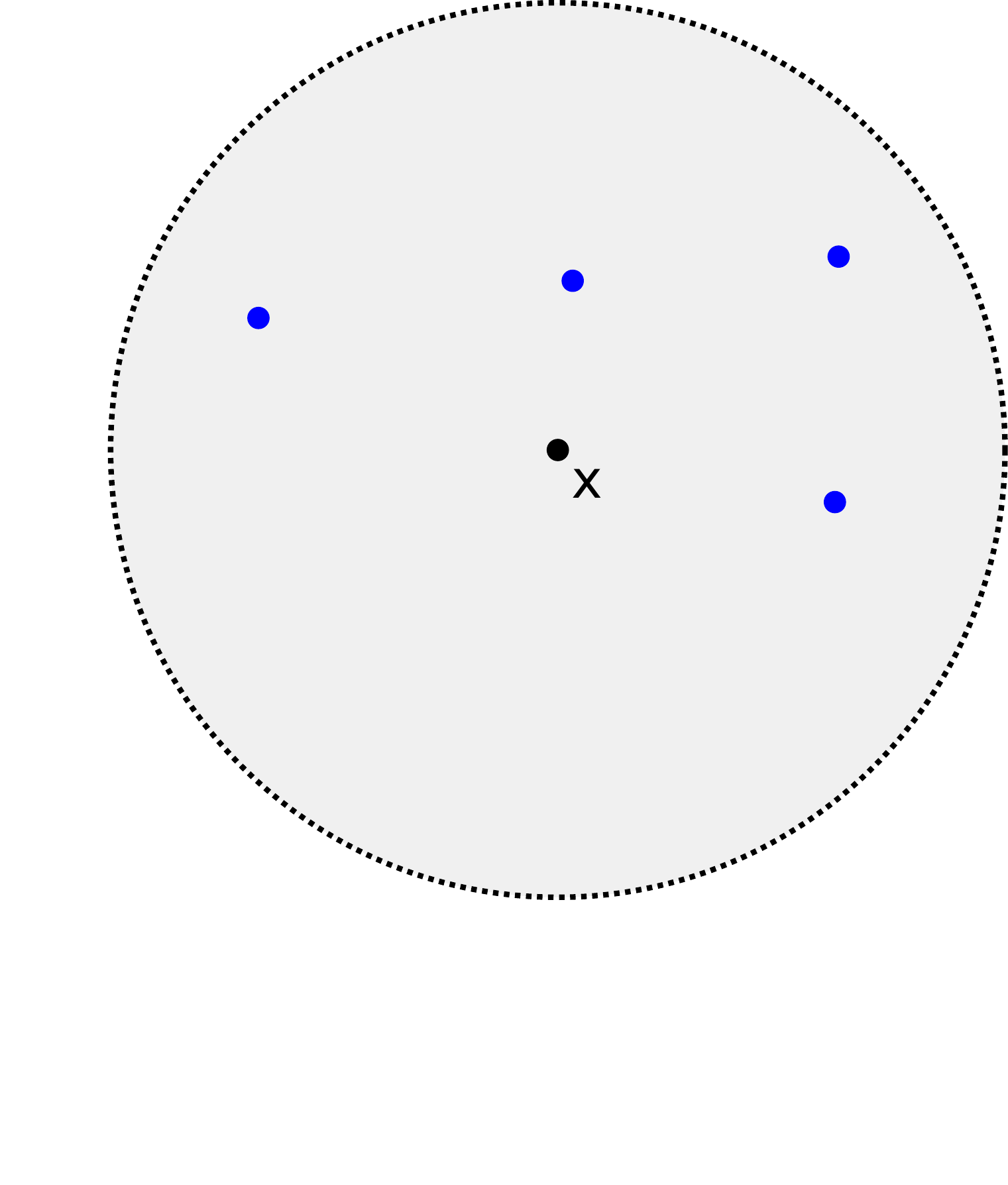}
  	\caption{Subset $H(\pmb{x})$ of critical high cost points.}
		\label{fig:algo3}
    \vspace*{9mm} 
	\end{subfigure} 
	
	\begin{subfigure}{.44\textwidth}
		\centering
		\includegraphics[width=.7\textwidth]{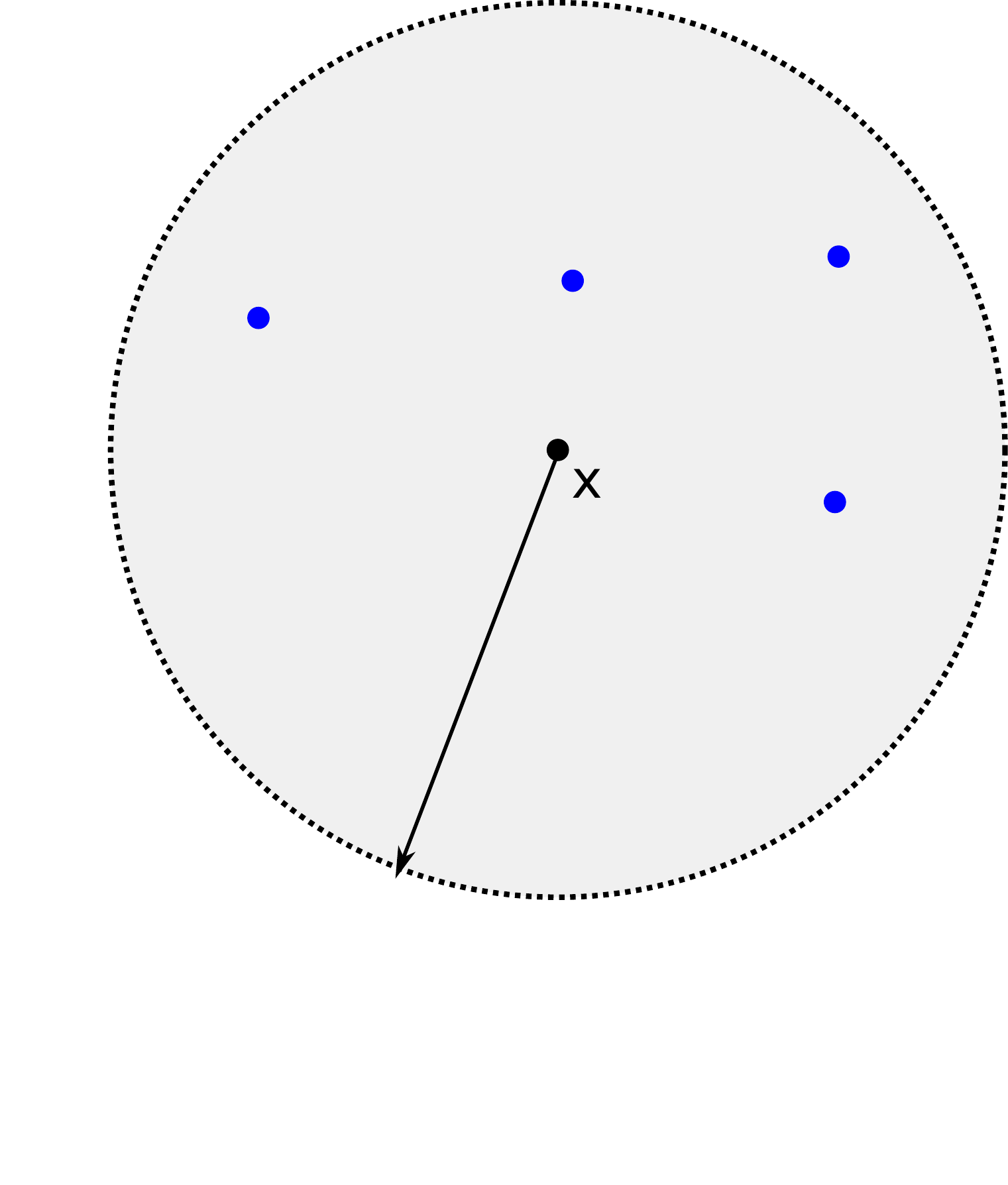}
		\caption{A descent direction is identified by solving a second order cone problem.}
		\label{fig:algo5}
	\end{subfigure}%
	\hspace{.1\textwidth}
	\begin{subfigure}{.44\textwidth}
		\centering
		\includegraphics[width=.7\textwidth]{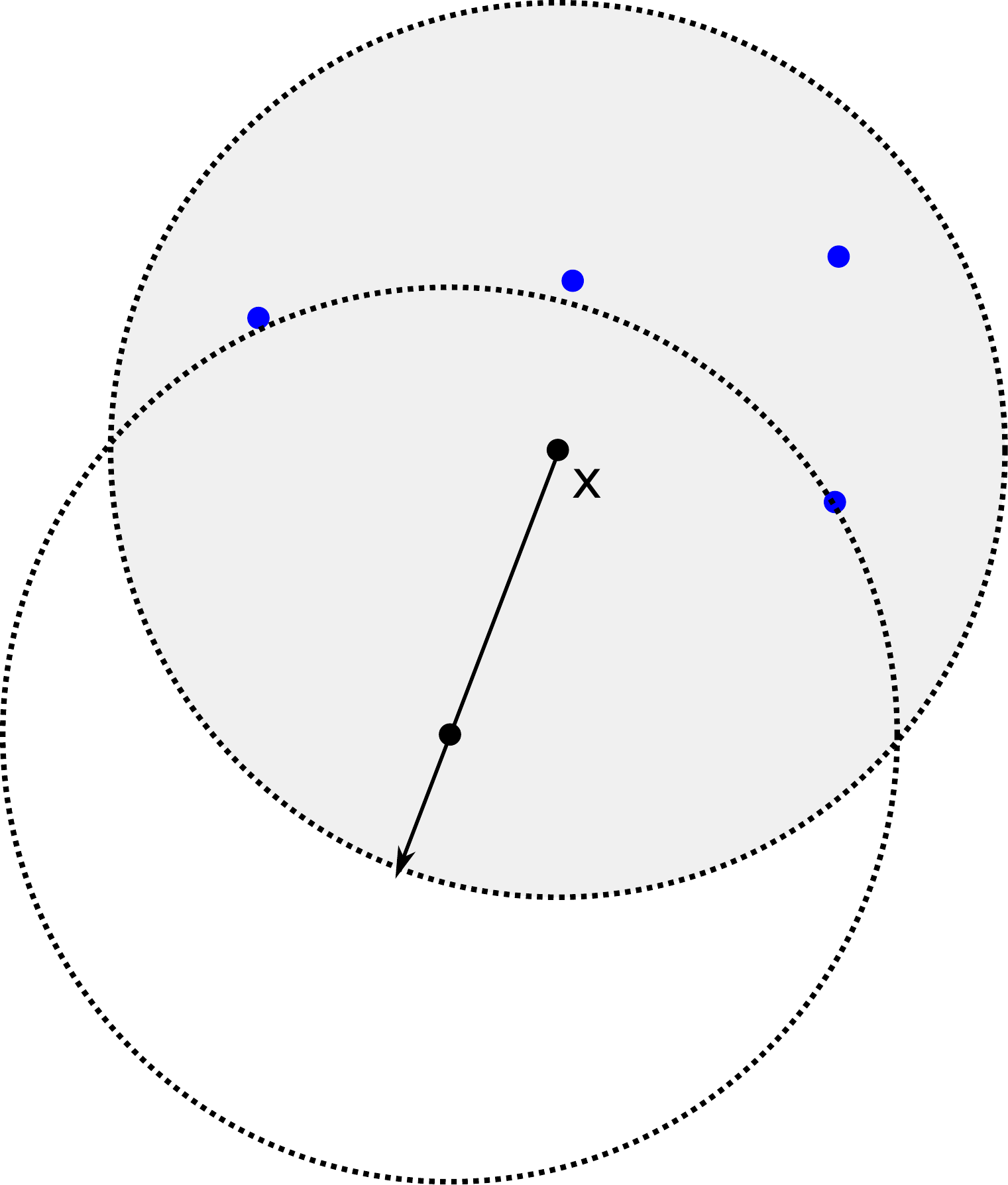}
  	\caption{The step size is determined.}
		\label{fig:algo6}
	\end{subfigure}	
	
	\caption{Description of the descent direction robust local search approach \cite{BertsimasNohadaniTeo2010}.}
	\label{fig:BertsimasDescription}
\end{figure}

We briefly summarise the local search approach for (ROP) that was developed in \cite{BertsimasNohadaniTeo2010}.
Here, (ROP) is solved using a local robust optimisation heuristic illustrated by Figure~\ref{fig:BertsimasDescription}. An initial decision variable vector $\hat{\pmb{x}}$ is randomly sampled. Then a series of gradient ascent searches are undertaken within the $\Gamma$-uncertainty neighbourhood of this candidate solution to identify hcps, see Figure~\ref{fig:algo2}. This approximates the inner maximisation problem $\max_{\Delta\pmb{x}} f(\hat{\pmb{x}} + \Delta\pmb{x})$. Using a threshold value that is dynamically adjusted during the algorithm, a subset $H(\hat{\pmb{x}})$ of all evaluated points is identified, see Figure~\ref{fig:algo3}. 

In the next step, a descent direction is identified that points away from the set $H(\hat{\pmb{x}})$, see Figure~\ref{fig:algo5}. To this end, a mathematical programming approach is used, minimising the angle between the hcps and the candidate solution. This leads to the following second order cone problem.
\begin{align}
\min_{\pmb{d},\beta} \ & \beta \label{soc1} \\
\text{s.t. } &\| \pmb{d}\| \le 1 \label{soc2} \\
& \pmb{d}^T \pmb{h} \le \beta & \forall \pmb{h}\in H(\hat{\pmb{x}}) \label{soc3}\\
& \beta \le -\varepsilon \label{soc4}
\end{align}
Here, $\pmb{d}$ is the descent direction, which is normalised by Constraint~\eqref{soc2}. Constraints~\eqref{soc3} ensure that $\beta$ is the maximum angle between $\pmb{d}$ and all high cost points $\pmb{h}$. Through Constraint~\eqref{soc4}, we require a feasible descent direction to point away from all points in $H(\hat{\pmb{x}})$. When an optimal direction cannot be found, the algorithm stops -- a robust minimum has been reached.

Next the size of the step to be taken is calculated, see Figure~\ref{fig:algo6}. A step size just large enough to ensure that all of the hcps are outside of the $\Gamma$-uncertainty neighbourhood of the next candidate solution is used. Using the identified descent direction and step size the algorithm moves to a new candidate point, and so the heuristic repeats iteratively until a robust minimum has been identified.

\section{A new largest empty hypersphere approach}
\label{sec:LargestEmptyHypersphere}

\subsection{Algorithm overview}
\label{sec:GlobalRobustSearchLEH}

\begin{figure}[H]
	\centering
	\begin{subfigure}{.48\textwidth}
		\centering
		\includegraphics[width=2.2in]{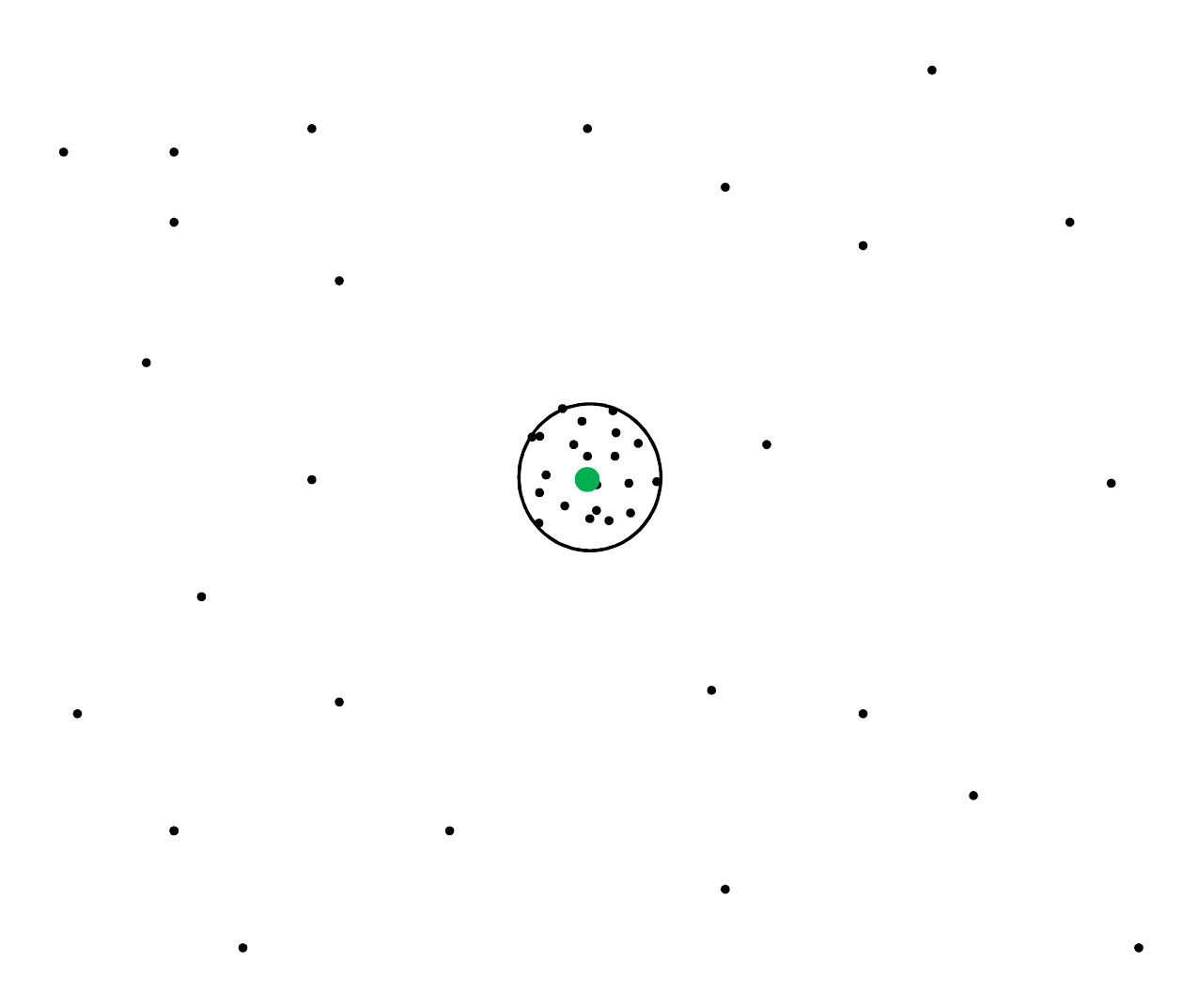} 	
  	\caption{The decision variable space is seeded randomly. Perform an inner search around one candidate point.}
		\label{fig:newalgo3}
    \vspace*{5mm} 
	\end{subfigure}%
			\hfill
	\begin{subfigure}{.48\textwidth}
		\centering
		\includegraphics[width=2.2in]{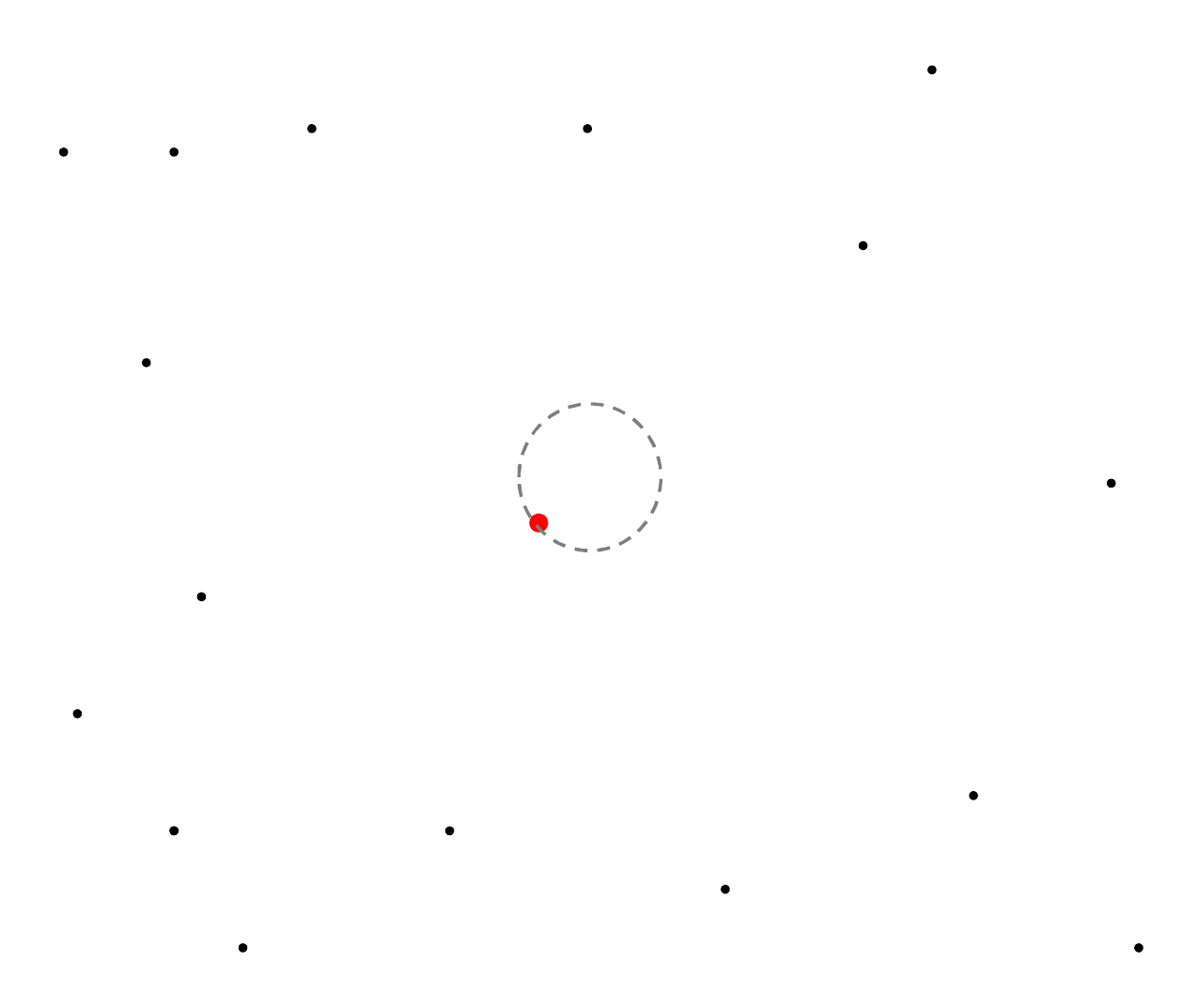} 
  	\caption{The current high cost set, including one point from the previous inner search and some of the seed points.}
		\label{fig:newalgo4}
    \vspace*{5mm} 
	\end{subfigure}
	
	\begin{subfigure}{.48\textwidth}
		\centering
		\includegraphics[width=2.2in]{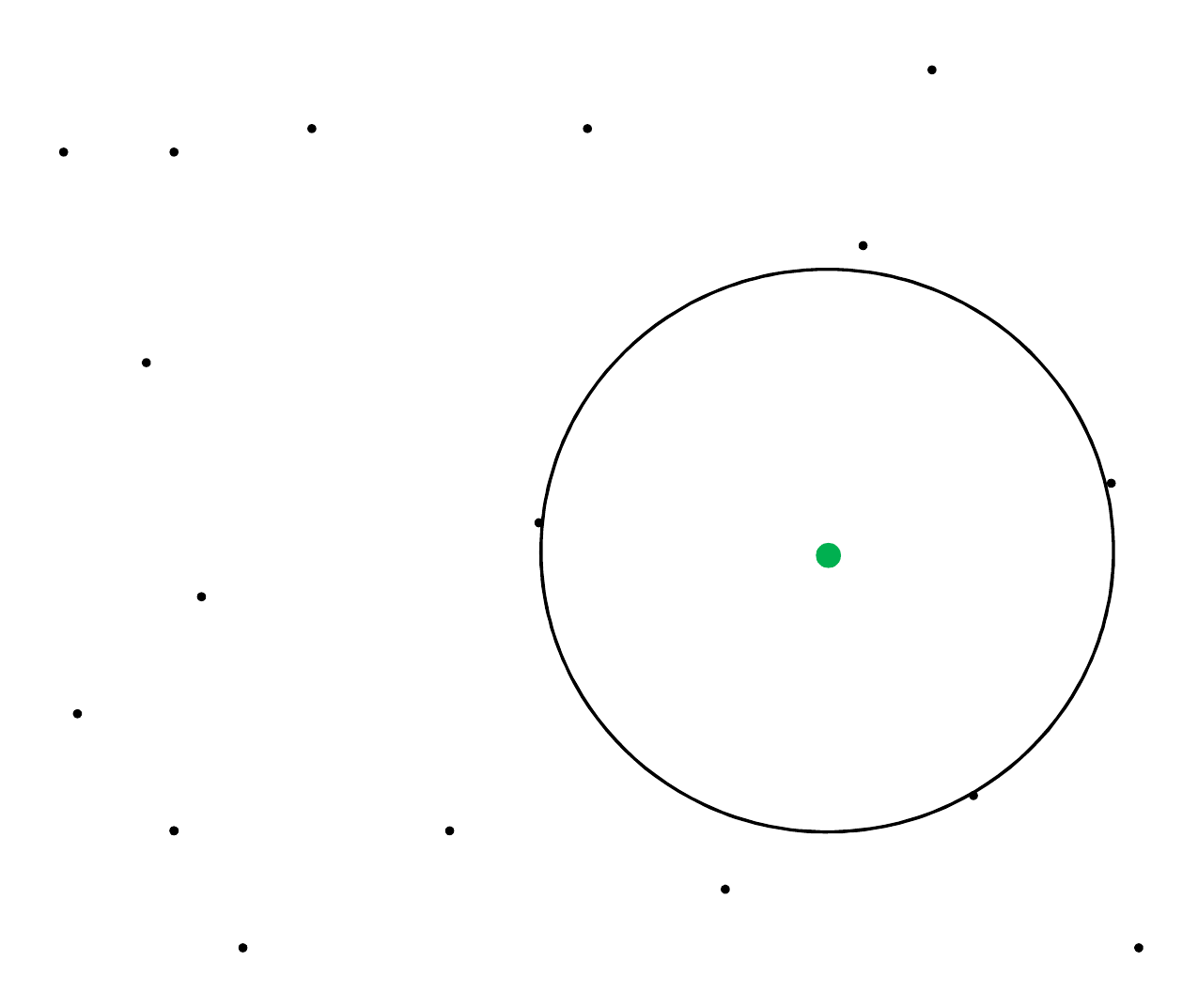} 	
		\caption{Identify the largest empty hypersphere, the centre of which is the next candidate point.}
		\label{fig:newalgo5}
    \vspace*{5mm} 
	\end{subfigure}%
	\hfill
	\begin{subfigure}{.48\textwidth}
		\centering
		\includegraphics[width=2.2in]{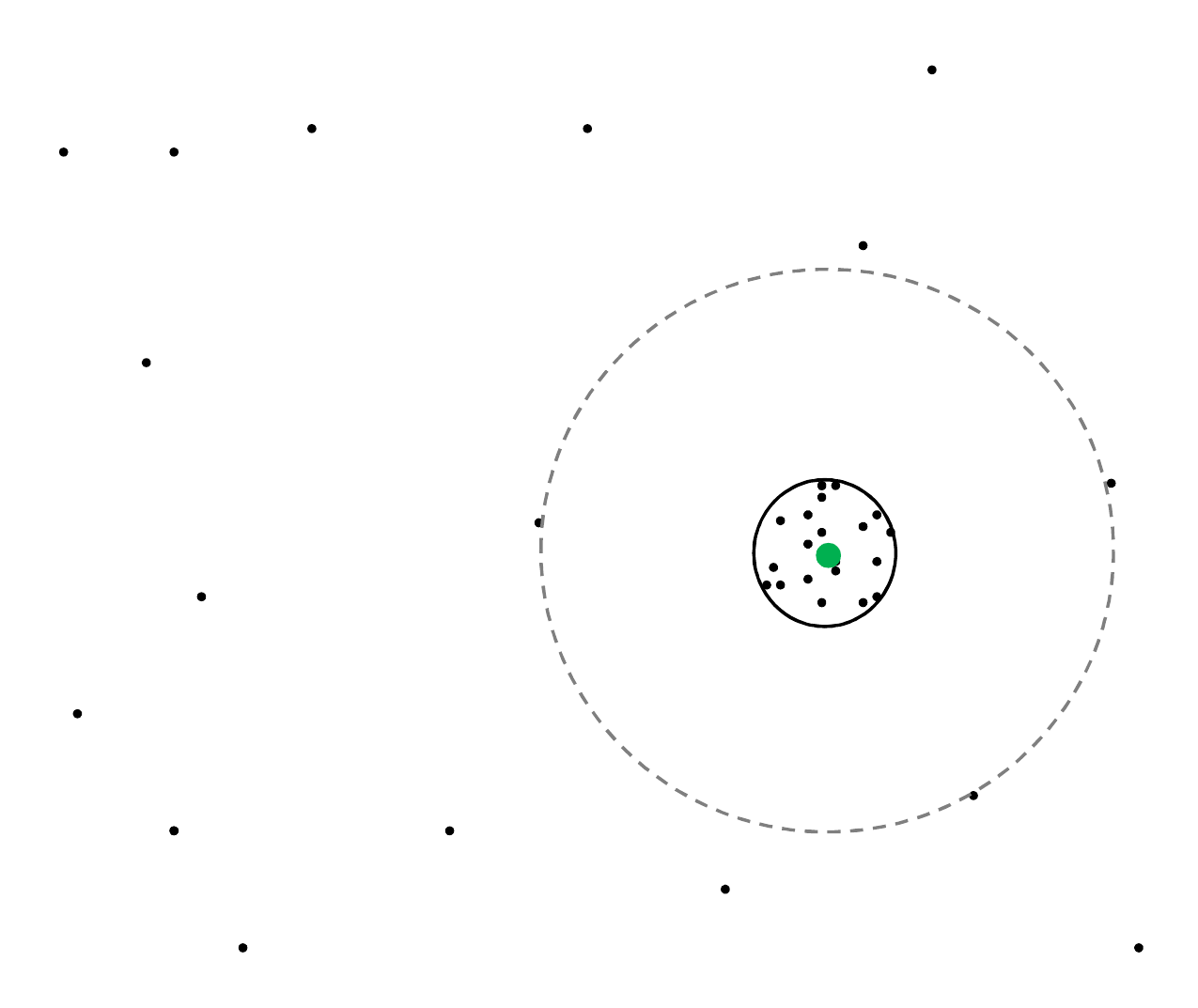} 
  	\caption{Inner search around the new candidate. The robust value here is less than the current global minimum.}
		\label{fig:newalgo6}
    \vspace*{5mm} 
	\end{subfigure}	
	
	\begin{subfigure}{.48\textwidth}
		\centering
		\includegraphics[width=2.2in]{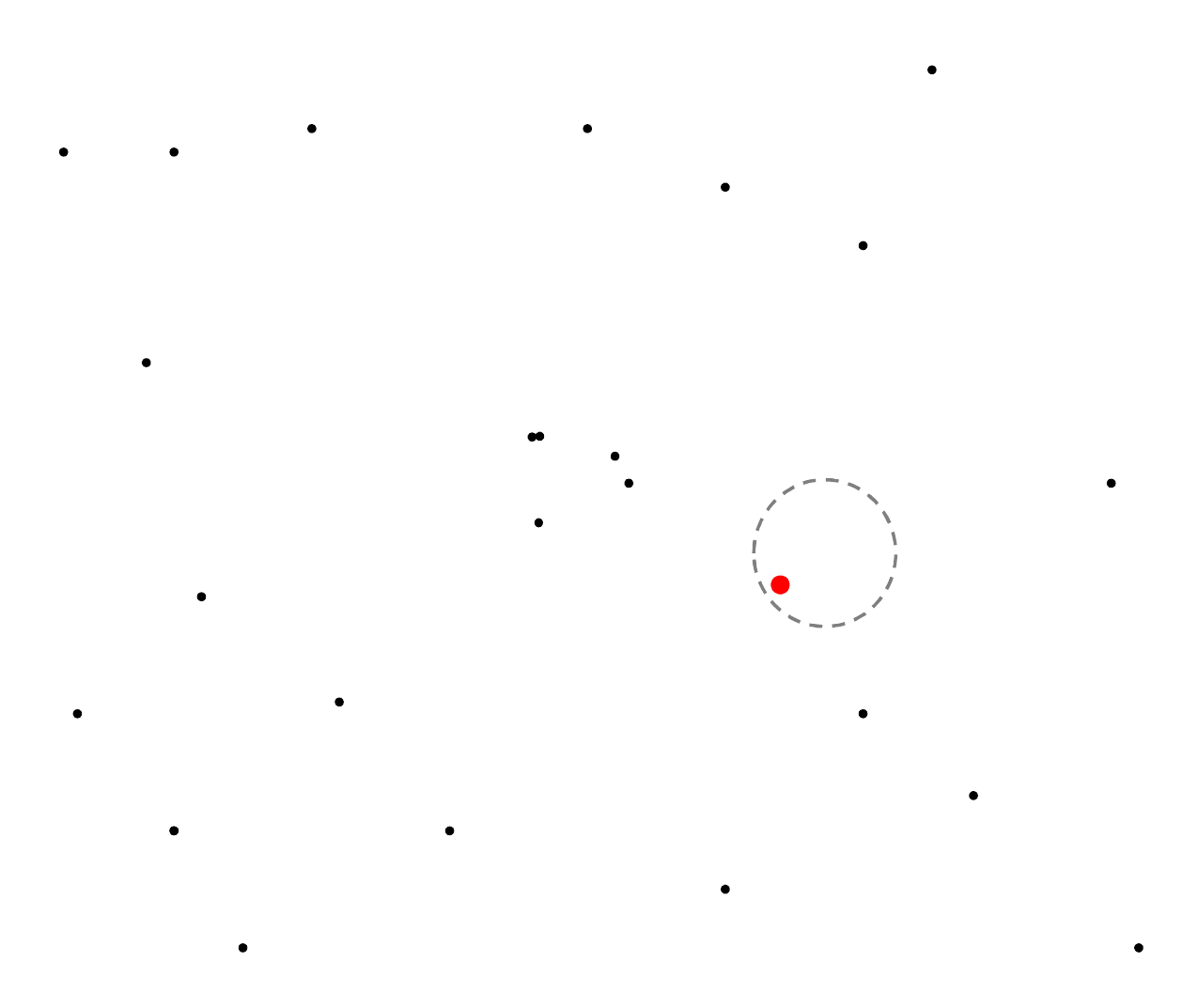} 
		\caption{The current high cost set, including more previously evaluated points due to the reduced high cost threshold.}
		\label{fig:newalgo7}
    \vspace*{5mm} 
	\end{subfigure}%
	\hfill
	\begin{subfigure}{.48\textwidth}
		\centering
		\includegraphics[width=2.2in]{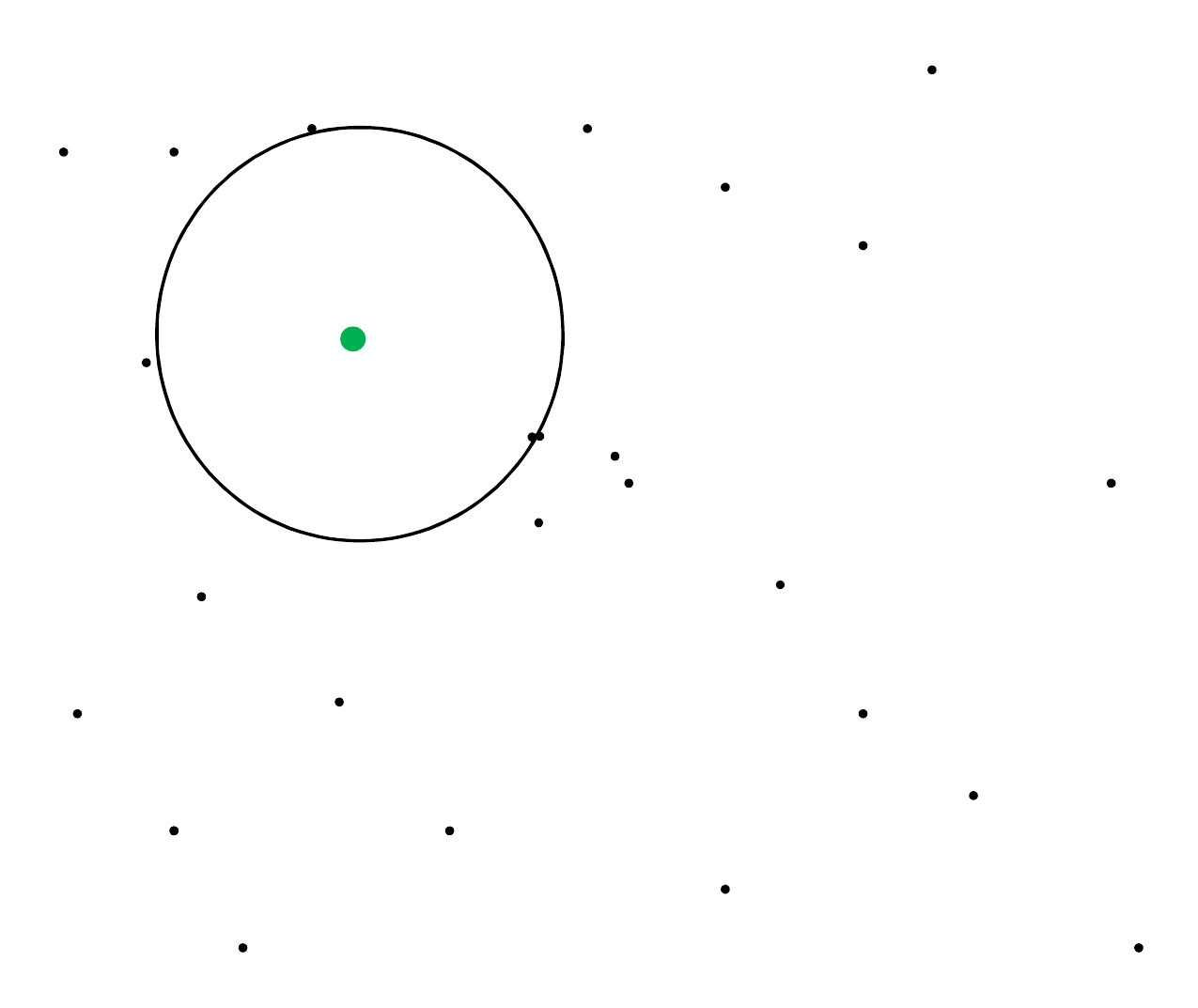} 		
  	\caption{Identify the largest empty hypersphere, the centre of which is the next candidate point.}
		\label{fig:newalgo8}
    \vspace*{5mm} 
	\end{subfigure}	
	
	\caption{Description of largest empty hypersphere (LEH) approach.}
	\label{fig:LEHDescription}
\end{figure}

\noindent Building on the notion of a search that progresses locally by moving away from already identified poor (high cost) points, we develop a global approach that iteratively moves to the region of the decision variable solution space furthest away from recognised hcps. This is an exploration-focused approach, although rather than concentrating on examining unvisited regions the intention here is to identify and visit regions devoid of hcps. Assuming uncertainty as considered previously in terms of a single value $\Gamma$ that defines a radius of uncertainty in all decision variables, we associate the idea of the largest empty region (empty of hcps) with the idea of the largest empty hypersphere (LEH), or largest empty circle in 2D. The approach is then to locate the next point in the search at the centre of the identified LEH, and to iteratively repeat this as more regions are visited and hcps identified. The approach is described in Figure~\ref{fig:LEHDescription}.

We start by randomly sampling one or more points and evaluating the objective function $f$ at each. From these start points a candidate point is selected and an inner analysis undertaken in the candidate's $\Gamma$-uncertainty neighbourhood with a view to identifying the local maximum, Figure~\ref{fig:newalgo3}. This local worst case cost for the candidate is the first estimate of a robust global minimum, that is a global min max, and is located at the candidate point. The aim is now to move to a point whose uncertainty neighbourhood has a lower worst case cost than the current global value. We seek to achieve this by identifying the largest hypersphere of radius at least $\Gamma$ within the defined feasibility bounds which is completely empty of hcps, and moving to the centre of that LEH, see Figures~\ref{fig:newalgo4} -~\ref{fig:newalgo5}.

All points evaluated are recorded in a history set, a subset of which forms the high cost set. The high cost set contains a record of all points evaluated so far with an objective function value greater or equal to a high cost threshold, and here the high cost threshold is set as the current estimate of the robust global minimum. Both the history set and the high cost set are updated as more points are visited and the high cost threshold reduces, see Figures~\ref{fig:newalgo6} -~\ref{fig:newalgo7}. On performing all inner searches after the first candidate, a candidate's robust value may be no better than the current estimate of the robust global minimum (and therefore the current high cost threshold), in which case at least one point will be added to the high cost set. Alternatively if a candidate's robust value is better than the current estimate of the robust global minimum, this current recorded optimum is overwritten and the high cost threshold reduced accordingly. Again this introduces at least one high cost point to the high cost set, but the reducing threshold may also introduce additional points from the history set; this is suggested in Figure~\ref{fig:newalgo7}.

The search stops when no LEH of radius greater than $\Gamma$ exists or some pre-defined resource limit has been reached. Then the candidate point around which the current estimate of the robust global minimum has been determined is deemed the robust global minimum. Otherwise the search repeats, performing analysis in the $\Gamma$-uncertainty neighbourhood around candidates to estimate the local (inner) max, updating the global minimum worst case cost if appropriate, and moving to the next identified LEH, Figure~\ref{fig:newalgo8}.

The critical feature of such an approach is the identification of regions of the solution space that are currently empty of, and furthest from, the undesirable hcps. As defined here this corresponds to identifying the largest hypersphere devoid of all hcps.


Given a discrete history set $H$ of all points evaluated so far, high cost points are those members of $H$ with objective value which is at least the current high cost threshold $\tau$, i.e., 
\[ H_\tau := \{ \pmb{h}\in H \mid f(\pmb{h}) \geq \tau \} \]
We denote $N_{\tau} = |H_\tau|$ as the cardinality of $H_\tau$, and write $H_\tau = \{\pmb{h}^1,\ldots,\pmb{h}^{N_\tau}\}$.
The identification of a point $\pmb{p}\in\X$ which is furthest from all $N_\tau$ high cost points in $H_\tau$ is a max min problem:
\[ \max_{\pmb{p}\in\X} \min_{i\in[N_\tau]} d(\pmb{p}, \pmb{h}^i), \tag{LEHP}\] 
where $d(\pmb{p}, \pmb{q})$ is the Euclidean distance between two points $\pmb{p}$ and $\pmb{q}$, see \cite{OkabeSuzuki1997}. 

In the following, we specify this general LEH approach by considering two aspects in more detail: The \textit{outer search} is concerned with placing the next candidate point $\pmb{x}$ by solving (LEHP). The \text{inner search} then evaluates this candidate by calculating $g(\pmb{x})$ approximately.

\subsection{Outer search methods}

Here we will introduce different approaches to identifying the largest empty hypersphere, given a set of high cost points $H_\tau$. It should be noted that none of these approaches requires additional function evaluations, which is usually considered the limiting resource in black-box settings.

\subsubsection{Randomly sampled LEH algorithm}
\label{secLEH1}
A very simple approach is to generate potential candidates randomly within the feasible region, then determine whether they are more than $\Gamma$ away from all hcps. If so they are a valid candidate, if not re-sample up to some defined maximum number of times beyond which it is assumed that no such candidate point can be found and the solution has converged on a robust global minimum. 
Rather than being a largest empty hypersphere approach this is just a valid empty hypersphere approach, and the size of the identified empty hypersphere might vary considerably from one candidate to the next. 

\subsubsection{Genetic Algorithm for LEH}
\label{secLEH2}
The solution of (LEHP) is an optimisation problem. Furthermore, given a point $\pmb{p}$ which is a potential candidate for the centre of the largest empty hypersphere, the inner minimisation calculation in (LEHP) involves just an enumeration over the $N_\tau$ Euclidean distance calculations between each hcp and $\pmb{p}$ to identify the minimum distance $d(\pmb{p}, \pmb{h}^k)$, where $\pmb{h}^k$ is the closest hcp. Therefore the focus for the solution of (LEHP) is the outer maximisation, for which we may consider an approximate heuristic approach. We employ a genetic algorithm (GA), a commonly cited evolutionary algorithm (EA) \cite{Ghazali2009}. Here each individual represents a point $\pmb{p}$ in the decision variable space, and the objective function $f_{LEH}(\pmb{p}):=\min_{\pmb{h}\in H_\tau} d(\pmb{p},\pmb{h})$ is the minimum distance between a given point $\pmb{p}$ and all hcps in $H_\tau$. We seek to maximise this minimal distance by evolving a population of points starting from randomly selected feasible points in the decision variable space $\X$. The best point generated by the GA is the next candidate point -- that is estimated centre of the LEH, for the current $H$, $\tau$ and $H_\tau$.

\subsubsection{Voronoi based LEH}
\label{secLEH3}

Within the literature a widely referenced approach for tackling low dimensional LEH problems is due to \cite{Toussaint1983}, and is based on the geometric Voronoi diagram approach, see \cite{Chazelle1993, OkabeSuzuki1997}. The Voronoi approach partitions a space into regions (cells). For a given set of points each cell corresponds to a single point such that no point in the cell is closer to any other point in the set. Points on the edges between cells are equidistant between the set points which lie on either side of that edge. For our LEH problem the set of points is $H_\tau$, and the Voronoi diagram approach corresponds to segmenting the feasible space $\X$ into $N_\tau$ separate cells, one for each hcp. The (Voronoi) vertices that lie at the intersection of these cell (Voronoi) edges maximise the minimum distance to the nearby set points, see \cite{Chazelle1993, OkabeSuzuki1997}. So for a given $H_\tau$ if we can determine the Voronoi diagram we can use the identified Voronoi vertices as potential candidate points $\pmb{p}$. The solution of (LEHP) is then simply a matter of enumeration, for each $\pmb{p}$ calculating the (inner) minimum Euclidean distance to all hcps, and then selecting the (outer) maximum such minimal distance.

The original approach due to \cite{Toussaint1983} includes the identification of vertices (candidate centres of LEHs) that can be sited outside of defined boundaries, in infeasible regions. This is not exactly as required here. To deal with this edges that cross feasibility boundaries are identified and the associated vertices which are outside of $\X$ are relocated to an appropriate point on the boundary of $\X$. Here any coordinate $i\in[n]$ of such an external vertex that is either less than $l_i$ or greater than $u_i$ is re-set to $l_i$ or $u_i$ as appropriate.

However the Voronoi approach has exponential dependence on $n$, as constructing the Voronoi diagram of $N_\tau$ points requires $O(N_\tau log N_\tau + N_\tau^{\ceil*{n/2}})$ time \cite{Chazelle1993}. This suggests that such an approach in not computationally viable for anything other than low dimensional problems. On the basis that a Voronoi diagram based approach is the primary recognised heuristic for identifying the largest empty \textit{circle} we will consider a Voronoi based robust LEH heuristic here only in the context that for 2D problems in our experimental analysis this approach will serve as a good direct comparator for our other robust LEH heuristics.

\subsection{Inner search methods}

Discussions of the LEH approach have so far focussed on the outer minimisation search, assuming some form of inner search that provides the inner robust maximum for each candidate point in the minimisation search. In \cite{BertsimasNohadaniTeo2010} a two-stage gradient ascent search is recommended for each inner search around a candidate point. This assumes gradient information is available and proposes $(n+1)$ individual two-stage gradient ascents for each candidate. For a 100-dimensional problem this would require several thousand function evaluations around each candidate point. In practical terms both the number of function evaluations required to undertake a global search and the requirement for gradient information may make such extensive inner searches prohibitive. Given, for example, budgetary restrictions on the number of function evaluations, some trade-off must be achieved between the extent of each inner $\Gamma$-radius uncertainty neighbourhood search and globally exploring the search space. But this trade-off between robustness in terms of the extent of the inner searches, and performance in terms of the outer global search, is complex, see \cite{MirjaliliLewisMostaghim2015, DiazHandlXu2017}. For example the determination of an appropriate inner approach -- type of search, extent of search and parameter settings -- may be both instance (problem and dimension) dependent and dependent on the outer approach. 

Here we do not propose to recommend a definitive inner search approach. From a theoretical point of view we assume the information is provided by some oracle. From an experimental point of view in the algorithm testing and comparisons below we assume the same basic inner $\Gamma$-radius uncertainty neighbourhood analysis for all heuristics, to ensure a consistency when comparing results for alternative search approaches.

There is, however, an aspect of our LEH approach that enables an additional feature, the forcing of an early end to an inner search.
The LEH approach is exploration-led, the objective being to locate and move to the candidate point in the decision variable space furthest from all hcps. Hcps are designated based on the determination of a high cost threshold $\tau$, set here as the current estimate of the robust global minimum (min max) value. The nature of this approach enables (inner) uncertainty neighbourhood searches around each candidate point to be restricted when appropriate. If an inner search identifies a local point with objective function value above $\tau$ the inner search can be immediately curtailed on the basis that the candidate is not distant from hcps. This equates to the recognition that the candidate point is not an improvement on the current estimated robust optima. Such regulating of inner searches has the potential to significantly reduce the number of function evaluations expended on local neighbourhood analysis. In the case of budgetary limitations on numbers of function evaluations this further enables more exploration of the decision variable space.

\subsection{Algorithm summary}

Given one of our three approaches to identifying the LEH devoid of hcps, random, GA or Voronoi, the overarching algorithm for the robust exploratory LEH heuristic is given in Algorithm~\ref{LEHAlgorithm}. Here one of these three approaches to the outer search
is applied in line~\ref{LEHSubRoutine} as $LEH\textunderscore Calculator(H_\tau)$, for a defined high cost set $H_\tau$. It is assumed that this routine will return a candidate point $\pmb{x}_{LEH}$ and an associated radius $r_{LEH}$, that is the minimal distance between $\pmb{x}_{LEH}$ and all points in $H_\tau$. The heuristic will halt if $r_{LEH}$ is not greater than $\Gamma$.

For a defined number of initialisation points, random points in $\X$ are selected and the function $f$ evaluated at these points. The points and their function evaluations are recorded in history sets $H$ and $F_H$, lines~\ref{InitialisationStart} -~\ref{InitialisationEnd}. Having randomly selected a candidate point $\pmb{x}_c$ from $H$ we perform an inner maximisation in the $\Gamma$-uncertainty neighbourhood around $\pmb{x}_c$, see line~\ref{InnerSubRoutine}. The description of the inner maximisation is given below as Algorithm~\ref{InnerAlgorithm}. 
If this is the first candidate point, or the local robust value for this candidate $\tilde{g}(\pmb{x}_c)$ is less than the current best solution $\tau$, this minimum is updated and the associated global minimum point $\pmb{x}_{Op}$ replaced by $\pmb{x}_c$, see lines~\ref{RobustGlobalStart} -~\ref{RobustGlobalEnd}. 

Next the high cost set $H_\tau$ is established as all members of $H$ with corresponding function values in $F_H$ that are greater than or equal to the current high cost threshold $\tau$, see line~\ref{HighCost}. Based on $H_\tau$, the next candidate point is identified via one of the outer search approaches, see line~\ref{LEHSubRoutine}. If the heuristic is halted at this stage due to an inability to identify a valid LEH
or at any stage due to the budget being exceeded, the extant estimate for the robust global minimum $\pmb{x}_{Op}$ is returned.


\begin{algorithm}[htbp] 
\caption{Robust global exploration using Largest Empty Hyperspheres} \label{LEHAlgorithm}
\vspace{2mm} 
\hspace*{\algorithmicindent} \textbf{Input:} $f$, $\X$, $\Gamma$ \\
\hspace*{\algorithmicindent} \textbf{Parameters:} $Num\textunderscore Initial$, $Budget$, $Max\textunderscore Search$
\vspace{2mm} 

\begin{algorithmic}[1]
\ForAll{$i$ in $[Num\textunderscore Initial]$}  \label{InitialisationStart}
	\State Choose random point $\pmb{x}^i\in\X$
	\State Calculate $f(\pmb{x}^i)$ and store in $F_H$
	\State $Budget \gets Budget-1$
	\State $H \gets H \cup \{\pmb{x}^i\}$
\EndFor  \label{InitialisationEnd} 

\State Select random point $\pmb{x}_c\in H$
\State $r_{LEH} \gets \infty$; $\tau \gets \infty$

\While{$r_{LEH} > \Gamma$}

	\State $\tilde{g}(\pmb{x}_c) \gets$ \textbf{CALL} Algorithm~\ref{InnerAlgorithm} \label{InnerSubRoutine}
	\If{$\tilde{g}(\pmb{x}_c) < \tau$}  \label{RobustGlobalStart} 
		\State $\pmb{x}_{Op} \gets \pmb{x}_c$
		\State $\tau \gets \tilde{g}(\pmb{x}_c)$		
	\EndIf  \label{RobustGlobalEnd} 
		
	\State $H_{\tau} \gets \{ \pmb{x}\in H : F_H(\pmb{x}) \ge \tau \}$ \label{HighCost}

	\State Find ($\pmb{x}_{LEH},r_{LEH})$ by calling LEH\textunderscore Calculator($H_\tau$)
	 \label{LEHSubRoutine}
	
		\State $\pmb{x}_c \gets \pmb{x}_{LEH}$
	
\EndWhile

\State \Return A robust solution $\pmb{x}_{Op}$ and robust objective estimate $\tau$ \label{EndAlgorithm}
\end{algorithmic}
\end{algorithm}


\begin{algorithm}[htbp] 
\caption{$\Gamma$-uncertainty neighbourhood inner maximisation} \label{InnerAlgorithm}
\vspace{2mm} 
\hspace*{\algorithmicindent} \textbf{Input:} $Budget$, $Max\textunderscore Search$, $\pmb{x}_c$, $\Gamma$, $\tau$
\vspace{2mm} 
 
\begin{algorithmic}[1]
\If{$\tau<\infty$} 
       \State Calculate $f(\pmb{x}_c)$ and store in $F_H$
	\State $H \gets H \cup \{\pmb{x}_c\}$
	\State $Budget \gets Budget-1$ 
	\If{$Budget==0$} 
		\State GOTO line~\ref{EndAlgorithm} of Algorithm 1
	\EndIf
\EndIf 

\State Set $Local\textunderscore Robust \gets f(\pmb{x_c})$

\ForAll{$i$ in $[Max\textunderscore  Search]$} \label{InnerStart}
	\State Choose $\Delta\pmb{x}^i_c \in \U$, set $\pmb{x}^i \gets \pmb{x}_c + \Delta\pmb{x}^i_c$ \label{RandomInHyper}
\State Calculate $f(\pmb{x}^i_c)$ and store in $F_H$
\State $H \gets H \cup \{\pmb{x}^i_c\}$
	\State $Budget \gets Budget-1$ 
	\If{$Budget==0$} 
		\State GOTO line line~\ref{EndAlgorithm} of Algorithm 1
	\EndIf
	
	\State $Local\textunderscore Robust \gets \max\{Local\textunderscore Robust,f(\pmb{x}^i_c)\}$ \label{UpdateRobust}
	
	\If{$Local\textunderscore Robust>\tau$}  \label{ThresholdStart} 
		\State GOTO line~\ref{DefineOutput}	
	\EndIf  \label{ThresholdEnd}	
\EndFor  \label{InnerEnd} 

\State $\tilde{g}(\pmb{x}_c) \gets Local\textunderscore Robust$  \label{DefineOutput}

\State \Return $\tilde{g}(\pmb{x}_c)$: estimated worst case cost at $\pmb{x}_c$
\end{algorithmic}
\end{algorithm}

\medskip

Algorithm~\ref{InnerAlgorithm}, the $\Gamma$-uncertainty neighbourhood inner maximisation called in line~\ref{InnerSubRoutine} of Algorithm~\ref{LEHAlgorithm}, requires several inputs: $Budget$ the current count of function evaluations completed, $Max\textunderscore Search$ the maximum number of function evaluations permitted in an inner search, $\pmb{x}_c$ the current candidate point (centre of an LEH) around which the inner search is to be performed, $\Gamma$ to define the uncertainty neighbourhood of $\pmb{x}_c$, and $\tau$ the high cost threshold for stopping the inner search if appropriate. 

Algorithm~\ref{InnerAlgorithm} proceeds by looping through up to $Max\textunderscore Search$ inner search points, identifying a point in the $\Gamma$-uncertainty neighbourhood of $\pmb{x}_c$ and evaluating the function at each point visited, lines~\ref{InnerStart} -~\ref{InnerEnd}. Here the point to be evaluated is determined by random sampling in the $\Gamma$-radius hypersphere centred on $\pmb{x}_c$, line~\ref{RandomInHyper}. Under other inner maximisation rules this would be determined by some explicit maximisation search heuristic. As the function is evaluated at the inner search points the local robust value (inner maximum) $Local\textunderscore Robust$ is updated as appropriate, line~\ref{UpdateRobust}. If $Local\textunderscore Robust$ exceeds the high cost threshold $\tau$ the inner maximisation is immediately terminated, lines~\ref{ThresholdStart} -~\ref{ThresholdEnd}. Algorithm~\ref{InnerAlgorithm} ends by returning an estimate for the worst case cost value at $\pmb{x}_c$, $\tilde{g}(\pmb{x}_c)$ into Algorithm~\ref{LEHAlgorithm}.

\subsection{Example LEH application}
\label{sec:NatureLEHSearch}

\begin{figure}[htbp]
	\centering
	\begin{subfigure}{.38\textwidth} 
		\centering
		\includegraphics[width=2.5in, height=2.6in]{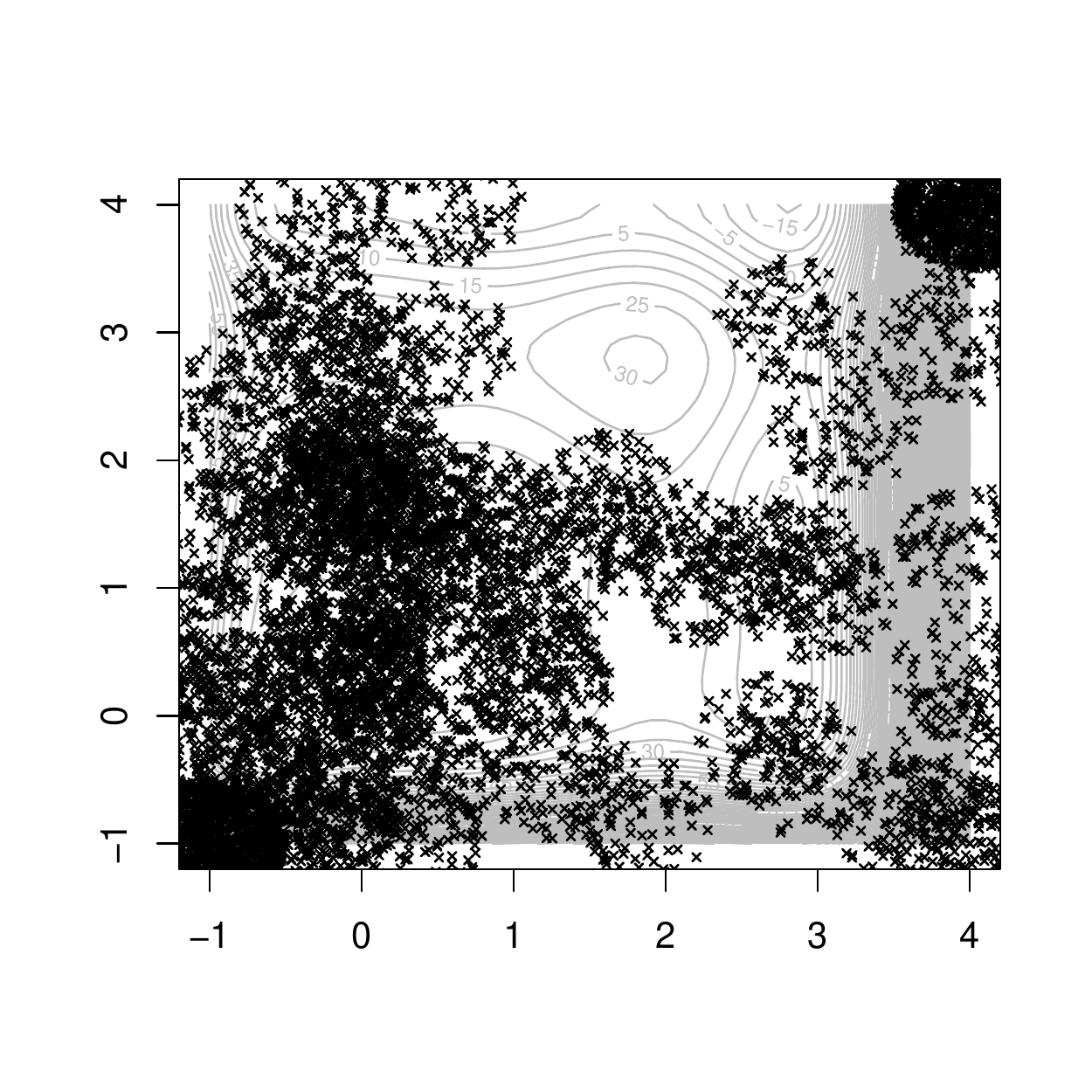}
		\vspace{-11mm} 
  	\caption{PSO points}
		\label{fig:PSOPoints}
	\end{subfigure}%
	\hspace{7mm} 
	\begin{subfigure}{.38\textwidth}
		\centering
		\includegraphics[width=2.5in, height=2.6in]{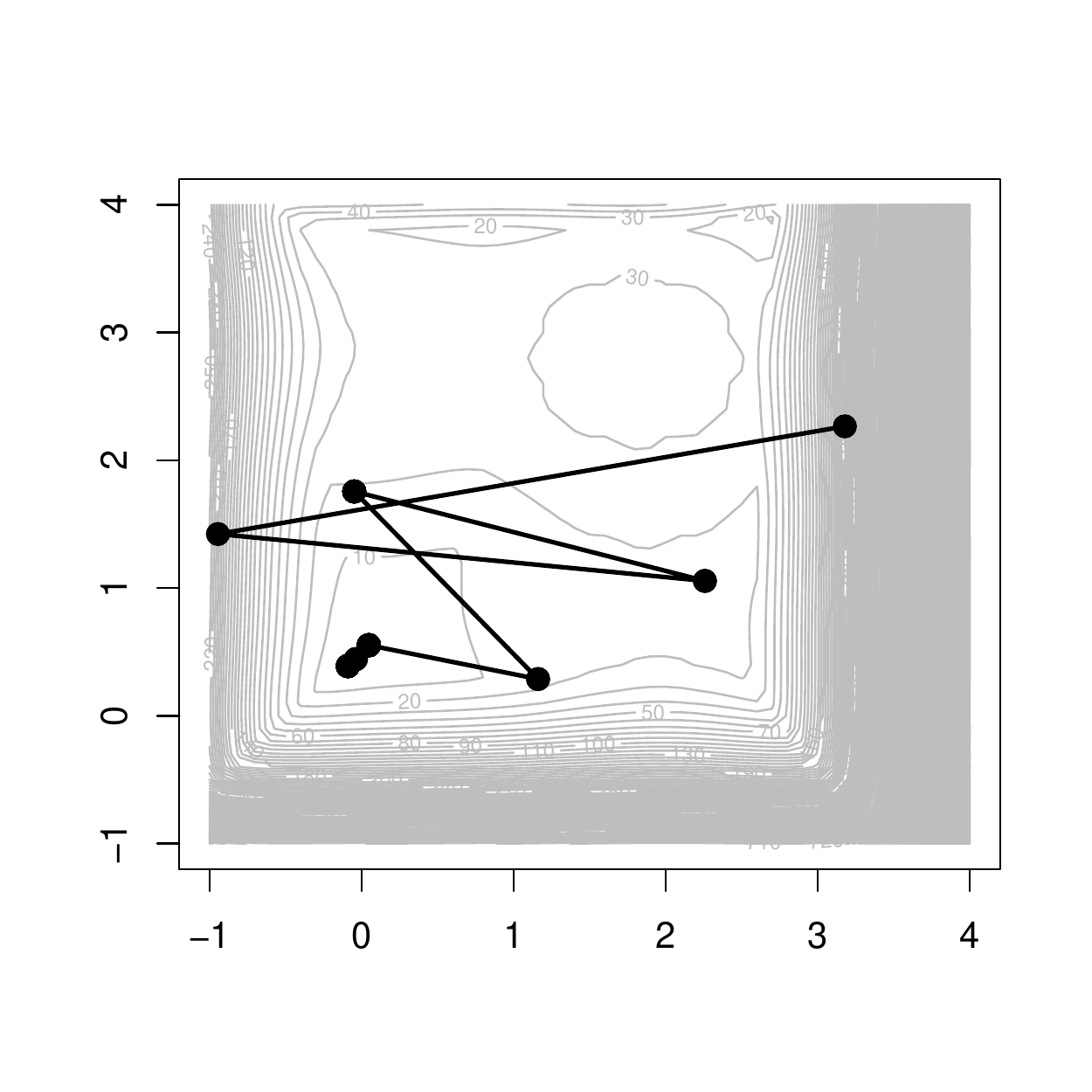}	
		\vspace{-11mm} 
  	\caption{PSO search}
		\label{fig:PSOSearch}
	\end{subfigure}
		
	\vspace{-5mm} 
	
	\begin{subfigure}{.38\textwidth}
		\centering
		\includegraphics[width=2.5in, height=2.6in]{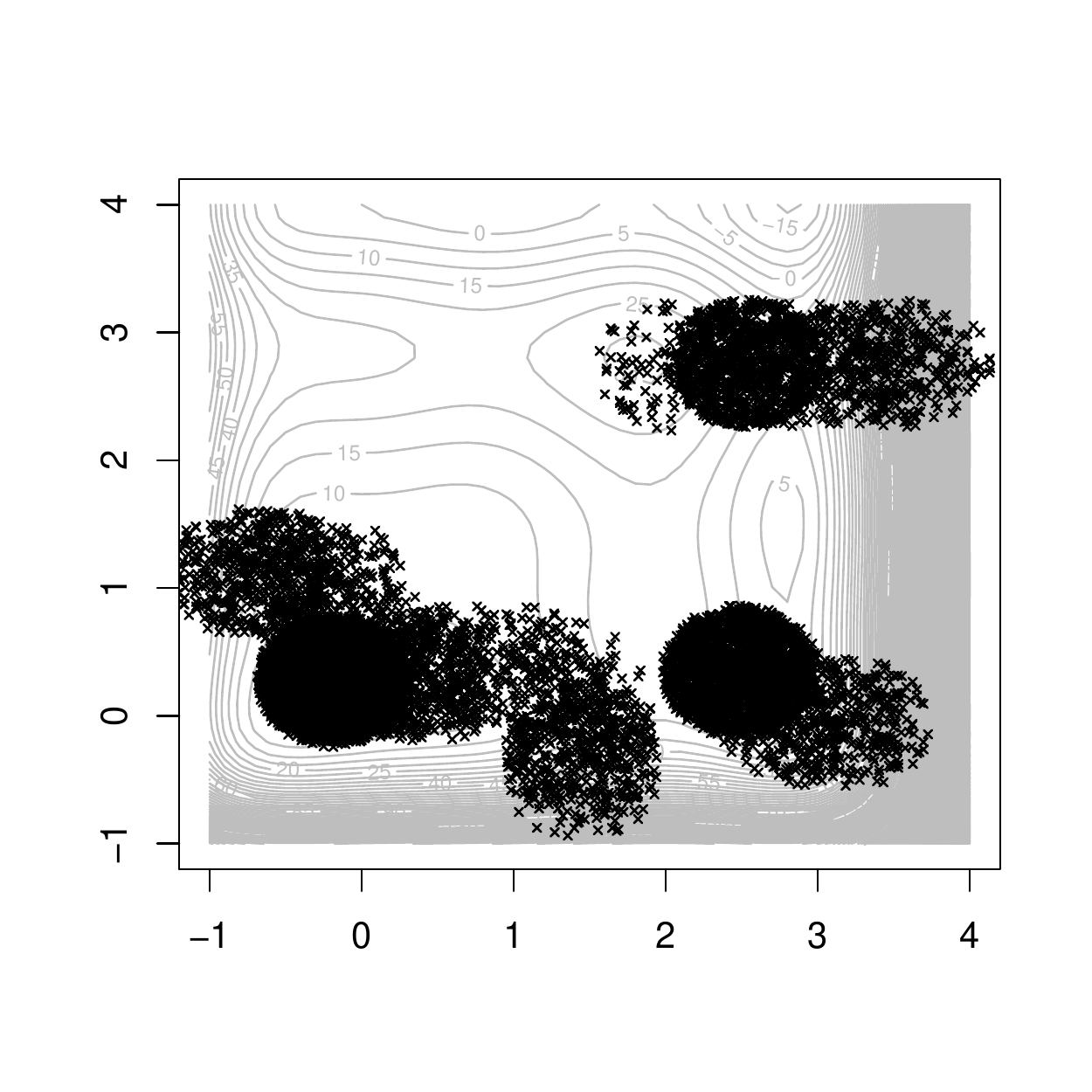}	
		\vspace{-11mm} 
		\caption{DD points}
		\label{fig:DesDirPoints}
	\end{subfigure}%
	\hspace{7mm} 
	\begin{subfigure}{.38\textwidth}
		\centering
		\includegraphics[width=2.5in, height=2.6in]{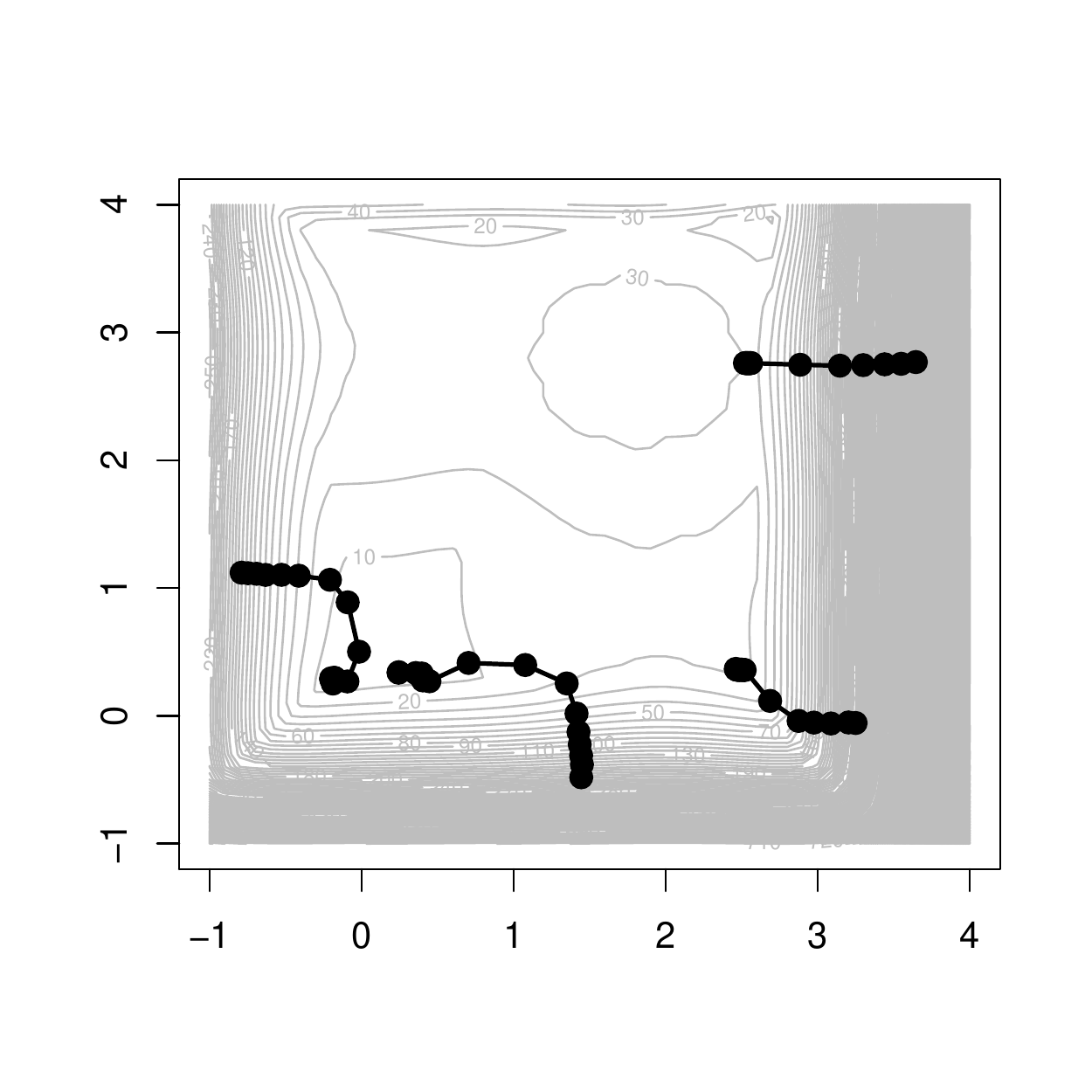}	
		\vspace{-11mm} 
  	\caption{DD search}
		\label{fig:DesDirSearch}
	\end{subfigure}	
		
	\vspace{-5mm} 
			
	\begin{subfigure}{.38\textwidth}
		\centering
		\includegraphics[width=2.5in, height=2.6in]{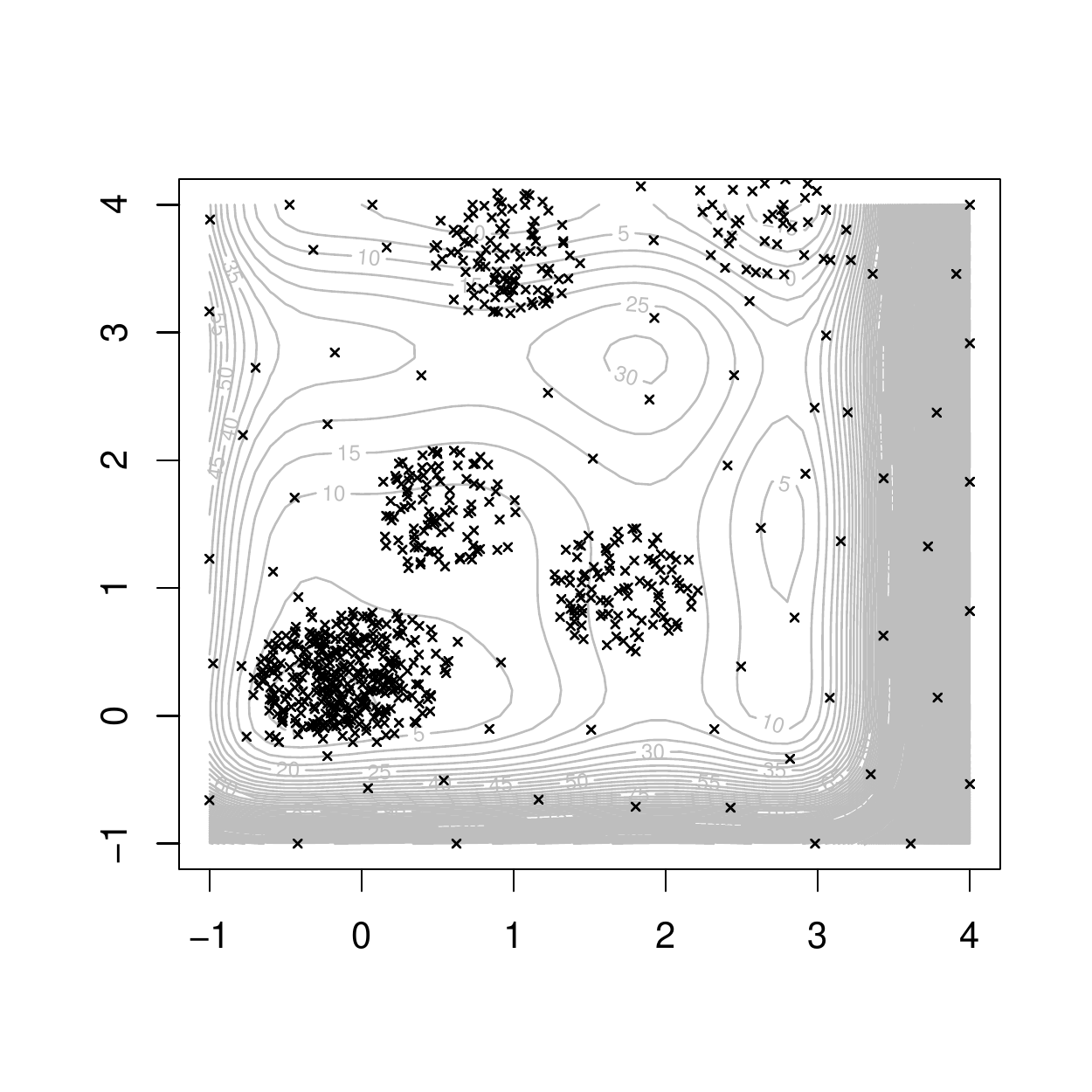}		
		\vspace{-11mm} 
	  \caption{LEH Vor points}	
		\label{fig:LEHVoronoiPoints}
	\end{subfigure}%
	\hspace{7mm} 
	\begin{subfigure}{.38\textwidth}
		\centering
		\includegraphics[width=2.5in, height=2.6in]{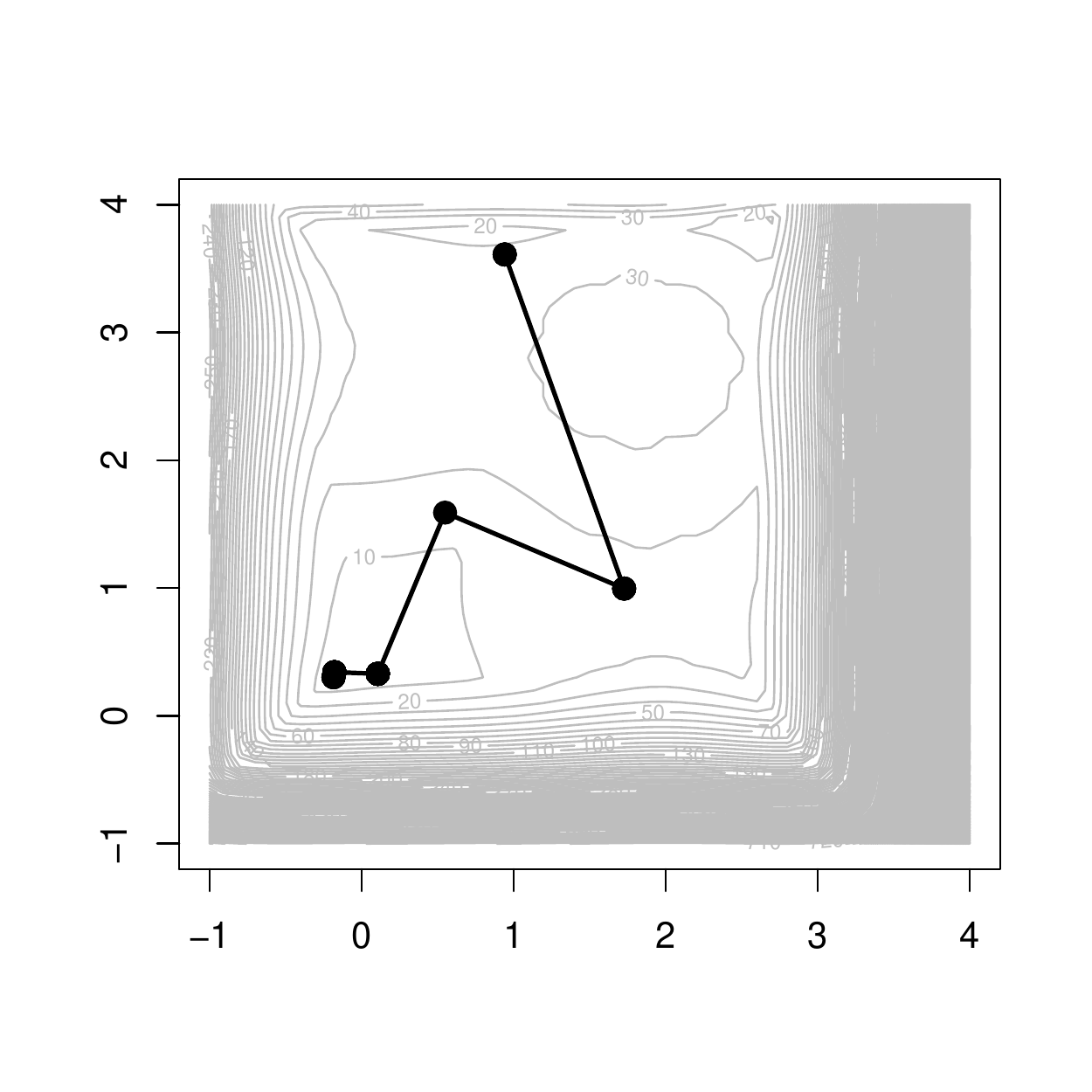}		
		\vspace{-11mm} 
  	\caption{LEH Vor search}
		\label{fig:LEHVoronoiSearch}
	\end{subfigure}	
	
	\caption{Contour plots of example searches of the 2-dimensional problem (poly2D), for $\Gamma$=0.5. Plots on the left show all points evaluated. Plots on the right show the progress of the current best robust solution. The heuristics used are: (top) outer PSO, (middle) outer descent directions with re-start, and (bottom) outer LEH using the Voronoi based approach.}
	\label{fig:AlternativeRobust2DSearches}
\end{figure}

In order to give some indication of the nature of our LEH search we have applied it to the 2-dimensional problem (poly2D) and plotted the points evaluated and associated search path of the current estimate of the robust global minimum in Figures~\ref{fig:LEHVoronoiPoints} and~\ref{fig:LEHVoronoiSearch}. Here the LEH Voronoi algorithm is used. For comparison we have also plotted corresponding results for two alternative heuristics, a robust Particle Swarm Optimisation (PSO) approach shown in Figures~\ref{fig:PSOPoints} and~\ref{fig:PSOSearch}, and the local descent directions approach from Section~\ref{sec:LocalRobustSearchDescentDirections} shown in Figures~\ref{fig:DesDirPoints} and~\ref{fig:DesDirSearch}. Here the robust PSO is used as a proxy to a brute force or co-evolutionary approach. The basic global PSO formulations have been used, as described in \cite{ShiEberhart1998}. The descent directions approach has been extended by using random re-starts, as a proxy to extending it to a global approach. In all cases inner random sampling in a hypersphere of 100 $\Gamma$-uncertainty neighbourhood points is used, and a maximum budget of 10,000 function evaluations employed. 

The plots shown in Figure~\ref{fig:AlternativeRobust2DSearches} are for only a single run of each heuristic, and as such should only be seen as exemplars intended to give some indication of the different natures of these outer search approaches. It can be seen that whilst the robust PSO explores the decision variable space somewhat, and the re-starting descent directions follows (exploits) a series of local paths, the LEH approach features both considerable exploration globally and more intense analysis of promising points. It is clear that the curtailing of the inner searches in the LEH approach enables much wider exploration for fewer function evaluations. In this example less than 1,000 function evaluations have been required before the LEH heuristic has stopped because an LEH of radius greater than $\Gamma$ cannot be found, but for larger (dimensional) problems such stopping prior to reaching the budgetary limit will not apply. One striking feature of Figure~\ref{fig:LEHVoronoiPoints} is how many of the inner searches stop immediately on the evaluation of a candidate point. This is because the objective value at these candidate points exceeds the current threshold $\tau$.

The Voronoi based search exemplified by Figures~\ref{fig:LEHVoronoiPoints} and~\ref{fig:LEHVoronoiSearch} is a good indicator of the nature of the searches due to all three LEH approaches, random, GA and Voronoi. However the radii of the LEH identified for each candidate will vary with the use of each of these algorithms. Figure~\ref{fig:AlternativeLEHSearchesRadii} in Appendix~\ref{sec:RadiiLEH} gives some indication of how the radii of the hyperspheres generated by each of these LEH heuristics progress as the exploration proceeds.


\section{Computational experiments}
\label{sec:ExperimentsResults}

\subsection{Set up}
\label{sec:SetUp}

In order to assess the effectiveness of the LEH approach the heuristic has been applied to eight test problems, and results compared against the two alternative search heuristics described in Section~\ref{sec:ComparatorHeuristics}. Experiments have been performed on 2D, 4D, 7D, 10D and 100D instances of these test problems; results have also been generated for (poly2D). Both the genetic algorithm and random forms of the LEH heuristic have been assessed for all instances. The LEH Voronoi has additionally been applied to the 2D instances, with the intention of giving some indication of the differences due to a `best' LEH identifier algorithm (Voronoi) versus the alternatives. All LEH approaches are initialised by randomly sampling a single point in $\X$. Assuming that for most real-world problems the optimisation analysis will be limited by resources, a fixed budget of 10,000 function evaluations (model runs) is assumed. The same inner approach is employed for all heuristics. A simple random sampling in a hypersphere of 100 points in a point's local $\Gamma$-uncertainty neighbourhood is used for all instances, and the local robust maximum is estimated as the maximum due to this sampling. For the LEH approaches this inner sampling is curtailed if a point is identified in the uncertainty neighbourhood that has objective value exceeding the current high cost threshold $\tau$.

All experiments have have been performed using Java, on an HP Pavilion 15 Notebook laptop computer, with 64 bit operating system, an Intel Core i3-5010U, 2.10GHz processor, and 8GB RAM. Each heuristic search has been applied to each test problem-dimension instance 50 times to reduce variability. For the solution of the Second Order Cone Problem as part of the descent directions algorithm \cite{BertsimasNohadaniTeo2010}, the IBM ILOG CPLEX Optimization Studio V12.6.3 package is called from Java.

\subsection{Comparator heuristics}
\label{sec:ComparatorHeuristics}

Our experiments have been conducted on LEH, a re-starting descent directions, and robust PSO metaheuristics. We have applied parameter tuning to 3 of the 5 comparator heuristics -- LEH Voronoi and LEH Random do no have tunable parameters -- employing an evolutionary tuning approach using a genetic algorithm to generate a single set of parameters for each heuristic, for all test problems. For each of the 3 tuned heuristics the same subset of the test instances was used, running each member of an evolving population on each of these instances multiple times to generate mean result for each member of a population on each test instance. The performance  of each individual in a population was ranked separately for each test instance, across the members of the population, leading to mean overall ranks which were used as the utility measure in tournament selection; see e.g. \cite{Eiben2012}. 
	
The effectiveness of the local descent directions approach \cite{BertsimasNohadaniTeo2010} suggests that extending this to a global search by using random re-starts will provide a reasonable comparator. A local descent directions search is undertaken from a random start point, and when this is complete it is repeated from another random start point. This is repeated until the function evaluations budget is reached. In descent directions a set of high cost points leads to the identification of an optimal stepping direction and step size, if a valid direction exists. However the algorithm includes a number of dynamically changing parameters which adapt the high cost set and enforce a minimum step size. Here we have tuned 5 parameters relating to these stages of the heuristic; see \cite{BertsimasNohadaniTeo2010} for further information. Labelled `\ddre' in the results section.

As a proxy to a brute force or co-evolutionary approach an outer particle swarm search is considered. The basic formulations for the global PSO approach have been used as described in \cite{ShiEberhart1998} and 5 parameters have been tuned: swarm size, number of iterations, and for the velocity equation the $C_1$ and $C_2$ acceleration parameters and inertia weight parameter $\omega$. The combined swarm size times number of iterations was limited to 100 in order to align with the budget of 10,000 function evaluations and the level of inner sampling. Labelled `PSO' in the results section.

Our robust LEH metaheuristic is considered for the three alternative ways of identifying the largest hypersphere that is empty of hcps:

\begin{itemize}[leftmargin=*] 
	\item Randomly sampled valid empty hypersphere, see Section~\ref{secLEH1}. This includes re-sampling up to 1,000 potential candidates in an attempt to identify a valid empty hypersphere, otherwise it is assumed that a valid point cannot be found and a robust global minimum has been reached. Labelled `LEH Rnd' in the results section.
	\item Genetic algorithm LEH, see Section~\ref{secLEH2}. Here we have tuned 6 parameters: the size of the population, number of generations, number of elites, tournament size, and mutation probability and size; we have fixed the use of tournament selection and the choice of mid-point crossover. The combined population size times number of generations was limited to 100, which is somewhat based on runtime considerations associated with the large value of $N_\tau$, the number of candidate points visited with a budget of 10,000 function evaluations. Labelled `LEH GA' in the results section.
	\item Voronoi based \cite{Toussaint1983} LEH, see Section~\ref{secLEH3}. Here the construction of the Voronoi diagram for the input points $H_\tau$ is performed using the Java library due to \cite{Nahr2017}. This generates geometric data, Voronoi vertices and edges, which are used to determine a set of potential candidate points -- Voronoi vertices, including those originally outside of $\X$ relocated to the boundary of $\X$ -- for the centre of the LEH. Labelled `LEH Vor' in the results section.
\end{itemize} 

\subsection{Test functions}
\label{sec:TestFunctions}

\noindent 
A large number of test functions are available for benchmarking optimisation algorithms, and posing a variety of difficulties, see \cite{Kruisselbrink2012, JamilYang2013}. Here eight are considered, plus (poly2D) as outlined in Section~\ref{sec:ProblemDescription}. In each case a single $\Gamma$-uncertainty value is used:

\begin{itemize}[leftmargin=*]  
	\item Ackleys: feasible region [-32.768, 32.768]; $\Gamma$=3.0.
	\item Multipeak F1: feasible region [0, 1]; $\Gamma$=0.0625.
	\item Multipeak F2: feasible region [0, 10]; $\Gamma$=0.5.	
	\item Rastrigin: feasible region [-5.12, 5.12]; $\Gamma$=0.5.
	\item Rosenbrock: feasible region [-2.048, 2.048]; $\Gamma$=0.25.
	\item Sawtooth: feasible region [-1, 1]; $\Gamma$=0.2.
	\item Sphere: feasible region [-5, 5]; $\Gamma$=1.0.
	\item Volcano: feasible region [-10, 10]; $\Gamma$=1.5.
\end{itemize} 

\begin{figure}[htbp]
	\captionsetup[subfigure]{font=scriptsize,labelfont=scriptsize}
	\centering
	\begin{subfigure}{.20\textwidth}
		\centering
		\includegraphics[width=1.5in, height=1.7in]{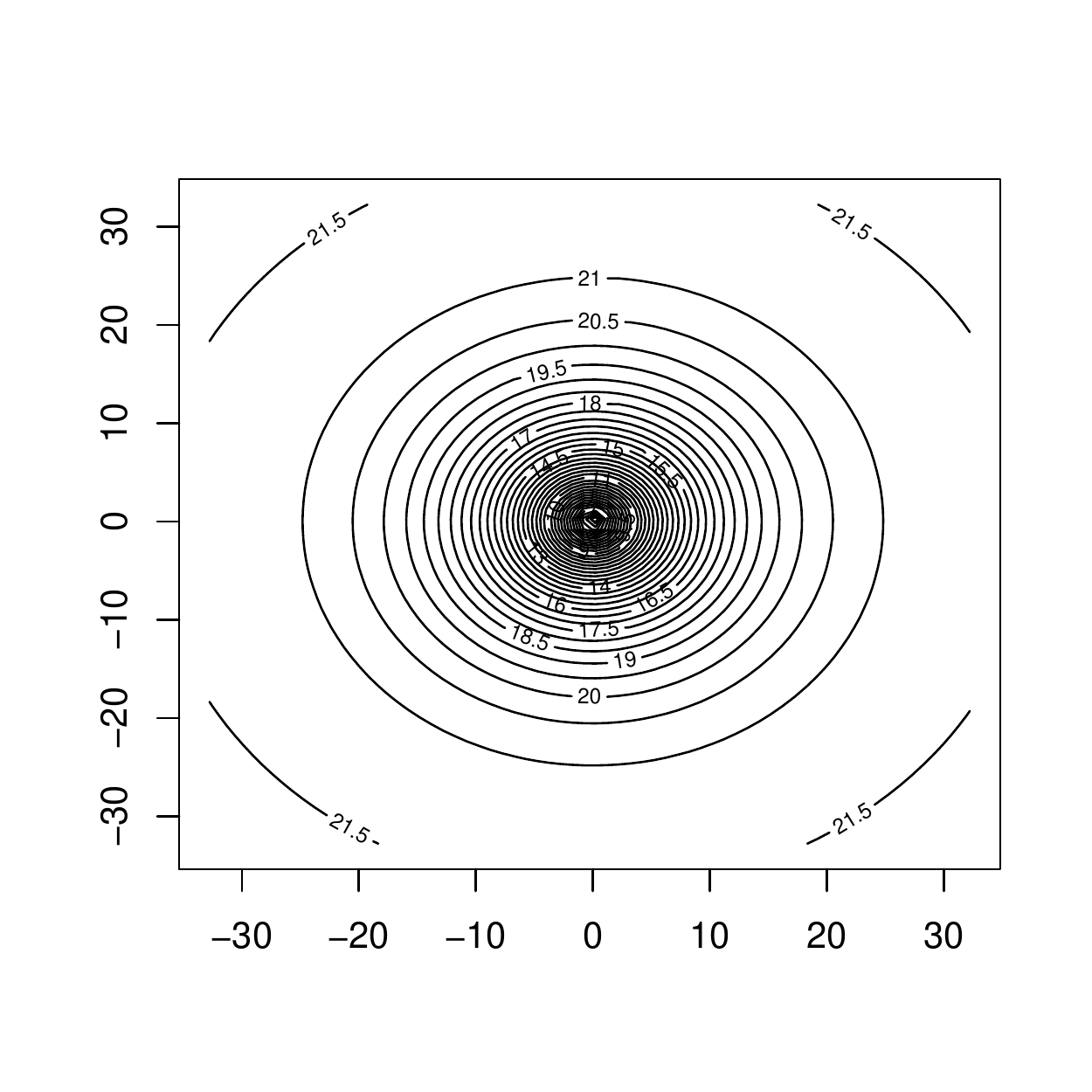}
		\label{fig:AckleyContour}
	\end{subfigure}
	\hspace{4mm} 
	\begin{subfigure}{.20\textwidth}
		\centering
		\includegraphics[width=1.5in, height=1.7in]{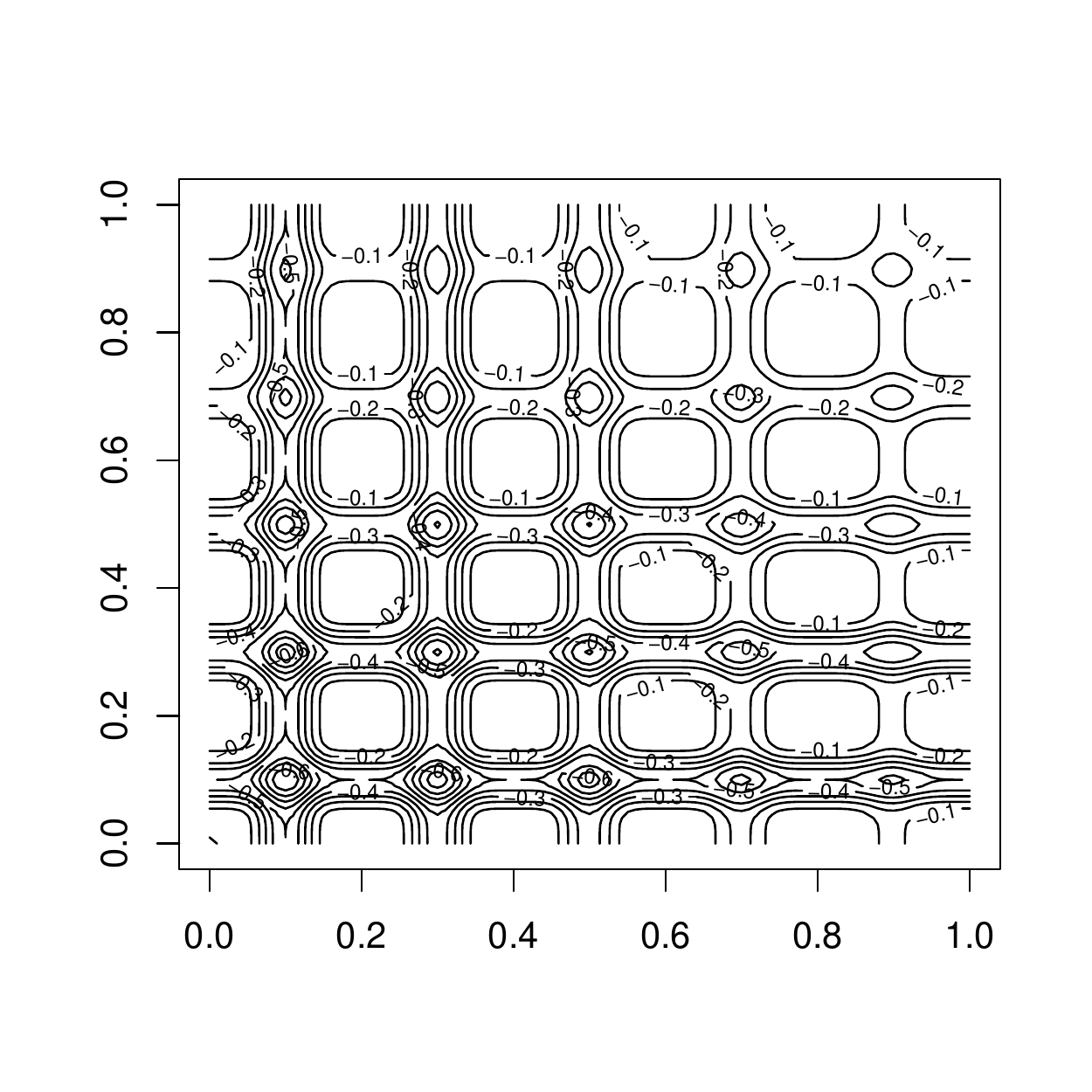}
		\label{fig:MultipeakF1Contour}
	\end{subfigure}
	\hspace{4mm} 
	\begin{subfigure}{.20\textwidth}
		\centering
		\includegraphics[width=1.5in, height=1.7in]{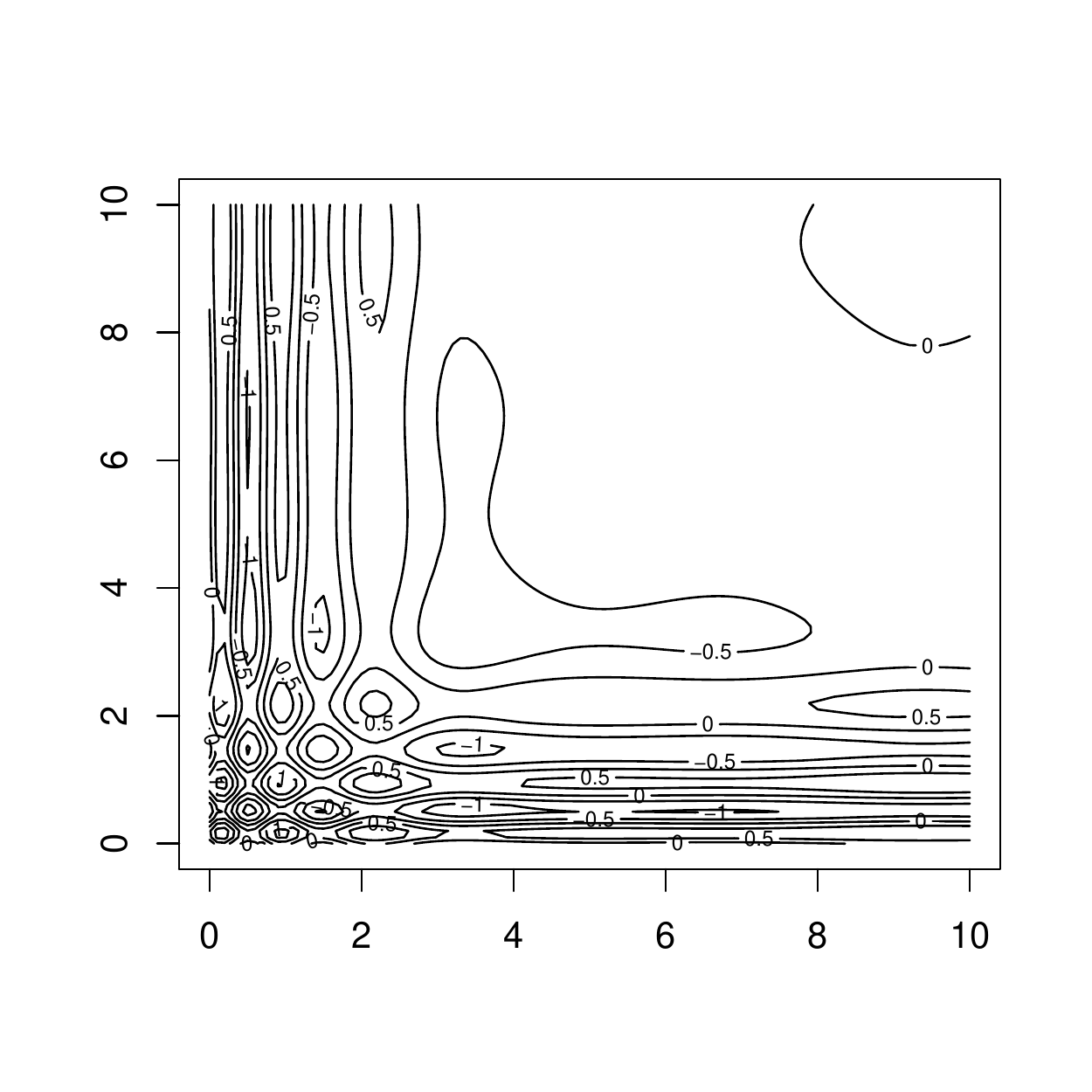}
		\label{fig:MultipeakF2Contour}
	\end{subfigure}	
	\hspace{4mm} 
	\begin{subfigure}{.20\textwidth}
		\centering
		\includegraphics[width=1.5in, height=1.7in]{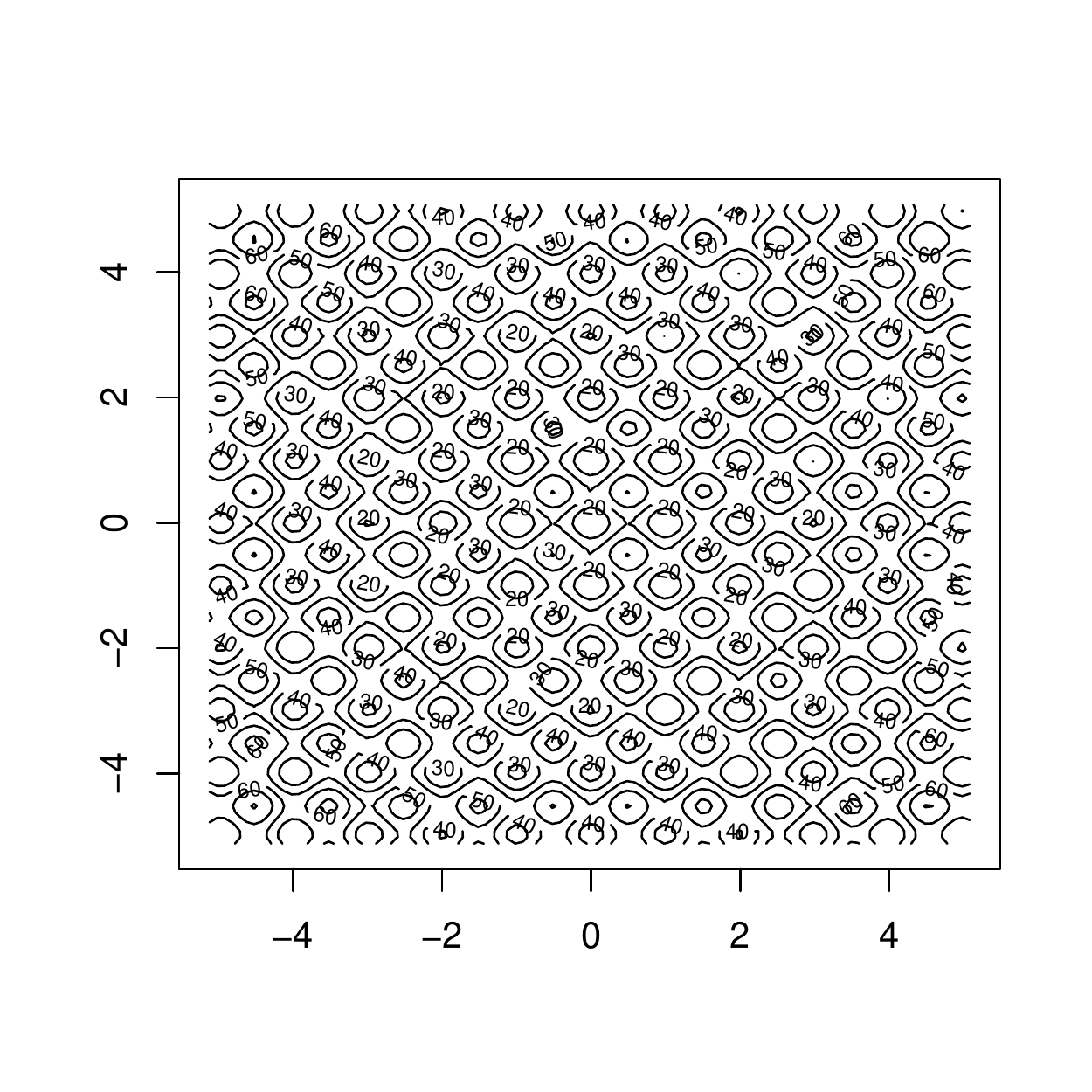}
		\label{fig:RastriginContour}
	\end{subfigure}%
	
	\vspace{-10mm} 
		
	\begin{subfigure}{.20\textwidth}
		\centering
		\includegraphics[width=1.5in, height=1.7in]{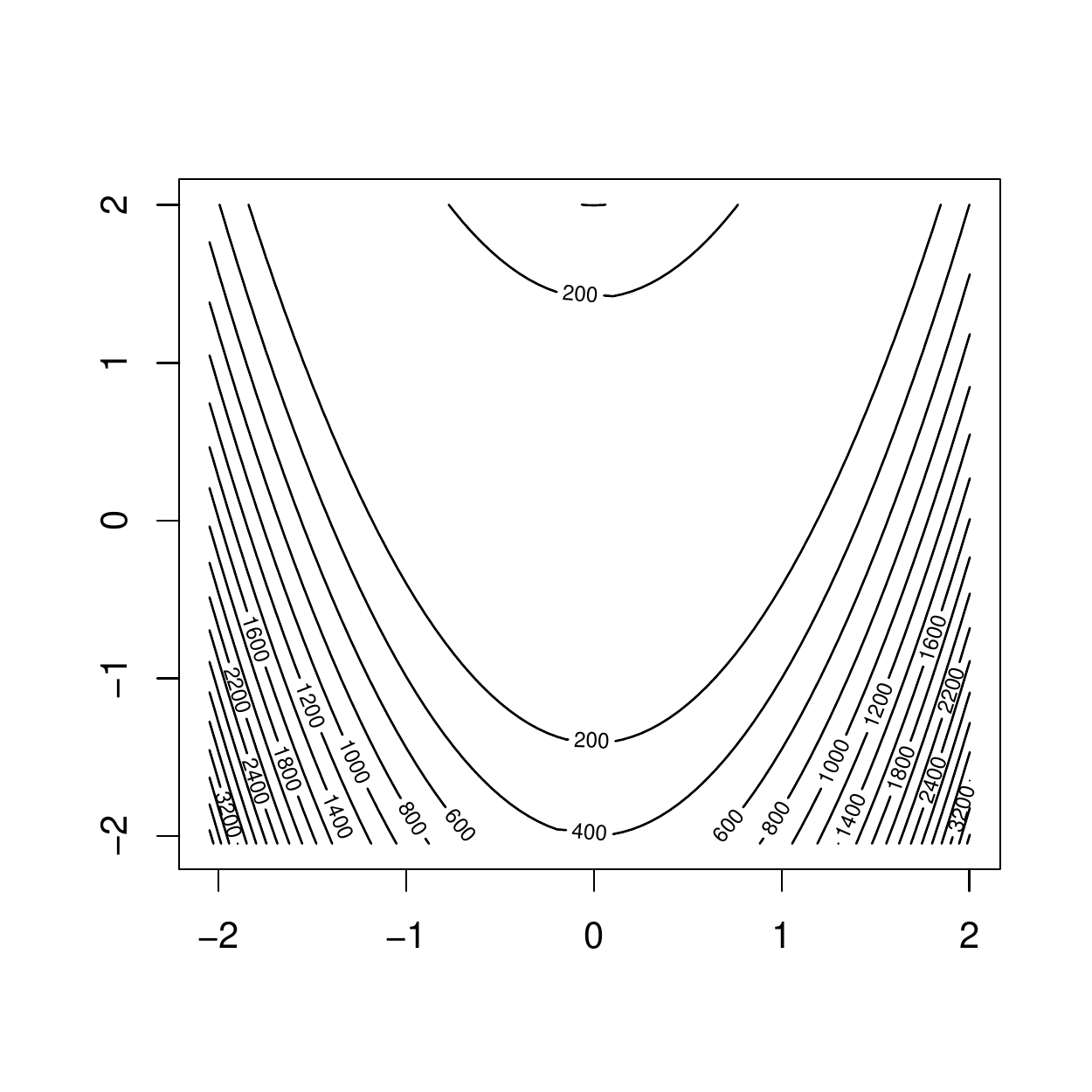}
		\label{fig:RosenbrockContour}
	\end{subfigure}	
	\hspace{4mm} 
	\begin{subfigure}{.20\textwidth}
		\centering
		\includegraphics[width=1.5in, height=1.7in]{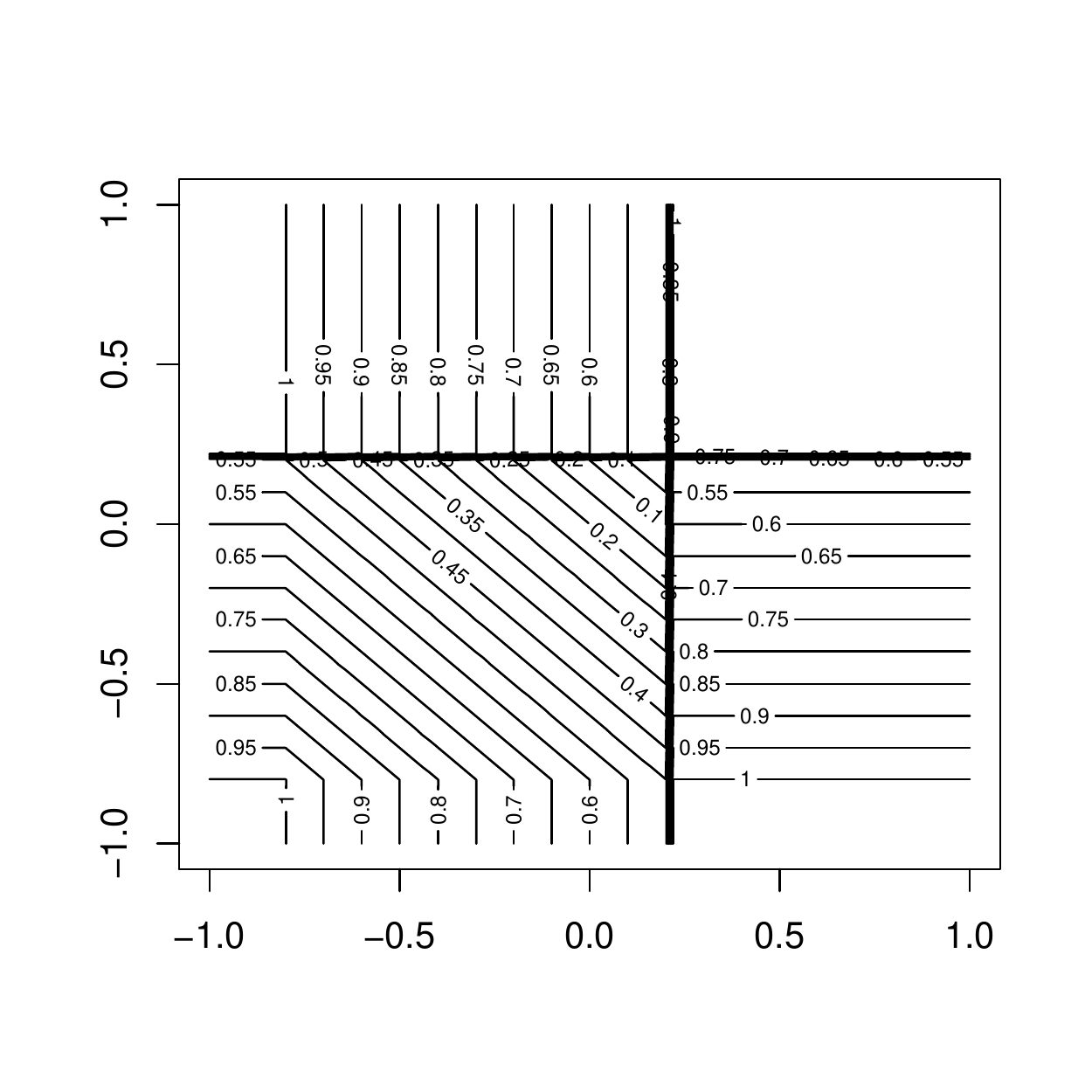}
		\label{fig:SawtoothContour}
	\end{subfigure}
	\hspace{4mm} 
	\begin{subfigure}{.20\textwidth}
		\centering
		\includegraphics[width=1.5in, height=1.7in]{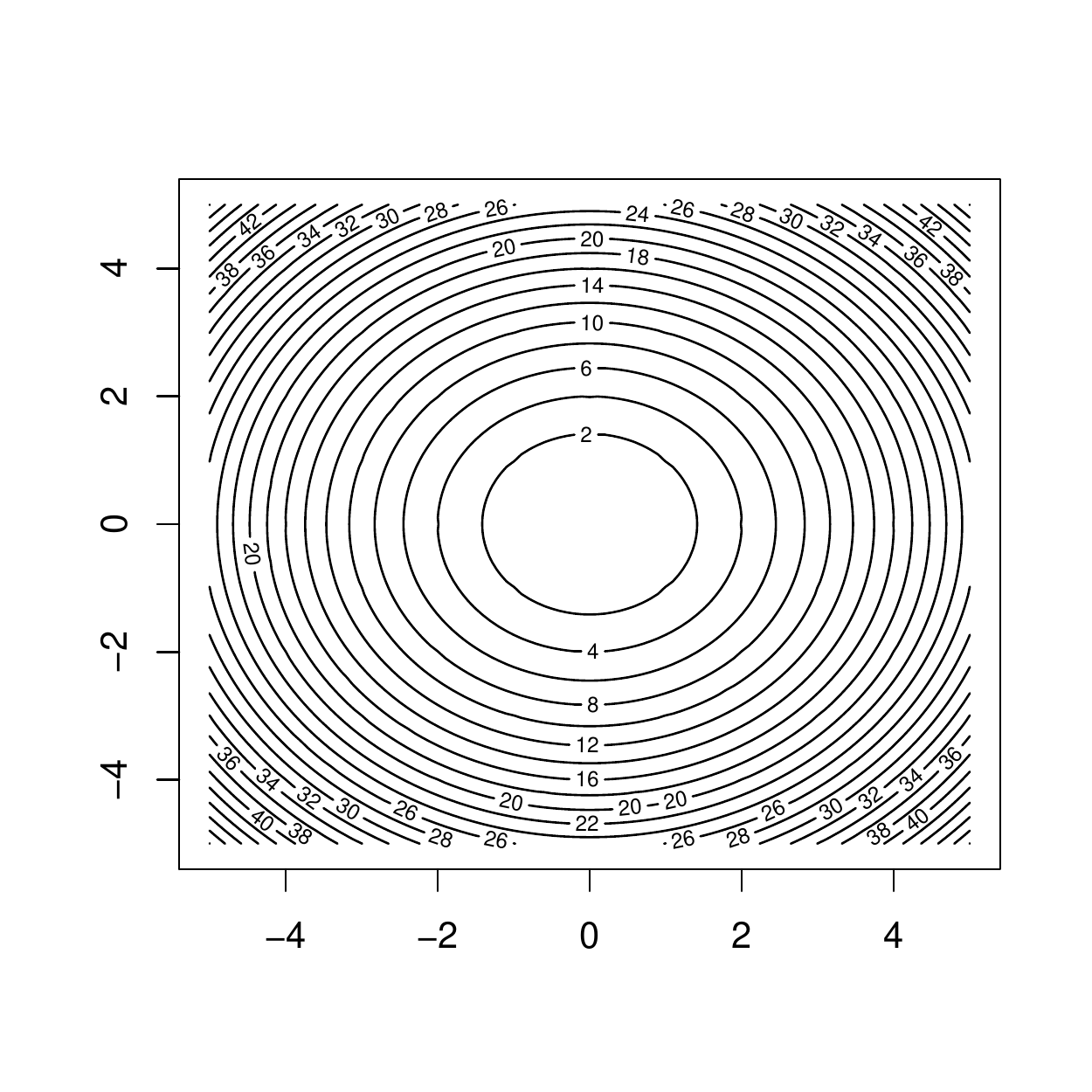}
		\label{fig:SphereContour}
	\end{subfigure}	
	\hspace{4mm} 
	\begin{subfigure}{.20\textwidth}
		\centering
		\includegraphics[width=1.5in, height=1.7in]{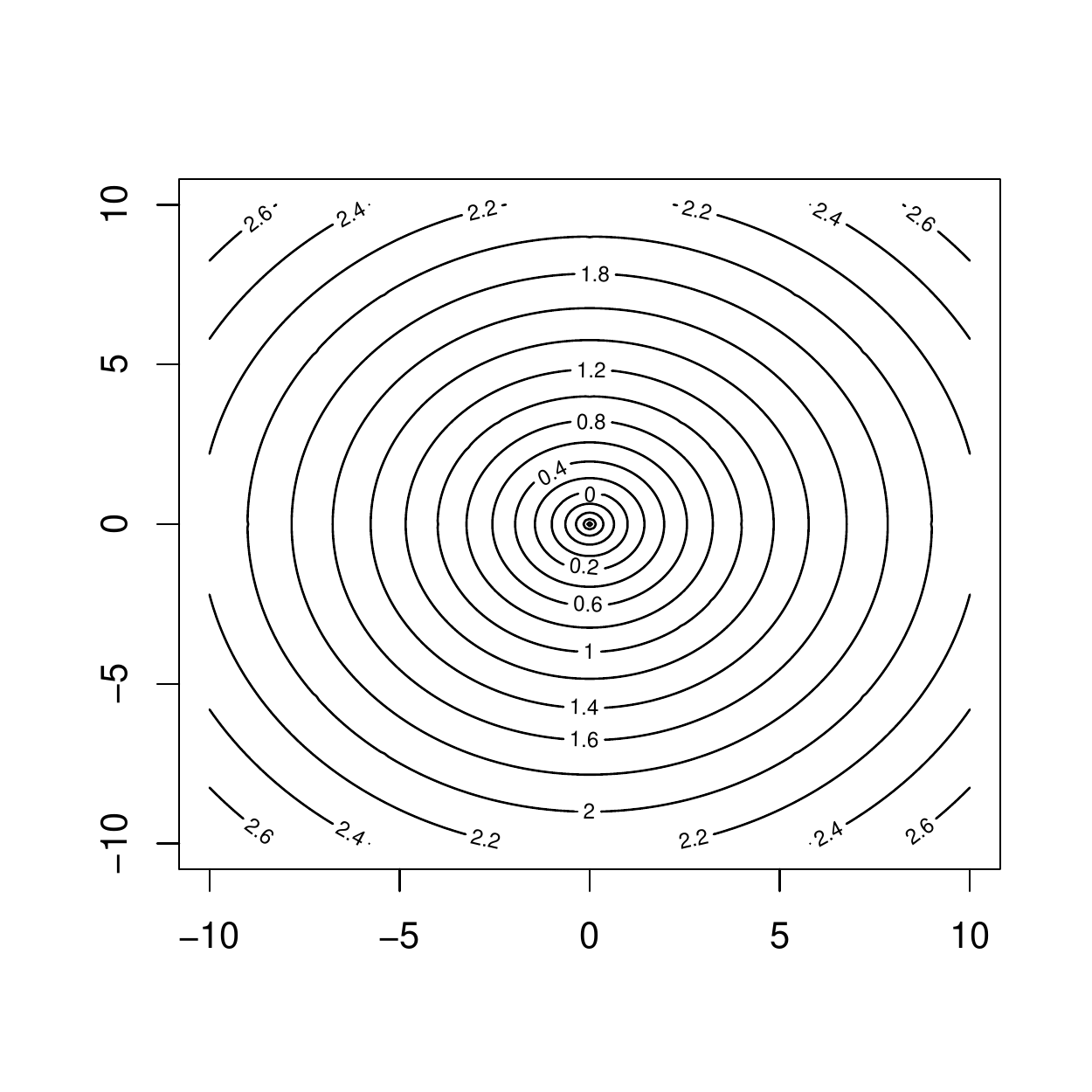}
		\label{fig:VolcanoContour}
	\end{subfigure}%
	
	\vspace{-10mm} 
	
	\begin{subfigure}{.20\textwidth}
		\centering
		\includegraphics[width=1.5in, height=1.7in]{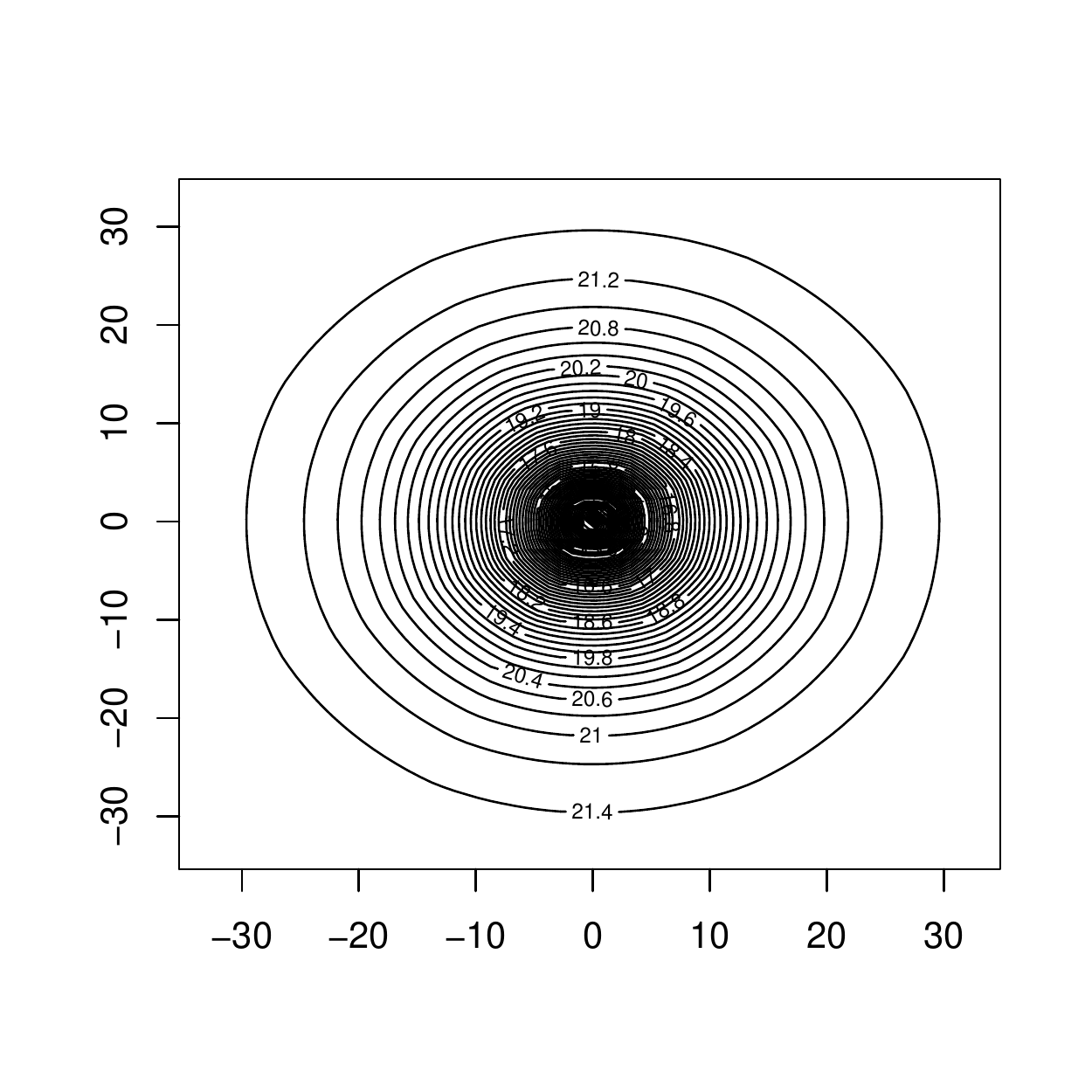}
		\label{fig:AckleyWorstContour}
	\end{subfigure}
	\hspace{4mm} 
	\begin{subfigure}{.20\textwidth}
		\centering
		\includegraphics[width=1.5in, height=1.7in]{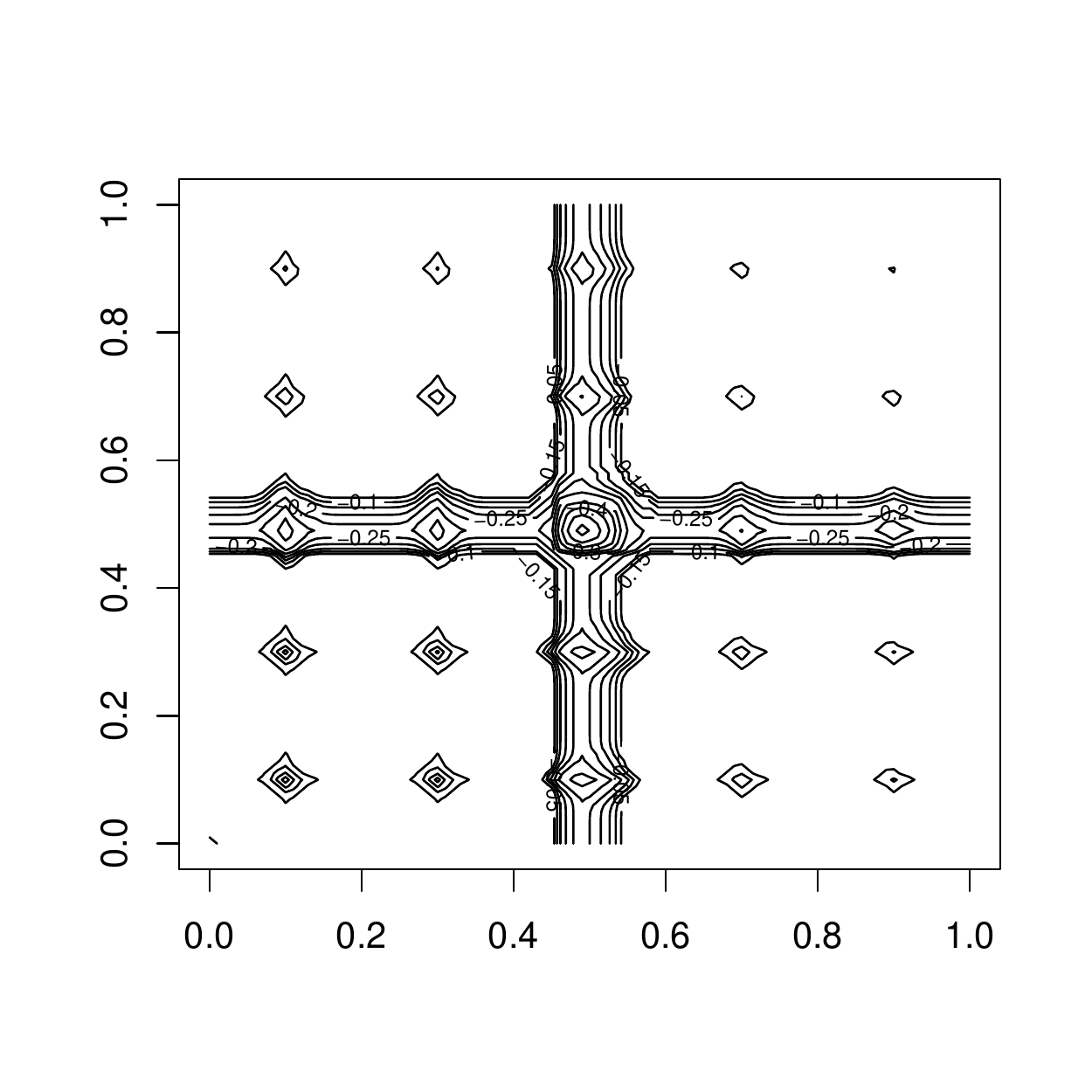}
		\label{fig:MultipeakF1WorstContour}
	\end{subfigure}
	\hspace{4mm} 
	\begin{subfigure}{.20\textwidth}
		\centering
		\includegraphics[width=1.5in, height=1.7in]{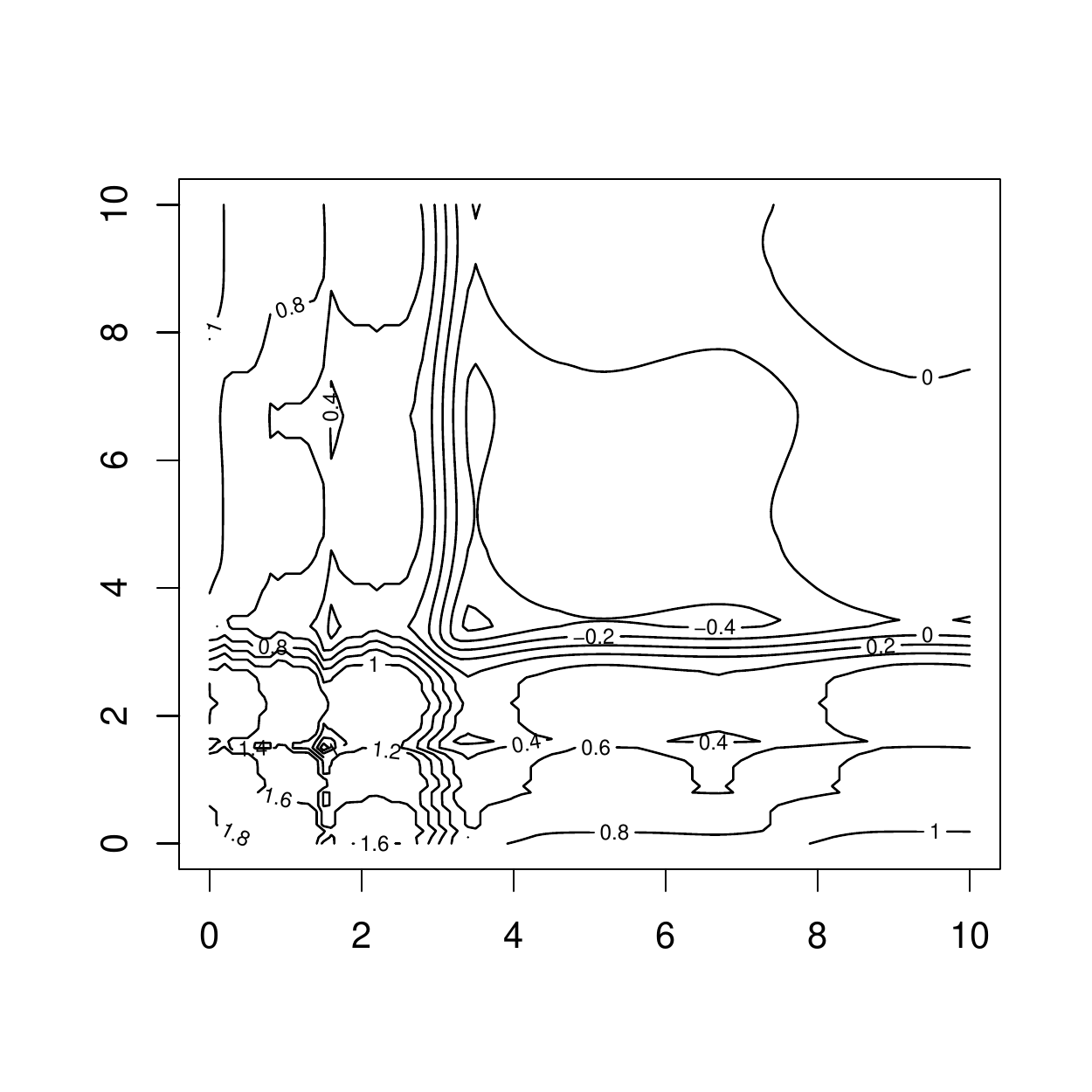}
		\label{fig:MultipeakF2WorstContour}
	\end{subfigure}
	\hspace{4mm} 
	\begin{subfigure}{.20\textwidth}
		\centering
		\includegraphics[width=1.5in, height=1.7in]{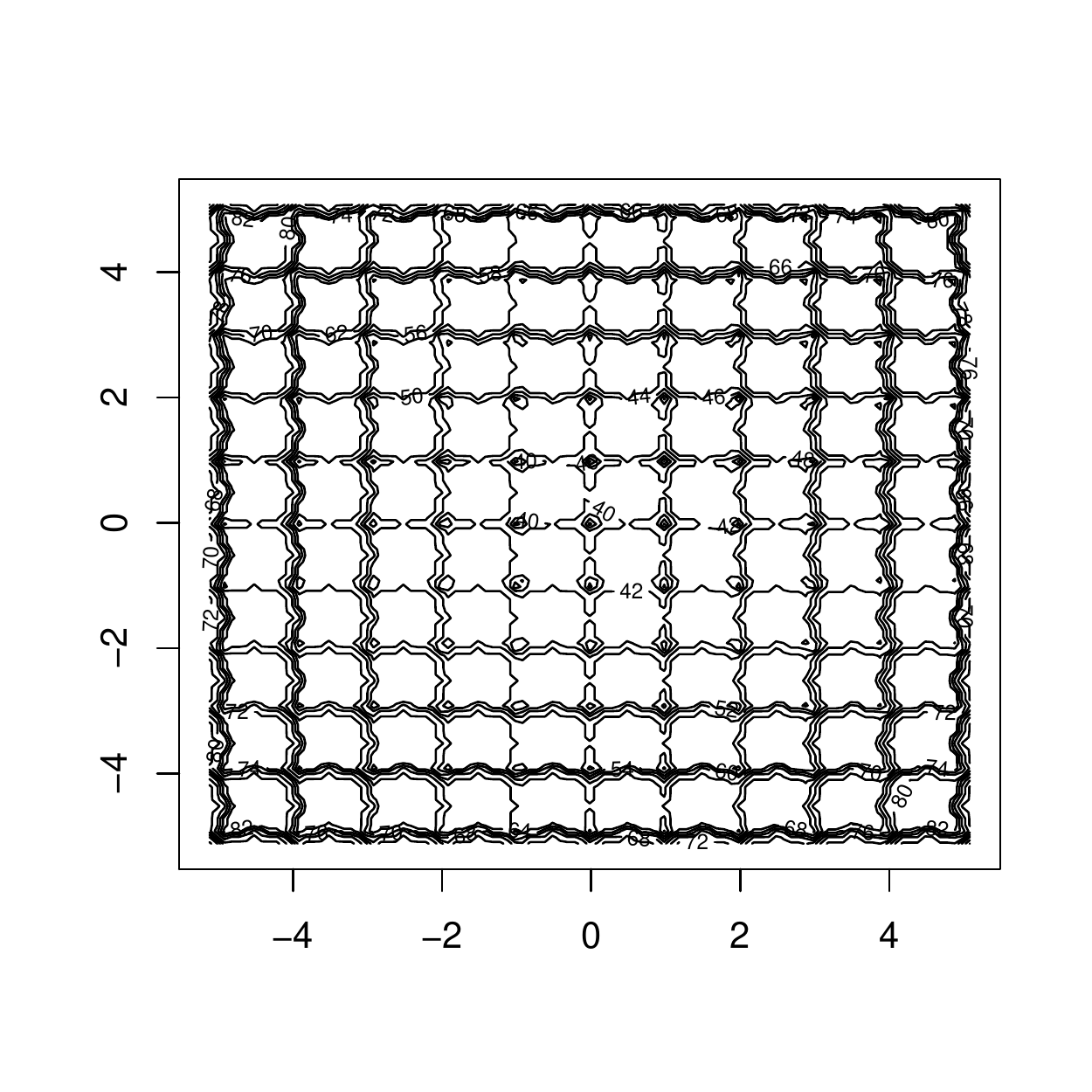}
		\label{fig:RastriginWorstContour}
	\end{subfigure}%
	
	\vspace{-10mm} 
	
	\begin{subfigure}{.20\textwidth}
		\centering
		\includegraphics[width=1.5in, height=1.7in]{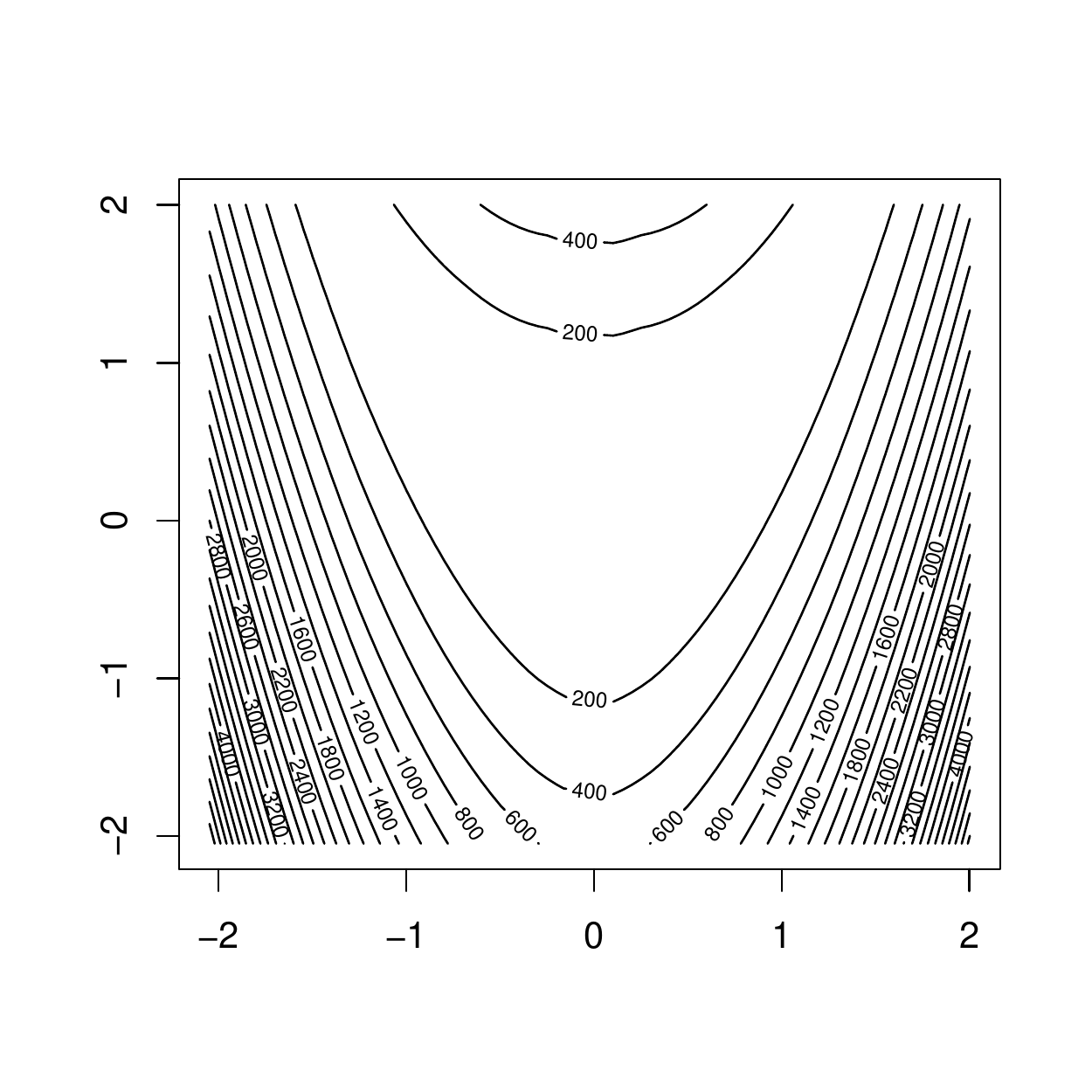}
		\label{fig:RosenbrockWorstContour}
	\end{subfigure}
	\hspace{4mm} 
	\begin{subfigure}{.20\textwidth}
		\centering
		\includegraphics[width=1.5in, height=1.7in]{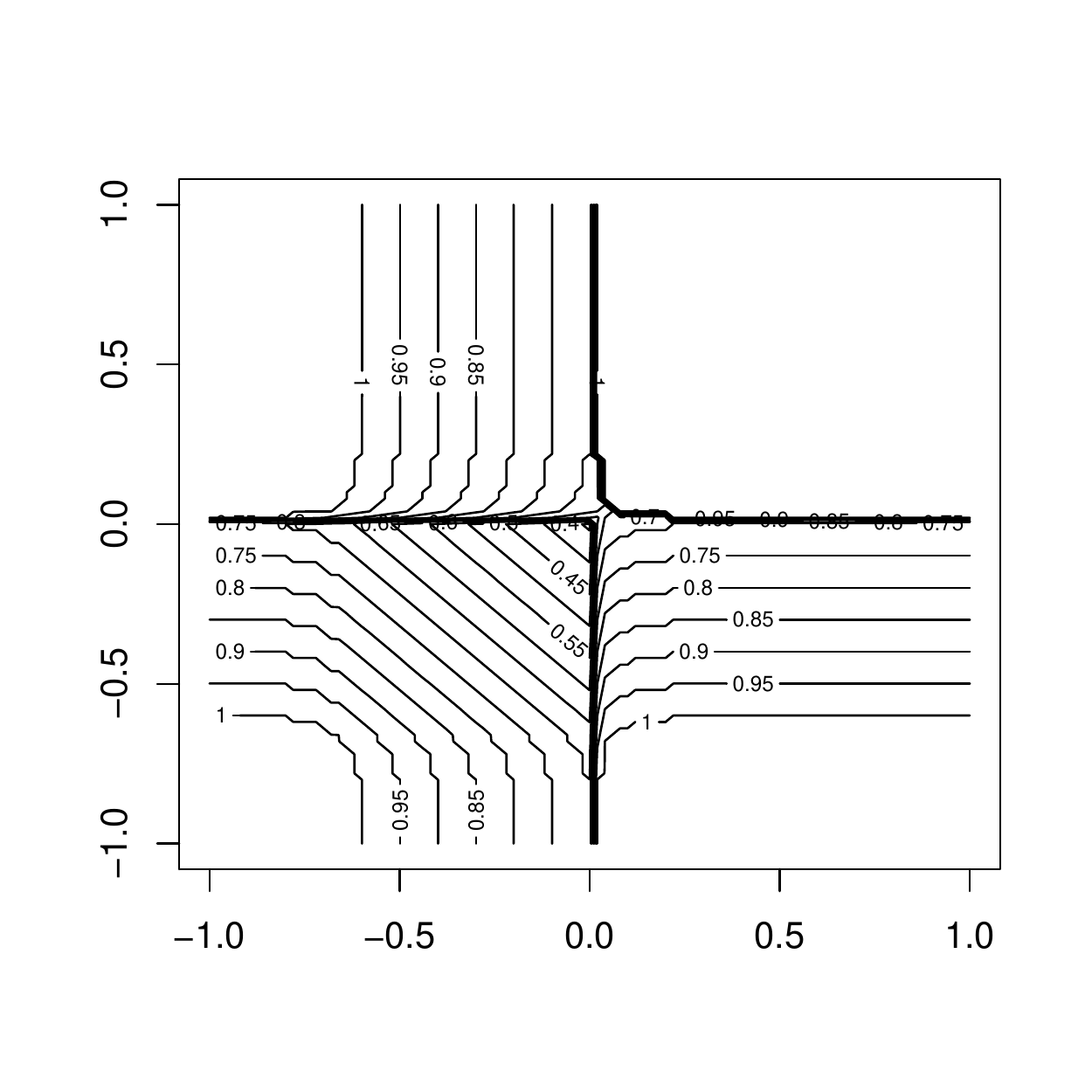}
		\label{fig:SawtoothWorstContour}
	\end{subfigure}
	\hspace{4mm} 
	\begin{subfigure}{.20\textwidth}
		\centering
		\includegraphics[width=1.5in, height=1.7in]{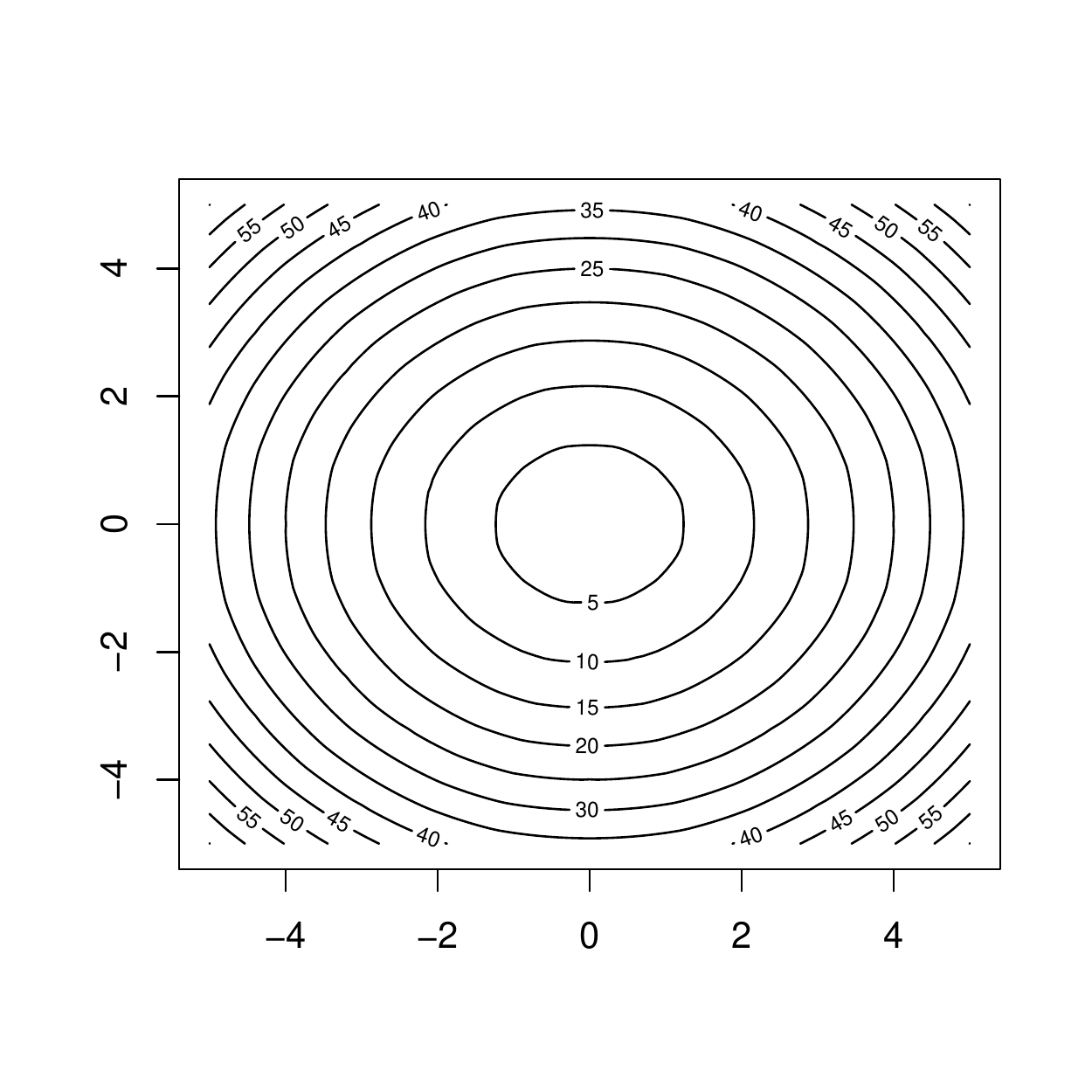}
		\label{fig:SphereWorstContour}
	\end{subfigure}
	\hspace{4mm} 
	\begin{subfigure}{.20\textwidth}
		\centering
		\includegraphics[width=1.5in, height=1.7in]{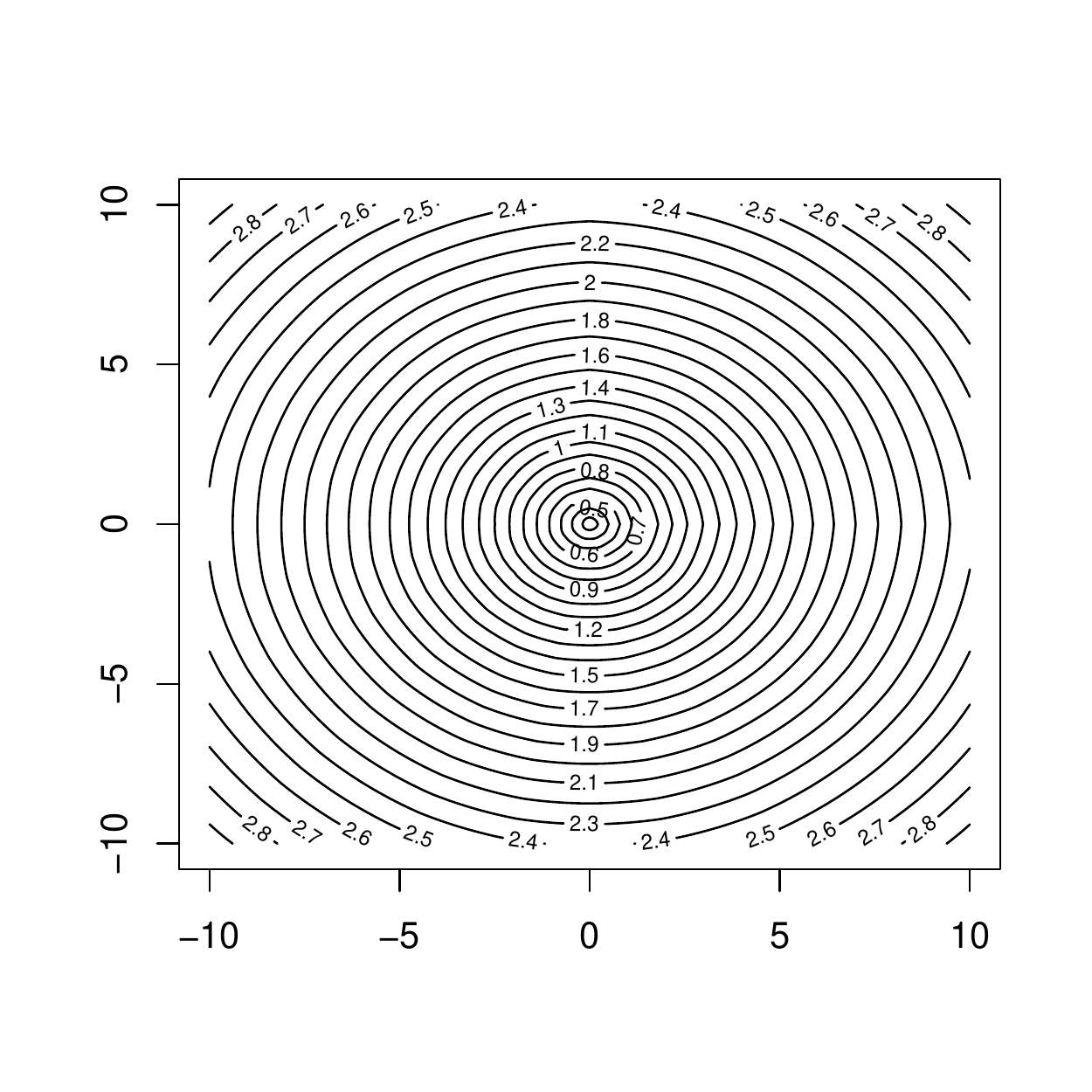}
		\label{fig:VolcanoWorstContour}
	\end{subfigure}%
	
	\caption{Contour plots of nominal (top 8) and worst case (bottom 8) 2D test functions. Left to right, top to bottom: Ackley, Multipeak F1, Multipeak F2, Rastrigin, Sawtooth, Sphere and Volcano.}
	\label{fig:NominalAndWorstCase2DTestFunctions}
\end{figure}

The full description of these eight test functions is given in Appendix~\ref{sec:TestFunctionFormulae}. To give some indication of the nature of these functions contour plots of the 2D instances are shown in Figure~\ref{fig:NominalAndWorstCase2DTestFunctions}, for both the nominal and worst cases.

\subsection{Results}
\label{sec:Results}

Results of the 50 samples runs for each heuristic applied to each test problem-dimension instance are presented here. In each run the best solution as identified by the heuristic is used. However the points in the decision variable space that have been identified as best have robust values generated using the simple inner random sampling approach, with a budget of up to 100 sample points. To better approximate the true robust values at these points their robust values have been re-estimated based on randomly sampling a large number of points (nominally 1,000,000) in the $\Gamma$-uncertainty neighbourhood of the identified robust point. This is a post processing exercise and does not affect the min max search.

\begin{figure}[htbp]
	\centering
			\begin{subfigure}{\textwidth}
				\centering
        \includegraphics[width=.5\textwidth]{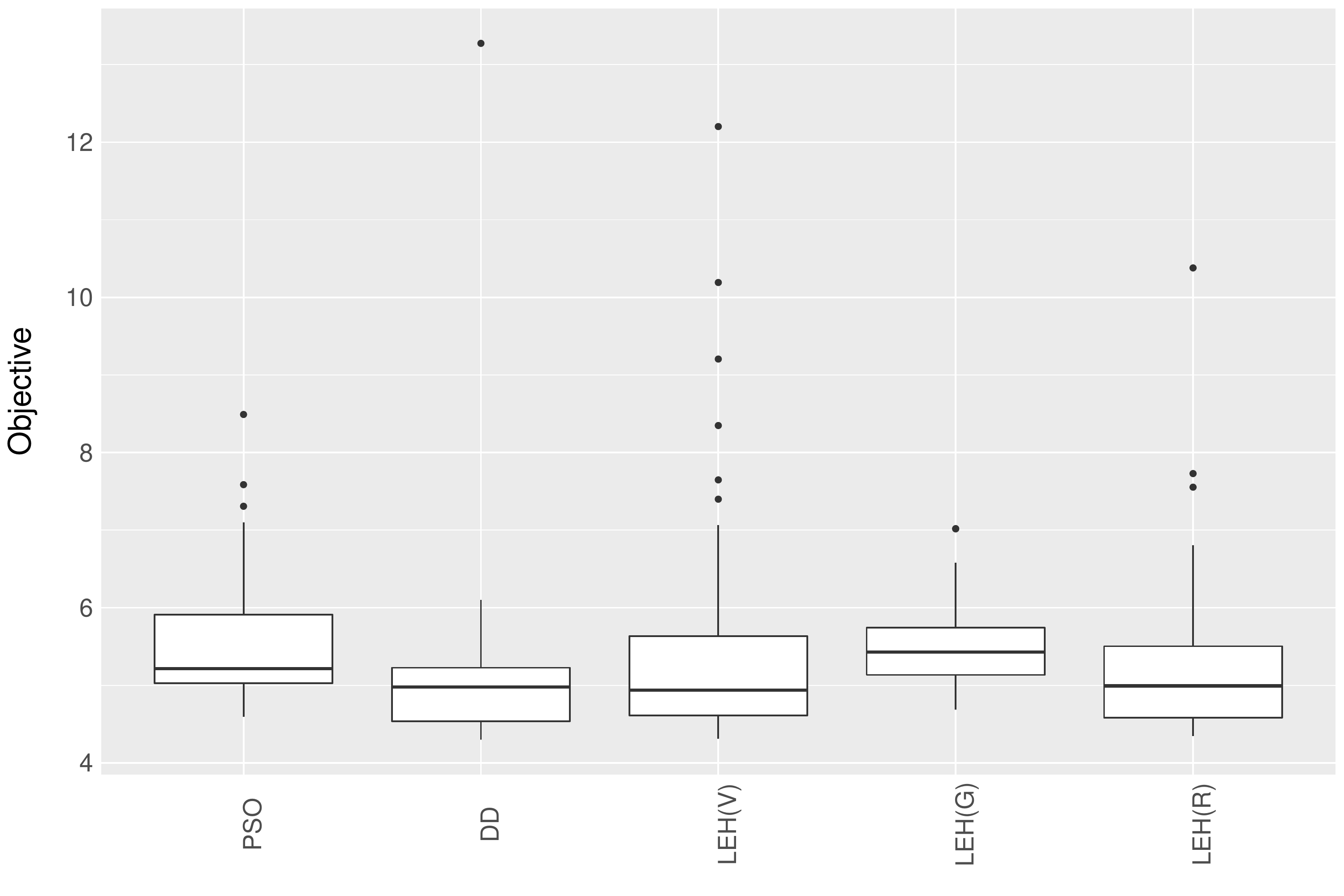}
				\label{fig:Bertsimas2Dboxplot}
				\vspace*{-1.0mm} 
			\end{subfigure} 
	\caption{Box plots of robust objective values due to multiple sample runs for the 2-dimensional polynomial function (poly2D) due to \cite{BertsimasNohadaniTeo2010}.}
	\label{fig:BertBoxPlotResults}
\end{figure}

\begin{figure}[htbp]
	\centering
			\begin{subfigure}{\textwidth}
				\centering
	      \includegraphics[width=1.0\textwidth]{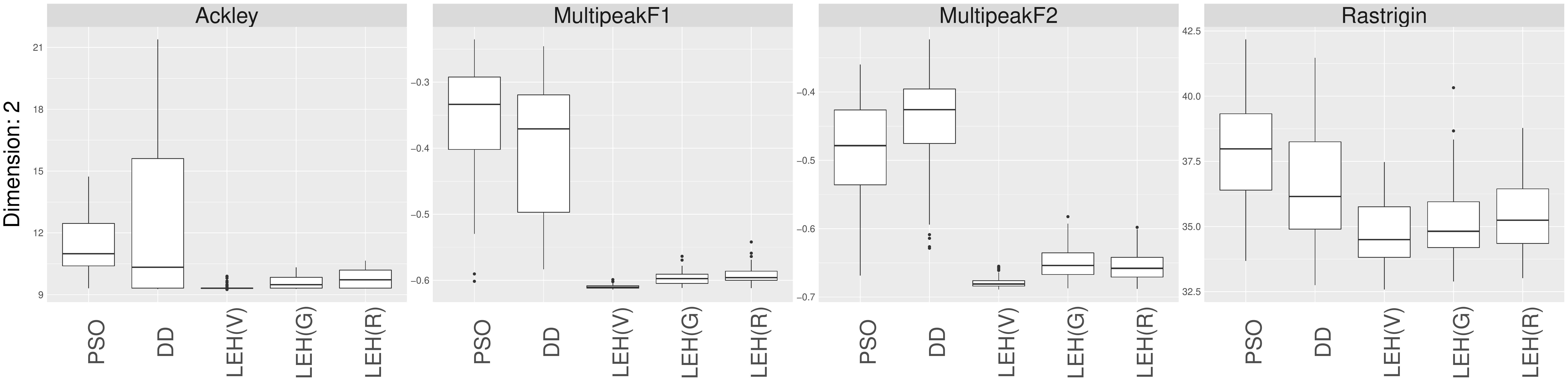}
				\label{fig:Aboxplots2D}
				\vspace*{-1.0mm} 
			\end{subfigure} 

			\begin{subfigure}{\textwidth}
				\centering
        \includegraphics[width=1.0\textwidth]{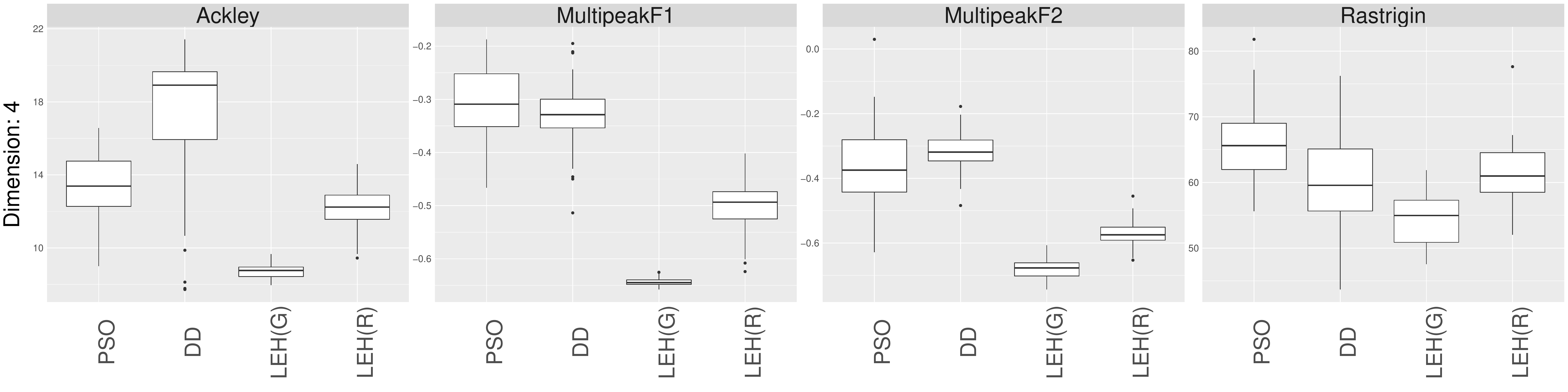}
				\label{fig:Aboxplots4D}
				\vspace*{-1.0mm} 
			\end{subfigure}	

			\begin{subfigure}{\textwidth}
				\centering
        \includegraphics[width=1.0\textwidth]{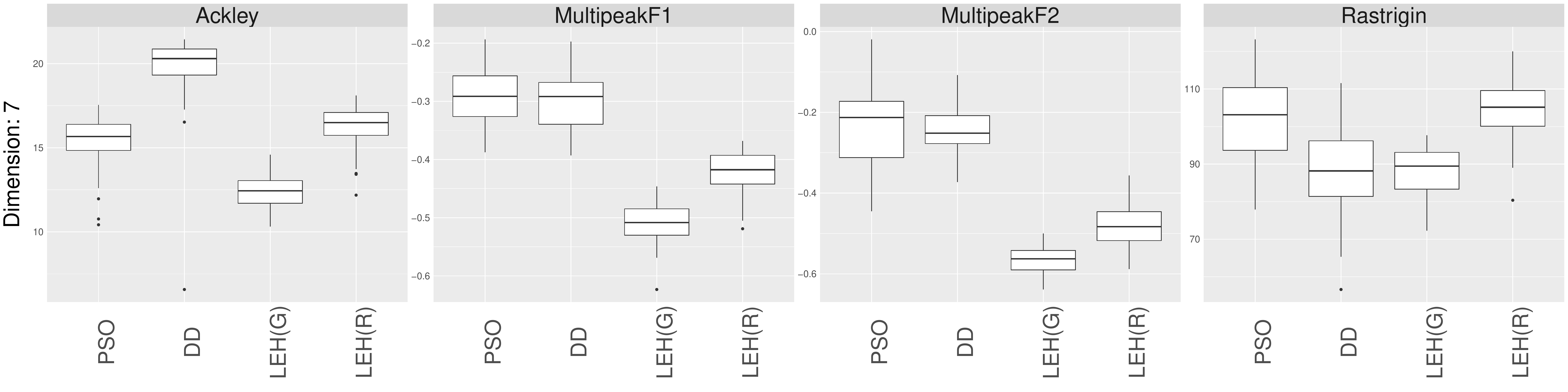}
				\label{fig:Aboxplots7D}
				\vspace*{-1.0mm} 
			\end{subfigure}		

			\begin{subfigure}{\textwidth}
				\centering
        \includegraphics[width=1.0\textwidth]{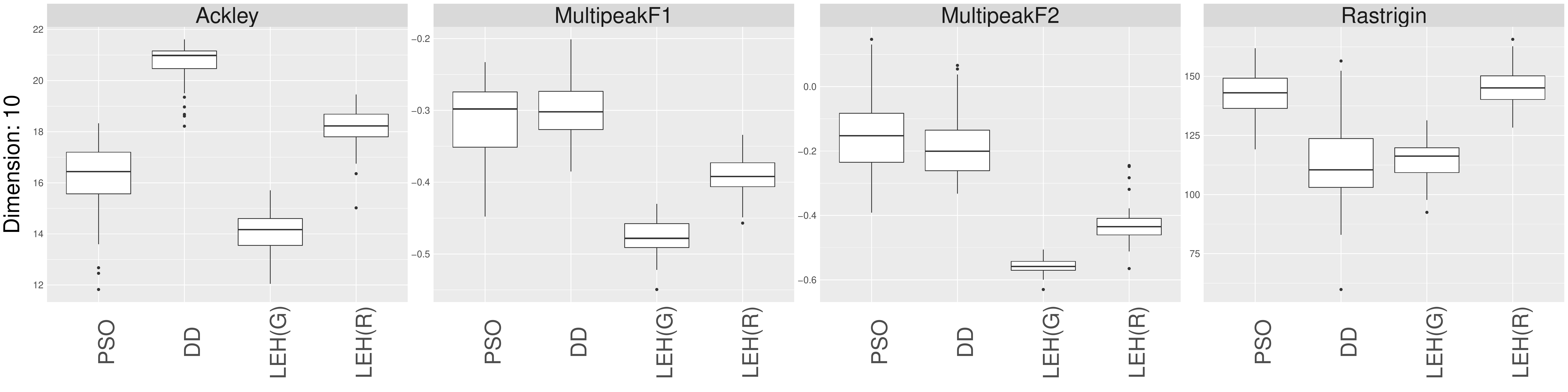}
				\label{fig:Aboxplots10D}
				\vspace*{-1.0mm} 
			\end{subfigure}	

			\begin{subfigure}{\textwidth}
				\centering
        \includegraphics[width=1.0\textwidth]{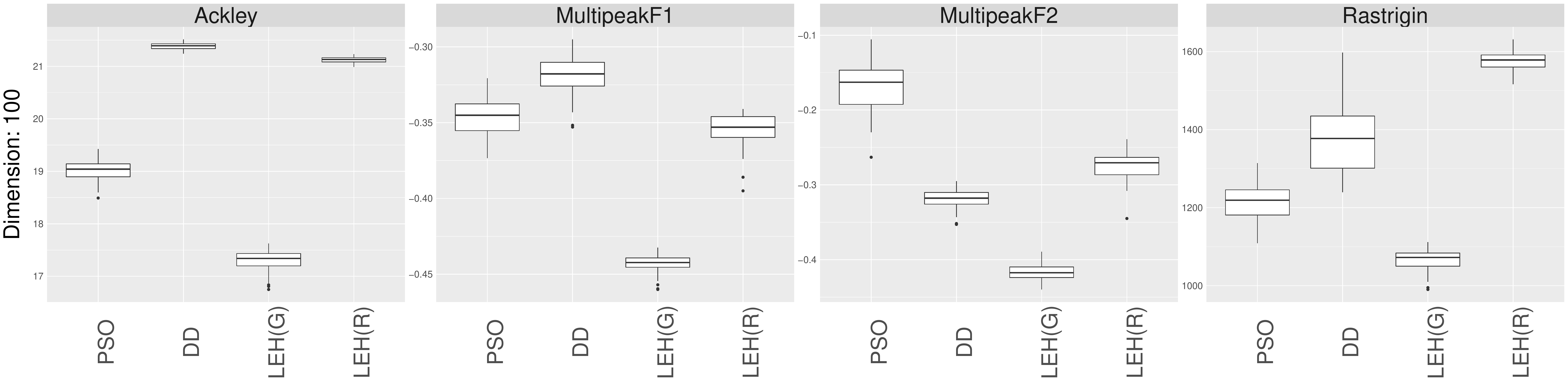}
				\label{fig:Aboxplots100D}
				\vspace*{-1.0mm} 
			\end{subfigure}	

	\caption{Box plots of robust objective values due to multiple sample runs. Left to right: Ackleys, Multipeak F1, Multipeak F2, Rastrigin; Top to bottom: 2D, 4D, 7D, 10D, 100D.}
	\label{fig:BoxPlotResultsa}
\end{figure}

\begin{figure}[htbp]
	\centering
			\begin{subfigure}{\textwidth}
				\centering
        \includegraphics[width=1.0\textwidth]{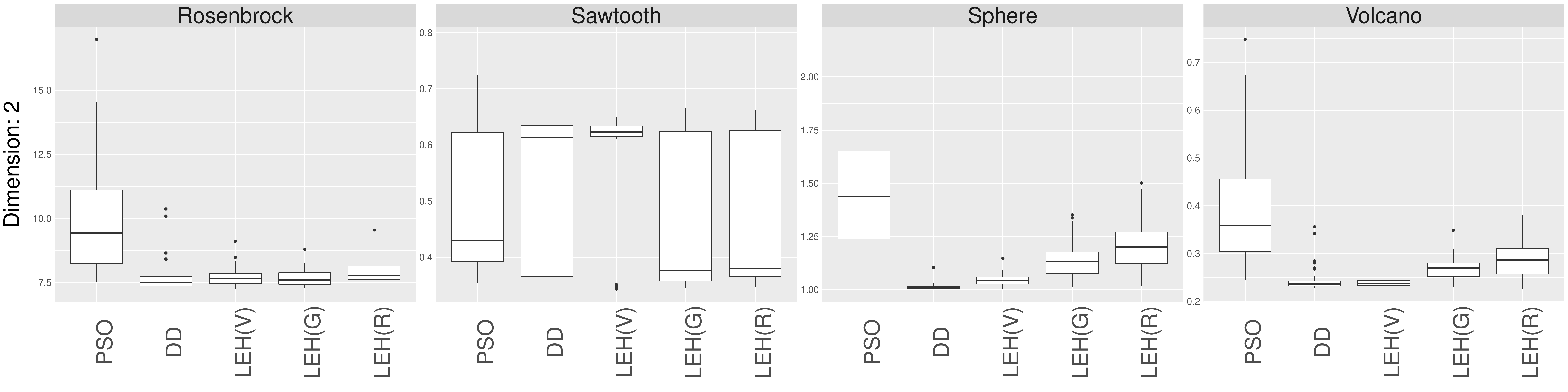}
				\label{fig:Bboxplots2D}
				\vspace*{-1.0mm} 
			\end{subfigure} 

			\begin{subfigure}{\textwidth}
				\centering
        \includegraphics[width=1.0\textwidth]{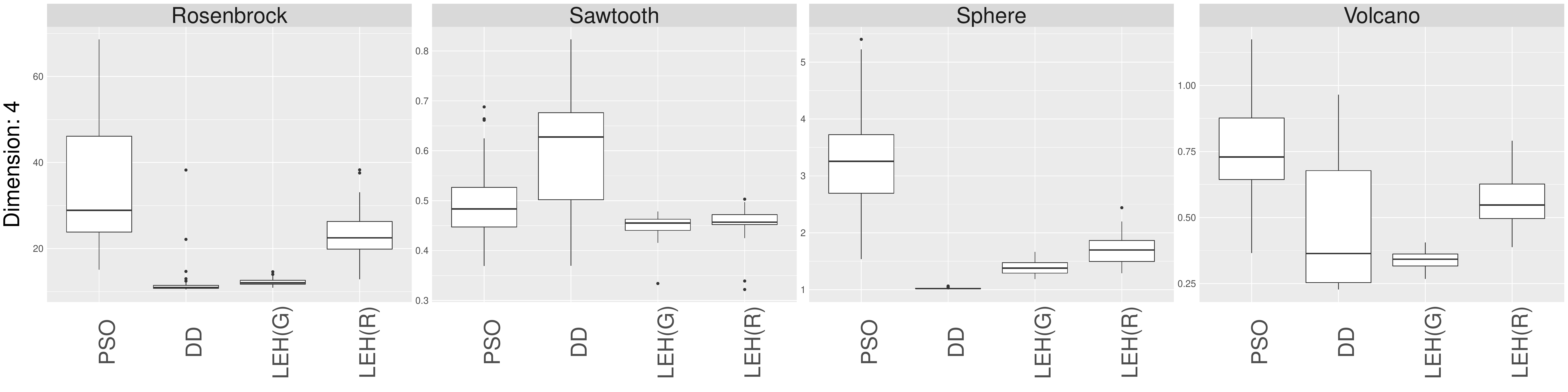}
				\label{fig:Bboxplots4D}
				\vspace*{-1.0mm} 
			\end{subfigure}	

			\begin{subfigure}{\textwidth}
				\centering
        \includegraphics[width=1.0\textwidth]{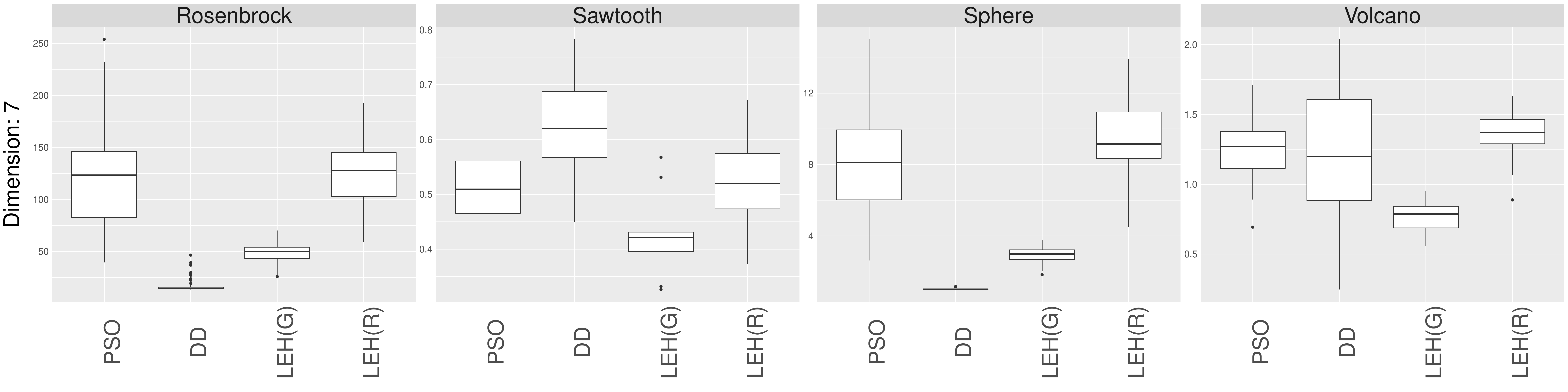}
				\label{fig:Bboxplots7D}
				\vspace*{-1.0mm} 
			\end{subfigure}		

			\begin{subfigure}{\textwidth}
				\centering
	      \includegraphics[width=1.0\textwidth]{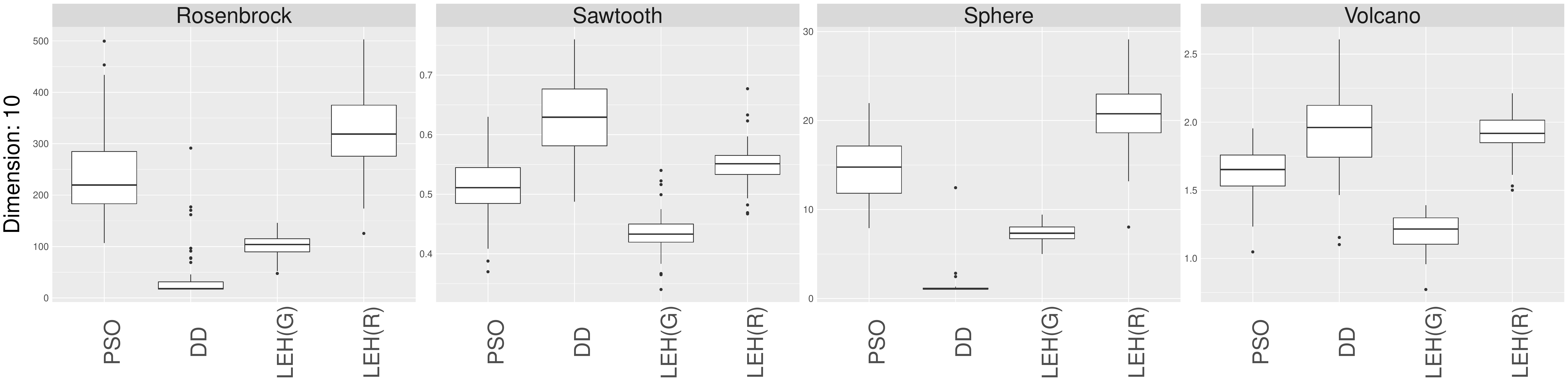}
				\label{fig:Bboxplots10D}
				\vspace*{-1.0mm} 
			\end{subfigure}	

			\begin{subfigure}{\textwidth}
				\centering
        \includegraphics[width=1.0\textwidth]{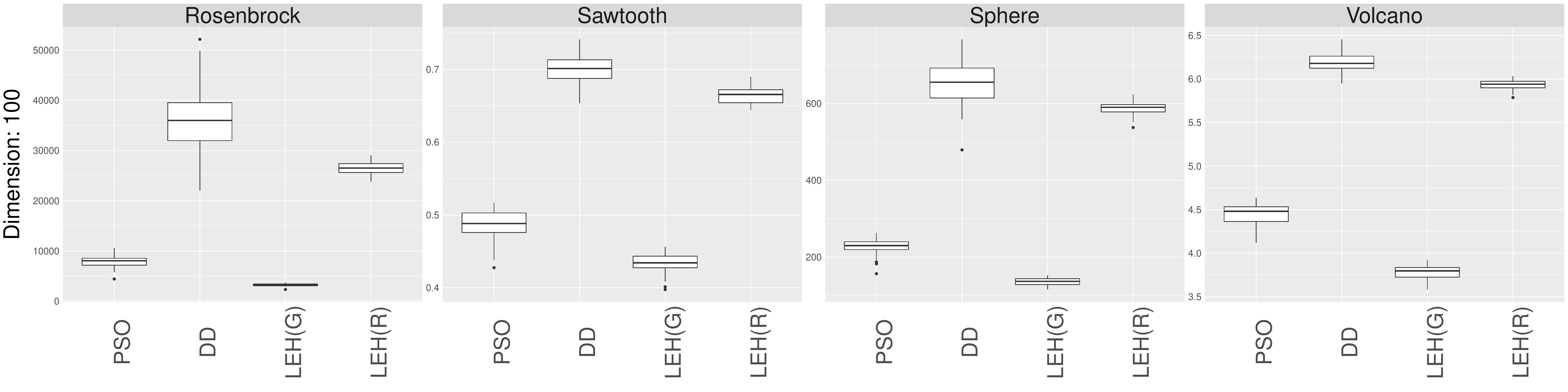}
				\label{fig:Bboxplots100D}
				\vspace*{-1.0mm} 
			\end{subfigure}	

	\caption{Box plots of robust objective values due to multiple sample runs. Left to right: Rosenbrock, Sawtooth, Sphere, Volcano; Top to bottom: 2D, 4D, 7D, 10D, 100D.}
	\label{fig:BoxPlotResultsb}
\end{figure}

Mean results due to each set of 50 sample runs are shown in Tables~\ref{fig:BertMeanResults} and~\ref{fig:MeanResults}. We have applied the Wilcoxon rank-sum test with 95\% confidence to identify the statistically best approaches. Results highlighted in bold indicate the approaches that are statistically equivalent to the best one observed, for a given problem-dimension instance. Corresponding box plots, giving some indication of how the results are distributed across the 50 samples, are shown in Figures~\ref{fig:BertBoxPlotResults},~\ref{fig:BoxPlotResultsa} and~\ref{fig:BoxPlotResultsb}. Additional results, the standard deviations due to each set of 50 sample runs, the average number of candidate points visited and average number of function evaluations undertaken, are given in Appendix~\ref{sec:AdditionalResults}.


\begin{table}[htbp]
\begin{center}
\begin{tabular}{rr|r}
 & & (poly2D) \\
\hline
 \multirow{5}{*}{2D} & PSO  & 5.57 \\
 & \ddre & \textbf{5.11} \\
 & LEH Vor & \textbf{5.52} \\
 & LEH GA  & 5.50 \\
 & LEH Rnd  & \textbf{5.26}
\end{tabular} 
\caption{Mean results due to 50 sample runs for the 2-dimensional polynomial function (poly2D) due to \cite{BertsimasNohadaniTeo2010}.}\label{fig:BertMeanResults}
\end{center}
\end{table}


\begin{sidewaystable}[htbp]
\begin{center}
\begin{tabular}{rr|r|r|r|r|r|r|r|r}
 & & Ackley's & MultipeakF1 & MultipeakF2 & Rastrigin & Rosenbrock & Sawtooth & Sphere & Volcano \\
\hline
 \multirow{5}{*}{2D} & PSO  & 11.44 & -0.36 & -0.49 & 38.04 & 10.00 & \textbf{0.49} & 1.47 & 0.39 \\
 & \ddre & 12.78 & -0.40 & -0.44 & 36.42 & \textbf{7.71} & 0.54 & \textbf{1.01} & \textbf{0.24} \\
 & LEH Vor & \textbf{9.36} & \textbf{-0.61} & \textbf{-0.68} & \textbf{34.67} & \textbf{7.71} & 0.59 & 1.05 & \textbf{0.24} \\
 & LEH GA  & 9.62 & -0.60 & -0.65 & \textbf{35.17} & \textbf{7.68} & \textbf{0.48} & 1.14 & 0.27 \\
 & LEH Rnd  & 9.77 & -0.59 & -0.65 & 35.52 & 7.92 & \textbf{0.47} & 1.21 & 0.29 \\
 \hline
\multirow{4}{*}{4D} & PSO  & 13.50 & -0.30 & -0.36 & 65.91 & 34.20 & 0.50 & 3.35 & 0.75 \\
 & \ddre & 17.32 & -0.33 & -0.32 & 60.43 & \textbf{11.94} & 0.60 & \textbf{1.02} & \textbf{0.46} \\
 & LEH GA  & \textbf{8.73} & \textbf{-0.64} & \textbf{-0.68} & \textbf{54.34} & 12.17 & \textbf{0.45} & 1.39 & \textbf{0.34} \\
 & LEH Rnd  & 12.21 & -0.50 & -0.57 & 61.39 & 23.18 & 0.46 & 1.70 & 0.57 \\
 \hline
\multirow{4}{*}{7D} & PSO  & 15.36 & -0.29 & -0.23 & 102.35 & 123.42 & 0.51 & 8.21 & 1.27 \\
 & \ddre & 19.72 & -0.30 & -0.24 & \textbf{88.44} & \textbf{17.47} & 0.63 & \textbf{1.03} & 1.21 \\
 & LEH GA  & \textbf{12.35} & \textbf{-0.51} & \textbf{-0.57} & \textbf{88.07} & 48.75 & \textbf{0.42} & 2.94 & \textbf{0.77} \\
 & LEH Rnd  & 16.19 & -0.42 & -0.48 & 104.31 & 126.28 & 0.52 & 9.49 & 1.37 \\
 \hline
\multirow{4}{*}{10D} & PSO  & 16.17 & -0.31 & -0.15 & 142.99 & 238.36 & 0.51 & 14.66 & 1.63 \\
 & \ddre & 20.69 & -0.30 & -0.19 & \textbf{112.61} & \textbf{41.12} & 0.63 & \textbf{1.40} & 1.93 \\
 & LEH GA  & \textbf{14.08} & \textbf{-0.48} & \textbf{-0.56} & \textbf{115.06} & 103.31 & \textbf{0.43} & 7.34 & \textbf{1.19} \\
 & LEH Rnd  & 18.11 & -0.39 & -0.43 & 145.52 & 322.27 & 0.55 & 20.62 & 1.92 \\
 \hline
\multirow{4}{*}{100D} & PSO  & 19.02 & -0.35 & -0.17 & 1,215.34 & 7,989.77 & 0.49 & 226.66 & 4.45 \\
 & \ddre & 21.38 & -0.32 & -0.32 & 1,386.77 & 36,141.80 & 0.70 & 656.86 & 6.18 \\
 & LEH GA  & \textbf{17.30} & \textbf{-0.44} & \textbf{-0.42} & \textbf{1,065.44} & \textbf{3,264.49} & \textbf{0.43} & \textbf{136.18} & \textbf{3.79} \\
 & LEH Rnd  & 21.12 & -0.36 & -0.28 & 1,577.84 & 26,526.42 & 0.66 & 588.03 & 5.93
\end{tabular}
\caption{Mean results due to 50 sample runs.}\label{fig:MeanResults}
\end{center}
\end{sidewaystable}



From Table~\ref{fig:MeanResults} we see that for 100D instances the LEH GA approach is best for all test problems, and in several cases the mean LEH GA result is substantially better than all of the alternative heuristics. From Tables~\ref{fig:BertMeanResults} and~\ref{fig:MeanResults} the LEH approach is among the best in at least 6 of the instances for all other dimensions.

For 2D instances the LEH Voronoi approach is among the best results for 7 of the 9 problems, whilst LEH GA and LEH Rnd are each amongst the best results for 3 and 2 problems respectively. It should also be noted that in 5 of the 7 instances where LEH Voronoi is among the best, LEH GA is either statistically equivalent or the mean value is second best. For the 2D Sphere instance {\ddre} is marginally better than LEH Voronoi, whilst {\ddre} and LEH heuristics are statistically equivalent for the (poly2D) and 2D Volcano and Rosenbrock instances. The robust PSO approach is statistically equivalent to LEH heuristics for the 2D Sawtooth instance.

For the 4D -- 10D instances {\ddre} is statistically equivalent to LEH GA in the 4D Volcano problem and the 7D and 10D instances of the Rastrigin problem, and better than LEH GA for the Rosenbrock and Sphere problems. For the 4D Rosenbrock and Sphere problems the differences between {\ddre} and LEH GA are reasonably small, however in the 7D and 10D instances {\ddre} is substantially better. Considering the shape of the Rosenbrock and Sphere functions it can be expected that a local search will perform particularly well for these problems.

LEH GA is better than LEH Rnd for all instances excluding (poly2D). In a number of instances LEH GA is substantially better than LEH Rnd. The number of candidate points that LEH can visit is substantially increased by the early stopping of inner searches as soon as the high cost threshold is exceeded, see Tables~\ref{fig:BertOtherResults} and~\ref{fig:CandidateResults} in Appendix~\ref{sec:AdditionalResults}. Although this feature must unquestionably play a role in the success of the LEH GA approach, the fact that LEH Rnd visits a comparable number of candidate points indicates that the additional pro active seeking of the largest hypersphere devoid of high cost points is also a significant factor in the success of LEH GA.


\section{Conclusions and further work}
\label{sec:SummaryConcusionsFurtherWork}

We have introduced a new metaheuristic for box-constrained robust optimisation problems with implementation uncertainty. We do not assume any knowledge on the structure of the original objective function, making the approach applicable to black-box and simulation-optimisation problems. We do assume that the solution is affected by uncertainty, and the aim is to find a solution that optimises the worst possible performance in this setting. This is the min max problem. Previously, few generic search methods have been developed for this setting.

We introduce a new approach for a global search based on distinguishing undesirable high cost -- high objective value -- points (hcps), identifying the largest hypersphere in the decision variable space that is completely devoid of hcps, and exploring the decision variable space by stepping between the centres of these largest empty hyperspheres.

We demonstrated the effectiveness of the approach using a series of test problems, considering instances of varying dimension, and comparing our LEH approach against one metaheuristic that employs an outer particle swarm optimisation and one from the literature that uses multiple re-starts of the local descent directions approach. For low and moderate dimensional instances the approach shows competitive performance; for high-dimensional problems the LEH approach significantly outperforms the comparator heuristics  for all problems.

There are several ways in which this work can be developed. Further consideration can be given to the inner maximisation search approach in order to better understand the trade-off between expending function evaluations on the local $\Gamma$-radius uncertainty neighbourhood search versus globally exploring the search space, in the context of our LEH approach.

The repeated calculation of large numbers of Euclidean distances each time a new LEH needs to be identified within the LEH GA heuristic is computationally expensive. Rather than only calculating a single next candidate point each time the GA is performed, identifying multiple points could speed up computation or alternatively enable the use of larger population-generation sizes to improve the estimation of the largest empty hypersphere.

Results of the mid-dimension experiments on the Rosenbrock and Sphere test problems suggest that an exploitation based approach works well in these instances, indicating a direction for extending our exploration focussed LEH approach.

It is clear that within the LEH algorithm the early stopping of the inner searches when it is established that the current robust global value cannot be improved upon has significant advantages. It is worth considering whether alternative search approaches could take advantage of this feature.


\appendix

\section*{Appendices}

\section{Test functions}
\label{sec:TestFunctionFormulae}

Functions used to assess the effectiveness of the Largest Empty Hypersphere robust metaheuristics taken from \cite{Kruisselbrink2012, JamilYang2013}.
\\
 
\noindent \textbf{Ackleys}
 
$$f(\pmb{x}) = -20 \exp \Bigg( -0.2 \sqrt{\frac{1}{n} \sum_{i=1}^n x_i^2} \Bigg) - \exp \Bigg( \frac{1}{n} \sum_{i=1}^n \cos(2\pi x_i) \Bigg) + 20 + \exp(1)$$

\vspace{2mm} 
\noindent The feasible region is the hypercube $x_{i} \in$ [-32.768, 32.768].
\\

\vspace{3mm} 
\noindent \textbf{MultipeakF1}

$$f(\pmb{x}) = - \frac{1}{n} \sum_{i=1}^n g(x_i) \ , \ g(x_i) =
\begin{cases}
	\exp(2 \ln 2 ( \frac{x_i-0.1}{0.8} )^2) \sqrt{ \left| \sin(5\pi x_i) \right| }  &\quad\text{if } 0.4 < x_i \le 0.6 \text{ ,} \\
  \exp(2 \ln 2 ( \frac{x_i-0.1}{0.8} )^2) \sin^{6}(5\pi x_i)  &\quad\text{otherwise} \\
\end{cases}$$

\vspace{2mm}
\noindent The feasible region is the hypercube $x_{i} \in$ [0, 1].
\\

\vspace{3mm}
\noindent \textbf{MultipeakF2}

$$f(\pmb{x}) = \frac{1}{n} \sum_{i=1}^n g(x_i) \ , \ g(x_i) =2\sin(10 \exp (-0.2x_i) x_i) \exp (-0.25x_i)$$

\vspace{2mm}
\noindent The feasible region is the hypercube $x_{i} \in$ [0, 10].
\\

\vspace{3mm}
\noindent \textbf{Rastrigin}
 
$$f(\pmb{x}) = 10 n + \sum_{i=1}^{n}{[x_{i}^{2} - 10 \cos(2\pi x_i)]}$$

\vspace{2mm} 
\noindent The feasible region is the hypercube $x_{i} \in$ [-5.12, 5.12].
\\

\vspace{3mm} 
\noindent \textbf{Rosenbrock} 

$$f(\pmb{x}) = \sum_{i=1}^{n-1} [100(x_{i+1} - x_i^2)^2 + (x_i-1)^2]$$

\vspace{2mm} 
\noindent The feasible region is the hypercube $x_{i} \in$ [-2.048, 2.048].
\\

\vspace{3mm} 
\noindent \textbf{Sawtooth} 

$$f(\pmb{x}) = 1 - \frac{1}{n} \sum_{i=1}^n g(x_i) \ , \ g(x_i) =
\begin{cases}
	x_i + 0.8  &\quad\text{if } -0.8 \le x_i < 0.2 \text{ ,} \\
  0  &\quad\text{otherwise} \\
\end{cases}$$

\vspace{2mm} 
\noindent The feasible region is the hypercube $x_{i} \in$ [-1, 1].
\\

\vspace{3mm}
\noindent \textbf{Sphere} 

$$f(\pmb{x}) = \sum_{i=1}^n x_i^2$$

\vspace{2mm} 
\noindent The feasible region is the hypercube $x_{i} \in$ [-5, 5].
\\

\vspace{3mm}
\noindent \textbf{Volcano} 

$$f(\pmb{x}) =
\begin{cases}
	\sqrt{\lVert \pmb{x} \rVert} - 1  &\quad\text{if } \lVert \pmb{x} \rVert > 1 \text{ ,} \\
  0  &\quad\text{otherwise} \\
\end{cases}$$

\vspace{2mm} 
\noindent The feasible region is the hypercube $x_{i} \in$ [-10, 10].

\section{Radii due to alternative LEH algorithms}
\label{sec:RadiiLEH}

\begin{figure}[htbp]
	\centering
	\includegraphics[width=1.0\textwidth]{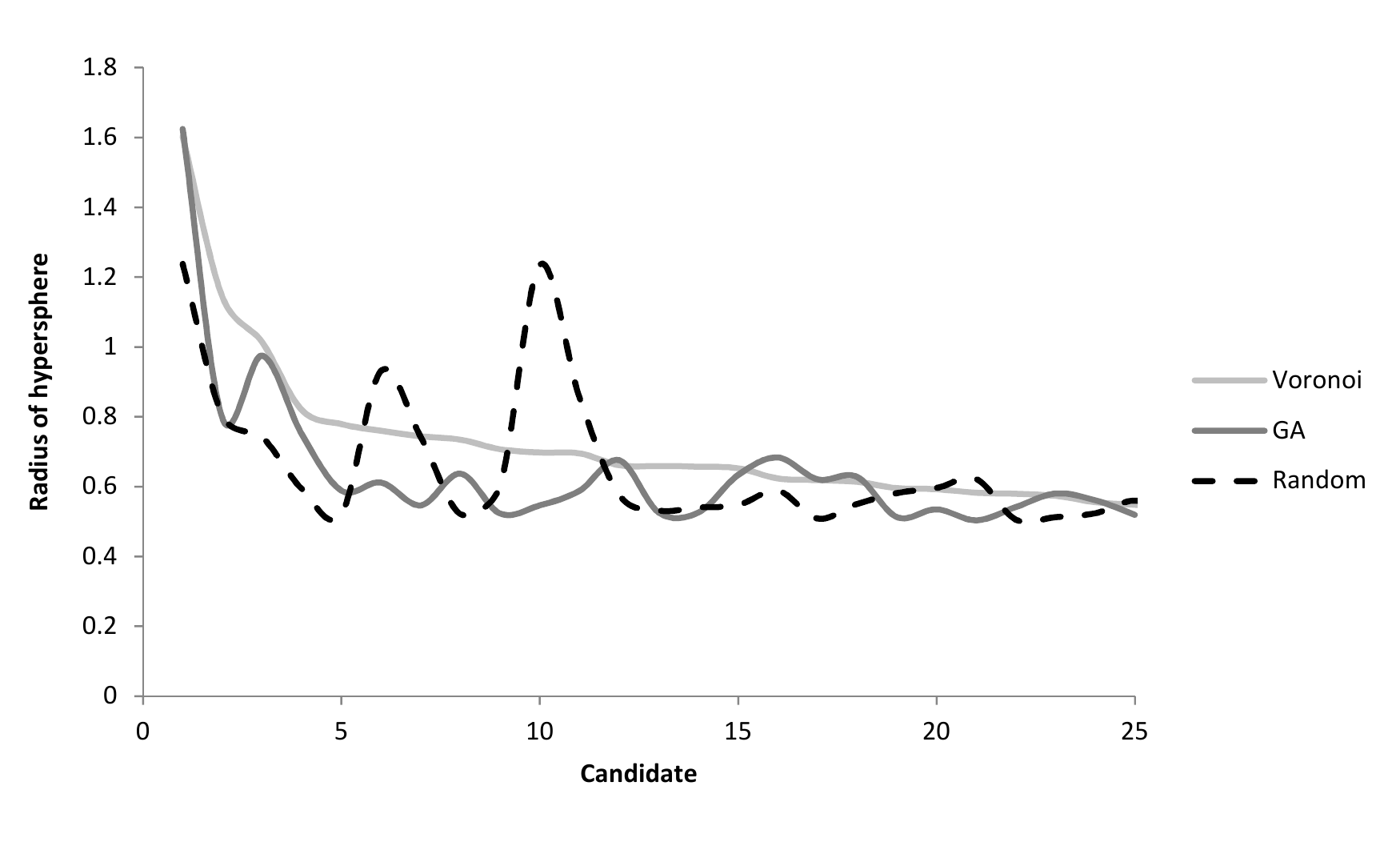} 
	\caption{Alternative LEH approaches applied to the same problem: variation in empty hypersphere radius with numbers of candidates evaluated for robustness.}
	\label{fig:AlternativeLEHSearchesRadii}
\end{figure}

Whilst the Voronoi based search exemplified by Figures~\ref{fig:LEHVoronoiPoints} and~\ref{fig:LEHVoronoiSearch} in Section~\ref{sec:NatureLEHSearch} is a good indicator of the nature of the searches due to all three alternative LEH approaches, random, GA and Voronoi, the radii of the LEH identified for each candidate will vary across these approaches. Here Figure~\ref{fig:AlternativeLEHSearchesRadii} gives some indication of how the radii of the hyperspheres generated by each these three LEH heuristics progress as the exploration proceeds. The three curves represent separate runs of the LEH algorithm when applied to (poly2D),and should be considered indicative.

As would be expected the general nature of the size of the radius of the LEH steadily decreases with increasing numbers of candidate points evaluated. However superimposed on this overall decrease are the indicative patterns due to the alternative heuristics. For the random algorithm the size of the LEH is quite variable, whilst for Voronoi the curve is smooth. The GA algorithm sits somewhere between the two.

\section{Additional results}
\label{sec:AdditionalResults}

The standard deviations due to each set of 50 sample runs, the average number of candidate points visited and average number of function evaluations undertaken, are shown in Tables~\ref{fig:BertOtherResults},~\ref{fig:SDResults},~\ref{fig:CandidateResults} and~\ref{fig:EvaluationResults} here respectively. Labelling of comparator heuristics in the tables is as follows:

\begin{itemize}[leftmargin=*]  
	\item PSO: particle swarm optimisation.
	\item {\ddre}: Multi re-start descent directions.	
	\item LEH Vor: LEH using a Voronoi \cite{Toussaint1983} approach; applied to 2D problems only. 
	\item LEH GA: LEH using a  genetic algorithm.
	\item LEH Rnd: LEH using random sampling.
\end{itemize} 


\begin{table}[htbp]
\begin{center}
\begin{tabular}{rr|r|r|r}
 & & Std. dev. & Candidates & Evaluations  \\
\hline
 \multirow{5}{*}{poly2D} & PSO & 0.83 & 100 &10,000 \\
 & \ddre & 1.27 & 100 & 10,000 \\
 & LEH Vor & 1.60 & 35 & 995 \\
 & LEH GA  & 0.56 & 30 & 727 \\
 & LEH Rnd  & 1.07 & 35 & 1,037
\end{tabular}
\caption{Standard deviations of results, average number of candidate points visited, and average number of points evaluated for the 50 sample runs for the 2-dimensional polynomial function (poly2D) due to \cite{BertsimasNohadaniTeo2010}.}\label{fig:BertOtherResults}
\end{center}
\end{table}


\begin{sidewaystable}[htbp]
\begin{center}
\begin{tabular}{rr|r|r|r|r|r|r|r|r}
 & & Ackley's & MultipeakF1 & MultipeakF2 & Rastrigin & Rosenbrock & Sawtooth & Sphere & Volcano \\
\hline
 \multirow{5}{*}{2D} & PSO  & 1.35 & 0.09 & 0.08 & 2.04 & 2.17 & 0.12 & 0.30 & 0.11 \\
& {\ddre} & 4.14 & 0.10 & 0.08 & 2.26 & 0.62 & 0.14 & 0.02 & 0.02 \\
& LEH Vor & 0.15 & 0.00 & 0.01 & 1.23 & 0.34 & 0.10 & 0.03 & 0.01 \\
& LEH GA  & 0.35 & 0.01 & 0.02 & 1.53 & 0.31 & 0.14 & 0.09 & 0.03 \\
& LEH Rnd  & 0.45 & 0.01 & 0.02 & 1.43 & 0.47 & 0.13 & 0.11 & 0.04 \\
\hline
 \multirow{4}{*}{4D} & PSO  & 1.66 & 0.06 & 0.13 & 5.79 & 14.21 & 0.08 & 0.89 & 0.18 \\
& {\ddre} & 3.55 & 0.06 & 0.06 & 7.29 & 4.17 & 0.11 & 0.01 & 0.25 \\
& LEH GA  & 0.38 & 0.01 & 0.03 & 4.04 & 0.84 & 0.02 & 0.13 & 0.03 \\
& LEH Rnd  & 1.19 & 0.05 & 0.04 & 4.47 & 5.30 & 0.03 & 0.26 & 0.08 \\
\hline
 \multirow{4}{*}{7D} & PSO  & 1.54 & 0.05 & 0.10 & 10.70 & 49.51 & 0.06 & 2.66 & 0.21 \\
& {\ddre} & 2.21 & 0.05 & 0.07 & 11.08 & 7.19 & 0.08 & 0.03 & 0.47 \\
& LEH GA  & 1.06 & 0.03 & 0.03 & 6.58 & 9.90 & 0.04 & 0.44 & 0.10 \\
& LEH Rnd  & 1.26 & 0.04 & 0.05 & 8.20 & 31.15 & 0.07 & 2.05 & 0.13 \\
\hline
 \multirow{4}{*}{10D} & PSO  & 1.43 & 0.05 & 0.13 & 10.06 & 89.74 & 0.06 & 3.36 & 0.19 \\
& {\ddre} & 0.83 & 0.04 & 0.09 & 16.48 & 53.56 & 0.06 & 1.63 & 0.29 \\
& LEH GA  & 0.82 & 0.03 & 0.02 & 7.95 & 20.73 & 0.04 & 1.08 & 0.13 \\
& LEH Rnd  & 0.79 & 0.03 & 0.06 & 7.93 & 80.92 & 0.04 & 3.77 & 0.14 \\
\hline
 \multirow{4}{*}{100D} & PSO  & 0.20 & 0.01 & 0.03 & 44.69 & 1284.24 & 0.02 & 20.84 & 0.13 \\
& {\ddre} & 0.06 & 0.01 & 0.01 & 97.02 & 6518.29 & 0.02 & 59.59 & 0.11 \\
& LEH GA  & 0.20 & 0.01 & 0.01 & 27.87 & 251.68 & 0.01 & 9.15 & 0.08 \\
& LEH Rnd  & 0.06 & 0.01 & 0.02 & 26.49 & 1289.96 & 0.01 & 16.83 & 0.06
\end{tabular} 
\caption{Standard deviations of results due to 50 sample runs.}\label{fig:SDResults}
\end{center}
\end{sidewaystable}


\begin{sidewaystable}[htbp]
\begin{center}
\begin{tabular}{rr|r|r|r|r|r|r|r|r}
 & & Ackley's & MultipeakF1 & MultipeakF2 & Rastrigin & Rosenbrock & Sawtooth & Sphere & Volcano \\
\hline
\multirow{5}{*}{2D} & PSO  & 100 & 100 & 100 & 100 & 100 & 100 & 100 & 100 \\
 & {\ddre} & 100 & 100 & 100 & 100 & 100 & 100 & 100 & 100 \\
 & LEH Vor & 245 & 107 & 205 & 309 & 129 & 26 & 24 & 63 \\
 & LEH GA  & 147 & 78 & 117 & 160 & 90 & 36 & 33 & 57 \\
 & LEH Rnd  & 213 & 98 & 183 & 256 & 109 & 26 & 23 & 53 \\
\hline
\multirow{4}{*}{4D} & PSO  & 100 & 100 & 100 & 100 & 100 & 100 & 100 & 100 \\
 & {\ddre} & 100 & 100 & 100 & 100 & 100 & 100 & 100 & 100 \\
 & LEH GA  & 9,054 & 7,543 & 8,355 & 3,624 & 8,433 & 1,304 & 1,334 & 3,890 \\
 & LEH Rnd  & 8,799 & 7,072 & 7,770 & 5,358 & 8,644 & 6,264 & 6,408 & 8,646 \\
\hline
\multirow{4}{*}{7D} & PSO  & 100 & 100 & 100 & 100 & 100 & 100 & 100 & 100 \\
 & {\ddre} & 100 & 100 & 100 & 100 & 100 & 100 & 100 & 100 \\
 & LEH GA  & 9,097 & 7,958 & 8,754 & 5,522 & 9,062 & 7,772 & 9,164 & 9,184 \\
 & LEH Rnd  & 8,799 & 7,341 & 7,819 & 6,821 & 8,781 & 8,229 & 8,789 & 8,816 \\
\hline
\multirow{4}{*}{10D} & PSO  & 100 & 100 & 100 & 100 & 100 & 100 & 100 & 100 \\
 & {\ddre} & 100 & 100 & 100 & 100 & 100 & 100 & 100 & 100 \\
 & LEH GA  & 9,063 & 7,822 & 8,716 & 6,526 & 9,039 & 7,977 & 9,037 & 9,086 \\
 & LEH Rnd  & 8,746 & 7,674 & 7,679 & 7,367 & 8,874 & 8,279 & 8,734 & 8,773 \\
\hline
\multirow{4}{*}{100D} & PSO  & 100 & 100 & 100 & 100 & 100 & 100 & 100 & 100 \\
 & {\ddre} & 100 & 100 & 100 & 100 & 100 & 100 & 100 & 100 \\
 & LEH GA  & 8,902 & 8,779 & 8,765 & 8,615 & 8,848 & 8,551 & 8,907 & 8,877 \\
 & LEH Rnd  & 7,805 & 8,679 & 8,533 & 8,708 & 8,954 & 8,708 & 8,889 & 8,827
\end{tabular} 
\caption{Average number of candidate points visited across 50 sample runs.}\label{fig:CandidateResults}
\end{center}
\end{sidewaystable}


\begin{sidewaystable}[htbp]
\begin{center}
\begin{tabular}{rr|r|r|r|r|r|r|r|r}
 & & Ackley's & MultipeakF1 & MultipeakF2 & Rastrigin & Rosenbrock & Sawtooth & Sphere & Volcano \\
\hline
\multirow{5}{*}{2D} & PSO  & 10,000 & 10,000 & 10,000 & 10,000 & 10,000 & 10,000 & 10,000 & 10,000 \\
 & {\ddre} & 10,000 & 10,000 & 10,000 & 10,000 & 10,000 & 10,000 & 10,000 & 10,000 \\
 & LEH Vor & 945 & 996 & 1,164 & 2,655 & 1,102 & 807 & 546 & 577 \\
 & LEH GA  & 830 & 890 & 958 & 1,835 & 871 & 780 & 686 & 696 \\
 & LEH Rnd  & 1,255 & 1,250 & 1,306 & 2,521 & 1,098 & 1,110 & 934 & 1,021 \\
\hline
\multirow{4}{*}{4D} &PSO  & 10,000 & 10,000 & 10,000 & 10,000 & 10,000 & 10,000 & 10,000 & 10,000 \\
 & {\ddre} & 10,000 & 10,000 & 10,000 & 10,000 & 10,000 & 10,000 & 10,000 & 10,000 \\
 & LEH GA  & 10,000 & 9,198 & 10,000 & 10,000 & 9,822 & 2,899 & 2,279 & 4,866 \\
 & LEH Rnd  & 10,000 & 10,000 & 10,000 & 10,000 & 10,000 & 8,477 & 8,031 & 10,000
\end{tabular}
\caption{Average number of points evaluated in the 50 sample runs. For all test problems of dimension 7 or higher the full budget of 10,000 function evaluations was used in all sample runs.}\label{fig:EvaluationResults}
\end{center}
\end{sidewaystable}

\clearpage 
\newpage



\begin{thebibliography}{VMdCM11}

\bibitem[ABV09]{AissiBazganVanderpooten2009}
H.~Aissi, C.~Bazgan, and D.~Vanderpooten.
\newblock Min–max and min–max regret versions of combinatorial optimization
  problems: A survey.
\newblock {\em European Journal of Operational Research}, 197(2):427 -- 438,
  2009.

\bibitem[BNT07]{BertsimasNohadaniTeo2007}
D.~Bertsimas, O.~Nohadani, and K.~M. Teo.
\newblock Robust optimization in electromagnetic scattering problems.
\newblock {\em Journal of Applied Physics}, 101(7):074507, 2007.

\bibitem[BNT10a]{BertsimasNohadaniTeo2010nonconvex}
D.~Bertsimas, O.~Nohadani, and K.~M. Teo.
\newblock Nonconvex robust optimization for problems with constraints.
\newblock {\em INFORMS journal on computing}, 22(1):44--58, 2010.

\bibitem[BNT10b]{BertsimasNohadaniTeo2010}
D.~Bertsimas, O.~Nohadani, and K.~M. Teo.
\newblock Robust optimization for unconstrained simulation-based problems.
\newblock {\em Operations Research}, 58(1):161--178, 2010.

\bibitem[BS04]{BertsimasSim2004}
D.~Bertsimas and M.~Sim.
\newblock The price of robustness.
\newblock {\em Operations Research}, 52(1):35--53, 2004.

\bibitem[BS07]{BeyerSendhoff2007}
H-G. Beyer and B.~Sendhoff.
\newblock Robust optimization--a comprehensive survey.
\newblock {\em Computer methods in applied mechanics and engineering},
  196(33):3190--3218, 2007.

\bibitem[BTBB10]{BenTalBertsimasBrown2010}
A.~Ben-Tal, D.~Bertsimas, and D.~Brown.
\newblock A soft robust model for optimization under ambiguity.
\newblock {\em Operations research}, 58(4-part-2):1220--1234, 2010.

\bibitem[BTGN09]{BenTalElGhaouiNemirovski2009}
A.~Ben-Tal, L.~El Ghaoui, and A.~Nemirovski.
\newblock {\em Robust Optimization}.
\newblock Princeton University Press, Princeton and Oxford, 2009.

\bibitem[BTN98]{BenTalNemirovski1998}
A.~Ben-Tal and A.~Nemirovski.
\newblock Robust convex optimization.
\newblock {\em Mathematics of Operations Research}, 23(4):769--805, 1998.

\bibitem[CG16]{ChasseinGoerigk2016}
A.~Chassein and M.~Goerigk.
\newblock A bicriteria approach to robust optimization.
\newblock {\em Computers \& Operations Research}, 66:181 -- 189, 2016.

\bibitem[Cha93]{Chazelle1993}
B.~Chazelle.
\newblock An optimal convex hull algorithm in any fixed dimension.
\newblock {\em Discrete {\&} Computational Geometry}, 10(4):377--409, Dec 1993.

\bibitem[CSZ09]{CramerSudhoffZivi2009}
A.~M. Cramer, S.~D. Sudhoff, and E.~L. Zivi.
\newblock Evolutionary algorithms for minimax problems in robust design.
\newblock {\em IEEE Transactions on Evolutionary Computation}, 13(2):444--453,
  2009.

\bibitem[DHX17]{DiazHandlXu2017}
J.~Esteban Diaz, J.~Handl, and D-L. Xu.
\newblock Evolutionary robust optimization in production planning –
  interactions between number of objectives, sample size and choice of
  robustness measure.
\newblock {\em Computers \& Operations Research}, 79:266 -- 278, 2017.

\bibitem[dMB14]{HomemdeMelloBayraksan2014}
T.~Homem de~Mello and G.~Bayraksan.
\newblock Monte carlo sampling-based methods for stochastic optimization.
\newblock {\em Surveys in Operations Research and Management Science}, 19(1):56
  -- 85, 2014.

\bibitem[EKS17]{eichfelder2017decision}
G.~Eichfelder, C.~Kr{\"u}ger, and A.~Sch{\"o}bel.
\newblock Decision uncertainty in multiobjective optimization.
\newblock {\em Journal of Global Optimization}, 69(2):485--510, 2017.

\bibitem[ES12]{Eiben2012}
A.~E. Eiben and S.~K. Smit.
\newblock {\em Evolutionary Algorithm Parameters and Methods to Tune Them},
  pages 15--36.
\newblock Springer Berlin Heidelberg, Berlin, Heidelberg, 2012.

\bibitem[GLT97]{GoldenLaporteTaillard1997}
B.~L. Golden, G.~Laporte, and É.~D. Taillard.
\newblock An adaptive memory heuristic for a class of vehicle routing problems
  with minmax objective.
\newblock {\em Computers \& Operations Research}, 24(5):445 -- 452, 1997.

\bibitem[GS10]{GohSim2010}
J.~Goh and M.~Sim.
\newblock Distributionally robust optimization and its tractable
  approximations.
\newblock {\em Operations Research}, 58(4-part-1):902--917, 2010.

\bibitem[GS16]{GoerigkSchobel2016}
M.~Goerigk and A.~Sch\"{o}bel.
\newblock Algorithm engineering in robust optimization.
\newblock In L.~Kliemann and P.~Sanders, editors, {\em Algorithm Engineering:
  Selected Results and Surveys}, volume 9220 of LNCS State of the Art of {\em
  Lecture Notes in Computer Science}, pages 245--279. Springer Berlin /
  Heidelberg, 2016.

\bibitem[Her99]{Herrmann1999}
J.~W. Herrmann.
\newblock A genetic algorithm for minimax optimization problems.
\newblock In {\em Evolutionary Computation, 1999. CEC 99. Proceedings of the
  1999 Congress on}, volume~2, pages 1099--1103. IEEE, 1999.

\bibitem[Jen04]{Jensen2004}
M.~T. Jensen.
\newblock {\em A New Look at Solving Minimax Problems with Coevolutionary
  Genetic Algorithms}, pages 369--384.
\newblock Springer US, Boston, MA, 2004.

\bibitem[JSW98]{JonesSchonlauWelch1998}
D.~R. Jones, M.~Schonlau, and W.~J. Welch.
\newblock Efficient global optimization of expensive black-box functions.
\newblock {\em J. of Global Optimization}, 13(4):455--492, December 1998.

\bibitem[JY13]{JamilYang2013}
M.~Jamil and X-S. Yang.
\newblock A literature survey of benchmark functions for global optimization
  problems.
\newblock {\em International Journal of Mathematical Modelling and Numerical
  Optimisation (IJMMNO)}, 4(2):150--194, 2013.

\bibitem[Kru12]{Kruisselbrink2012}
J.~W. Kruisselbrink.
\newblock {\em Evolution strategies for robust optimization}.
\newblock PhD thesis, Leiden Institute of Advanced Computer Science (LIACS),
  Faculty of Science, Leiden university, 2012.

\bibitem[KY97]{KouvelisYu1997}
P.~Kouvelis and G.~Yu.
\newblock {\em Robust Discrete Optimization and Its Applications}.
\newblock Kluwer Academic Publishers, 1997.

\bibitem[MKA11]{MasudaKuriharaAiyoshi2011}
K.~Masuda, K.~Kurihara, and E.~Aiyoshi.
\newblock A novel method for solving min-max problems by using a modified
  particle swarm optimization.
\newblock In {\em Systems, Man, and Cybernetics (SMC), 2011 IEEE International
  Conference on}, pages 2113--2120. IEEE, 2011.

\bibitem[MLM15]{MirjaliliLewisMostaghim2015}
S.~Mirjalili, A.~Lewis, and S.~Mostaghim.
\newblock Confidence measure: A novel metric for robust meta-heuristic
  optimisation algorithms.
\newblock {\em Information Sciences}, 317:114 -- 142, 2015.

\bibitem[MWPL13]{MarzatWalterPietLahanier2013}
J.~Marzat, E.~Walter, and H.~Piet-Lahanier.
\newblock Worst-case global optimization of black-box functions through kriging
  and relaxation.
\newblock {\em Journal of Global Optimization}, 55(4):707--727, Apr 2013.

\bibitem[MWPL16]{MarzatWalterPietLahanier2016}
J.~Marzat, E.~Walter, and H.~Piet-Lahanier.
\newblock A new expected-improvement algorithm for continuous minimax
  optimization.
\newblock {\em Journal of Global Optimization}, 64(4):785--802, 2016.

\bibitem[Nah17]{Nahr2017}
C.~Nahr.
\newblock Tektosyne library for java, 2017.
\newblock Available at http://www.kynosarges.org.

\bibitem[OS97]{OkabeSuzuki1997}
A.~Okabe and A.~Suzuki.
\newblock Locational optimization problems solved through voronoi diagrams.
\newblock {\em European Journal of Operational Research}, 98(3):445 -- 456,
  1997.

\bibitem[PBJ06]{PaenkeBrankeJin2006}
I.~Paenke, J.~Branke, and Y.~Jin.
\newblock Efficient search for robust solutions by means of evolutionary
  algorithms and fitness approximation.
\newblock {\em IEEE Transactions on Evolutionary Computation}, 10(4):405--420,
  2006.

\bibitem[SE98]{ShiEberhart1998}
Y.~Shi and R.~Eberhart.
\newblock A modified particle swarm optimizer.
\newblock In {\em 1998 IEEE International Conference on Evolutionary
  Computation Proceedings. IEEE World Congress on Computational Intelligence
  (Cat. No.98TH8360)}, pages 69--73, 1998.

\bibitem[SK02]{ShiKrohling2002}
Y.~Shi and R.~Krohling.
\newblock Co-evolutionary particle swarm optimization to solve min-max
  problems.
\newblock In {\em Evolutionary Computation, 2002. CEC'02. Proceedings of the
  2002 Congress on}, volume~2, pages 1682--1687. IEEE, 2002.

\bibitem[Tal09]{Ghazali2009}
E-G. Talbi.
\newblock {\em Metaheuristics: From Design to Implementation}.
\newblock Wiley Publishing, 2009.

\bibitem[Tou83]{Toussaint1983}
G.~T. Toussaint.
\newblock Computing largest empty circles with location constraints.
\newblock {\em International Journal of Computer {\&} Information Sciences},
  12(5):347--358, 1983.

\bibitem[uRL17]{urRehmanLangelaar2017}
S.~ur~Rehman and M.~Langelaar.
\newblock Expected improvement based infill sampling for global robust
  optimization of constrained problems.
\newblock {\em Optimization and Engineering}, 18(3):723--753, Sep 2017.

\bibitem[uRLvK14]{urRehmanLangelaarvanKeulen2014}
S.~ur~Rehman, M.~Langelaar, and F.~van Keulen.
\newblock Efficient kriging-based robust optimization of unconstrained
  problems.
\newblock {\em Journal of Computational Science}, 5(6):872 -- 881, 2014.

\bibitem[VDYL16]{VuDAmbrosioHamadiLiberti2016}
K.~Khac Vu, C.~D'Ambrosio, Y.Hamadi, and L.~Liberti.
\newblock Surrogate‐based methods for black‐box optimization.
\newblock {\em International Transactions in Operational Research},
  24(3):393--424, 2016.

\bibitem[VMdCM11]{ValleMartinezdaCunhaMateus2011}
C.~Arbex Valle, L.~Conegundes Martinez, A.~Salles da~Cunha, and G.~R. Mateus.
\newblock Heuristic and exact algorithms for a min–max selective vehicle
  routing problem.
\newblock {\em Computers \& Operations Research}, 38(7):1054 -- 1065, 2011.

\bibitem[ZZ10]{ZhouZhang2010}
A.~Zhou and Q.~Zhang.
\newblock A surrogate-assisted evolutionary algorithm for minimax optimization.
\newblock In {\em IEEE Congress on Evolutionary Computation}, pages 1--7, July
  2010.

\end{thebibliography}
\end{document}